\documentclass[11pt,english]{article}
\usepackage{lmodern}
\usepackage[T1]{fontenc}
\usepackage[latin9]{inputenc}
\usepackage{geometry}
\usepackage{amsmath}
\usepackage{amssymb}
\usepackage{stmaryrd}
\usepackage{graphicx}
\usepackage{wasysym}

\makeatletter
\date{}

\usepackage{amsthm}
\usepackage{mathrsfs}
\usepackage{enumerate}
\usepackage{bbm}

\usepackage[title]{appendix}

\usepackage{url}

\renewcommand{\subparagraph}{%
 \@startsection {subparagraph}{5}{\z@ }{3.25ex \@plus 1ex
 \@minus .2ex}{-1em}{\normalfont \normalsize \bfseries }}%

\usepackage{graphicx}
\usepackage[shortlabels]{enumitem}
\usepackage{tikz-cd}
\usetikzlibrary{arrows}

\theoremstyle{plain}

\theoremstyle{definition}
\newtheorem{main def}[thm]{Main Definition}\theoremstyle{remark}

\usepackage{fancyhdr}

\makeatother

\usepackage{babel}
\begin{document}
\title{Classes of ODE solutions: smoothness, covering numbers, implications
for noisy function fitting, and the curse of smoothness phenomenon\thanks{
numbers; local complexity; smooth functions; kernel methods}}
\author{Ying Zhu\thanks{Assistant Professor, Department of Economics, UC San Diego. Corresponding
author. yiz012@ucsd.edu.} ~~~and~~~Mozhgan Mirzaei\thanks{PhD, former student at Department of Mathematics, UC San Diego. momirzae@ucsd.edu}\\
\\
First draft on ArXiv: November 2020. This draft: March 2021.}
\maketitle
\begin{abstract}
Many numerical methods for recovering ODE solutions from data rely
on approximating the solutions using basis functions or kernel functions
under a least square criterion. The accuracy of this approach hinges
on the smoothness of the solutions. This paper provides a theoretical
foundation for these methods by establishing novel results on the
smoothness and covering numbers of ODE solution classes (as a measure
of their ``size''). Our results provide answers to ``how do the
degree of smoothness and the ``size'' of a class of ODEs affect
the ``size'' of the associated class of solutions?'' In particular,
we show that: (1) for $y^{'}=f\left(y\right)$ and $y^{'}=f\left(x,\,y\right)$,
if the absolute values of all $k$th ($k\leq\beta+1$) order derivatives
of $f$ are bounded by $1$, then the solution can end up with the
$(k+1)$th derivative whose magnitude grows factorially fast in $k$
-- ``a curse of smoothness''; (2) our upper bounds for the covering
numbers of the $(\beta+2)-$degree smooth solution classes are greater
than those of the ``standard'' $(\beta+2)-$degree smooth class
of univariate functions; (3) the mean squared error of least squares
fitting for noisy recovery has a convergence rate no larger than $\left(\frac{1}{n}\right)^{\frac{2\left(\beta+2\right)}{2\left(\beta+2\right)+1}}$
if $n=\Omega\left(\left(\beta\sqrt{\log\left(\beta\vee1\right)}\right)^{4\beta+10}\right)$,
and under this condition, the rate $\left(\frac{1}{n}\right)^{\frac{2\left(\beta+2\right)}{2\left(\beta+2\right)+1}}$
is minimax optimal in the case of $y^{'}=f\left(x,\,y\right)$; (4)
more generally, for the higher order Picard type ODEs, $y^{\left(m\right)}=f\left(x,\,y,\,y^{'},\,...,y^{\left(m-1\right)}\right)$,
the covering number of the solution class is bounded from above by
the product of the covering number of the class $\mathcal{F}$ that
$f$ ranges over and the covering number of the set where initial
values lie.
\end{abstract}

\section{Introduction}

Ordinary differential equations (ODE) enjoy a long standing history
in mathematics and have numerous applications in science and engineering.
ODE also play a crucial role in social science and business, for example,
the famous Solow growth model in economics, and the famous Bass product
diffusion model in marketing. Since the COVID-19 pandemic, lots of
attention has been given to the compartmental models in epidemiology
(see, e.g., \cite{Vynnycky and white 2010}). The compartmental models
have also been used in forecasting election outcomes (\cite{Volkenning et al. 2020})
as well as modeling the information diffusion (see, e.g., \cite{Cinelli et al. 2020}).

As analytical solutions for ODEs are often unavailable, developing
numerical methods for solving initial value or boundary value problems
has been an active research area in ODEs. An important technique (belonging
in the family of collocation methods) tries to approximate the solutions
using basis functions or kernel functions (see \cite{Deuflhard and Bornemann 2000,Meade and Fernadez 1994,Mehrkanoon et al 2012}).
For example, Mehrkanoon, et al. \cite{Mehrkanoon et al 2012} consider
approximate solutions in the form $\hat{y}\left(x\right)=w^{T}\varphi\left(x\right)+b$
where $\varphi:\mathbb{R}\mapsto\mathbb{R}^{p}$ is the feature map
and $w\in\mathbb{R}^{p}$ are the unknown weights (solved from a least
squares support vector machines formulation at collocation points). 

Another active area in ODEs involves the recovery of solutions from
their noisy measurements:

\begin{equation}
Y_{i}=y(x_{i})+\varepsilon_{i},\quad i=1,...,n.\label{eq:model}
\end{equation}
In the above, $y$ is a solution to some ODE, $\left\{ Y_{i}\right\} _{i=1}^{n}$
are noisy measurements of $\left\{ y\left(x_{i}\right)\right\} _{i=1}^{n}$
evaluated at a collection of fixed design points $\left\{ x_{i}\right\} _{i=1}^{n}$
(both $Y_{i}$ and $x_{i}$ are observed), and $\left\{ \varepsilon_{i}\right\} _{i=1}^{n}$
are unobserved noise terms. Noisy recovery of solutions can be useful
in situations where prediction is the ultimate goal; it is also used
as an intermediate step or part of an algorithm to recover parameters
in the ODEs. In the literature, function fitting with noisy measurements
of ODE solutions is typically performed via the least squares procedure:
\begin{equation}
\hat{y}\in\arg\min_{\tilde{y}\in\Pi}\frac{1}{2n}\sum_{i=1}^{n}\left(Y_{i}-\tilde{y}\left(x_{i}\right)\right)^{2}\label{eq:least squares}
\end{equation}
where $\Pi$ is a suitably chosen function class that contains $y$.
For example, Volkening et al. \cite{Volkenning et al. 2020} minimize
the sum of squared deviations between the averaged polling data and
the fitted solutions to a compartmental model. In this application,
analytical solutions to the ODEs are available so $\Pi$ would be
the class of functions with the same form as the ODE solution, indexed
by finitely dimensional parameters. In many problems, analytical solutions
for ODEs are unavailable. While integration based nonlinear least
squares overcome this issue, they often suffer many computational
issues (see discussions in \cite{Varah 1982}); for more discussions
on computational challenges in the standard nonlinear least squares
methods, see \cite{Liang and Wu 2008,Ramsay et al. 2007}. Therefore,
researchers often resort to approximate the underlying solutions using
polynomials and spline bases (see \cite{Liang and Wu 2008,Poyton et al. 2006,ramsay and silverman 2005,Ramsay 1996,Ramsay et al. 2007,Varah 1982}).
For example, in clinical studies of AIDs, plasma viral load and CD4+
T cell counts are measured with additive noise. With these measurements,
Liang and Wu \cite{Liang and Wu 2008} apply local polynomial regressions
with the least squares criterion to recover the ODE solutions and
their first derivatives. In the framework proposed by \cite{Liang and Wu 2008}
(see their Section 2.1), $\Pi$ would correspond to the class of three
times differentiable functions with bounded derivatives. 

In either of the two areas described above, many existing methods
rely on approximating the underlying solutions using basis functions
or kernel functions. This strategy assumes some ambient space of smooth
functions that contains the solutions and seeks a close enough estimator
from the restricted class of functions. The choice of the function
class and the accuracy of the approximations hinge on the smoothness
of $y$. For example, in terms of the local polynomial estimators
used in \cite{Liang and Wu 2008} for recovering $y$ and $y^{'}$
in the equation $y^{'}\left(x\right)=f\left(y\left(x\right)\right)$,
for such estimators to have well behaved biases, boundedness on the
second and third derivatives of $y$ is needed, respectively. As suggested
by \cite{Liang and Wu 2008}, higher degree local polynomials can
be applied and if so, boundedness on higher order derivatives of $y$
would be needed. Because the smoothness of $y$ ultimately depends
on the smoothness of $f$, this prompts us to ask the following fundamental
question: 
\begin{itemize}
\item \textbf{\textit{(Question 1)}}\textit{ in general, how does the degree
of smoothness of $f$ affect the smoothness of $y$ in $y^{'}\left(x\right)=f\left(y\left(x\right)\right)$
and $y^{'}\left(x\right)=f\left(x,\,y\left(x\right)\right)$; for
example, if $f$ in $y^{'}\left(x\right)=f\left(y\left(x\right)\right)$
is $\beta+1$ times differentiable with all derivatives bounded by
$1$ (where $\beta$ is a non-negative integer), can one simply reason
that all derivatives of $y$ also stay ``nicely'' bounded? }
\end{itemize}
To our knowledge, it appears that no rigorous answers have been provided
to the question above. Our results regarding \textbf{\textit{Question
1}} can be summarized as follows: 
\begin{itemize}
\item \textbf{\textit{(Lemma 3.1)}} \textit{For $y^{'}\left(x\right)=f\left(y\left(x\right)\right)$,
for all $k=0,...,\beta+1$, if $\left|f^{(k)}\left(x\right)\right|\leq1$,
then $\left|y^{\left(k+1\right)}\left(x\right)\right|\leq k!$. These
factorial bounds can be attained (and hence tight). }
\item \textbf{\textit{(Lemma 3.2)}} \textit{For $y^{'}\left(x\right)=f\left(x,\,y\left(x\right)\right)$,
if the absolute values of all $k$th order partial derivatives are
bounded by $1$, then $\left|y^{\left(k+1\right)}\left(x\right)\right|\leq2^{k}k!$
for $k=0,...,\beta+1$.}\footnote{The bound ``$1$'' is assumed to avoid notation cluttering and can
be relaxed easily. }\textit{ }
\end{itemize}
Noisy recovery is generally more difficult than noiseless recovery.
In terms of (\ref{eq:least squares}), the worst case bound for the
average squared error of $\hat{y}$ depends on the maximum ``correlation''
between the noise vector $\left\{ \varepsilon_{i}\right\} _{i=1}^{n}$
and the random error vector $\left\{ \hat{y}(x_{i})-y(x_{i})\right\} _{i=1}^{n}$.
To see this, note that since $y\in\Pi$ ($y$ is feasible) and $\hat{y}$
is optimal, we have 
\[
\frac{1}{2n}\sum_{i=1}^{n}\left(Y_{i}-\hat{y}\left(x_{i}\right)\right)^{2}\leq\frac{1}{2n}\sum_{i=1}^{n}\left(Y_{i}-y\left(x_{i}\right)\right)^{2},
\]
which yields 
\begin{equation}
\frac{1}{n}\sum_{i=1}^{n}\left(\hat{y}(x_{i})-y(x_{i})\right)^{2}\leq\frac{2}{n}\sum_{i=1}^{n}\varepsilon_{i}\left(\hat{y}(x_{i})-y(x_{i})\right).\label{eq:44-1-1}
\end{equation}
Controlling this right-hand-side (RHS) term in (\ref{eq:44-1-1})
can be reduced to controlling the ``local complexity'' (one of the
most important notions in machine learning theory) associated with
$\Pi$, which depends ultimately on the ``size'' of $\Pi$. In noisy
recovery problems, researchers typically have prior knowledge on the
smoothness structure of $f$ and use such information to guide the
choice of $\Pi$ in (\ref{eq:least squares}); therefore, a natural
choice for $\Pi$ is 
\begin{equation}
\mathcal{Y}=\left\{ y:\,y^{\left(m\right)}\left(x\right)=f\left(x,\,y\left(x\right),\,y^{'}\left(x\right),\,...,y^{\left(m-1\right)}\left(x\right)\right),\,f\in\mathcal{F}\right\} \label{eq:higher order}
\end{equation}
with $y=y^{(0)}$. In many contexts, $\mathcal{Y}$ is a set consisting
of infinitely many smooth functions. The notion ``covering number'',
dated back to the seminal work of Kolmogorov, Tikhomirov, and others
(see, e.g., \cite{Kolmogorov and Tikhomirov 1959}), provides a way
to measure the size of a class with an infinite number of elements.
In this paper, we derive \textit{nonasymptotic} bounds on the covering
numbers of $\mathcal{Y}$ to address the following novel question:
\begin{itemize}
\item \textbf{\textit{(Question 2)}}\textit{ how does the ``size'' of
a class of $f$ affect the ``size'' of the associated solution class
of $y$? }
\end{itemize}
Like \textbf{\textit{Question 1}}, we are not aware of any rigorous
results regarding\textbf{\textit{ Question 2}} from the existing literature.
Below provides an overview of our results regarding \textbf{\textit{Question
2}}.
\begin{itemize}
\item \textbf{\textit{(Theorem 2.1)}} \textit{We establish a general upper
bound on the covering number of solution classes associated with (\ref{eq:higher order}).
This result implies, the covering number of the underlying solution
class $\mathcal{Y}$ is bounded from above by the product of the covering
number of the class $\mathcal{F}$ that $f$ ranges over and the covering
number of the set where initial values lie, as long as $f$ satisfies
a Lipschitz condition with respect to the $y,...,y^{\left(m-1\right)}$
coordinates in $l_{2}-$norm (i.e., the Picard Lipschitz condition).
This general bound yields a sharp scaling in some problems (}\textbf{\textit{Corollary
2.1}}\textit{). }
\item \textbf{\textit{(Theorems 3.1-3.2)}} \textit{For the first order ODEs,
if the absolute values of all $k$th }($k\leq\beta+1$)\textit{ order
(partial) derivatives of $f$ are bounded by $1$, then the general
bound in }\textbf{\textit{Theorem 2.1}}\textit{ may be improved by
exploiting the factorial bounds in }\textbf{\textit{Lemmas 3.1 and
3.2}}\textit{.}\textbf{\textit{ }}\textit{Our upper bounds for the
covering numbers of the $(\beta+2)-$degree smooth solution classes,
$\mathcal{Y}$, are greater than those of the ``standard'' $(\beta+2)-$degree
smooth class, $\mathcal{S}_{\beta+2}$, of univariate functions (where
all derivatives are bounded by $1$). We also discuss a lower bound
on the covering numbers of solution classes associated with $y^{'}\left(x\right)=f\left(x,\,y\left(x\right)\right)$
in }\textbf{\textit{Lemma 3.3}}\textit{.}
\end{itemize}
The abovementioned bounds have important implications for analyzing
techniques for solving ODEs as well as noisy recovery. In this paper,
we focus on the noisy recovery problem and examine the implications
on the convergence rate of the mean squared error $\mathbb{E}\left[\frac{1}{n}\sum_{i=1}^{n}\left(\hat{y}(x_{i})-y(x_{i})\right)^{2}\right]$.
Let us first take a digression from discussing these implications.
For noisy recovery of a function in the standard $(\gamma+1)-$degree
smooth class of $m-$variate functions, the common wisdom is that
large $m$ brings a curse (resulting in a large convergence rate of
the least squares) and large $\gamma$ brings a bless (resulting in
a small convergence rate). While there is no ambiguity from the ``curse
of dimension'' here, this paper discovers that the ``bless of smoothness''
comes at the price of assuming a condition on the sample size $n$.
In particular, for noisy recovery of a function in the standard $(\beta+2)-$degree
smooth class, $\mathcal{S}_{\beta+2}$, the well-known minimax optimal
rate $\left(\frac{1}{n}\right)^{\frac{2\left(\beta+2\right)}{2\left(\beta+2\right)+1}}$
(which decreases in $\beta$) essentially coincides with the convergence
rate of the least squares when $n=\Omega\left(\left(\beta\vee1\right)^{2\beta+5}\right)$.
We view this requirement on the sample size, $n=\Omega\left(\left(\beta\vee1\right)^{2\beta+5}\right)$,
a ``curse of smoothness'' for least squares estimations of elements
in $\mathcal{S}_{\beta+2}$.\footnote{We also discover that the ``curse of smoothness'' is exacerbated
as the dimension $m$ increases. This observation can have important
practical implications. } A more detailed discussion is given in Section 3.1. 

Coming back to the ODE problems, with the assistance of \textbf{\textit{Theorems
3.1-3.2}} and \textbf{\textit{Lemma 3.3}}, we ask the following novel
questions regarding $y^{'}\left(x\right)=f\left(y\left(x\right)\right)$
and $y^{'}\left(x\right)=f\left(x,\,y\left(x\right)\right)$: 
\begin{itemize}
\item \textbf{\textit{(Question 3)}}\textit{ If $f\in\mathcal{S}_{\beta+1}$,
what is the condition on $n$ for the least squares $\hat{y}$ to
have $\mathbb{E}\left[\frac{1}{n}\sum_{i=1}^{n}\left(\hat{y}(x_{i})-y(x_{i})\right)^{2}\right]$
bounded above by $\left(\frac{1}{n}\right)^{\frac{2\left(\beta+2\right)}{2\left(\beta+2\right)+1}}$,
and what can we say about the optimality of $\left(\frac{1}{n}\right)^{\frac{2\left(\beta+2\right)}{2\left(\beta+2\right)+1}}$?}
\end{itemize}
Here are the answers to \textbf{\textit{Question 3}}.
\begin{itemize}
\item \textbf{\textit{(Theorem 3.3)}} \textit{The mean squared error of
$\hat{y}$ has a convergence rate essentially no larger than $\left(\frac{1}{n}\right)^{\frac{2\left(\beta+2\right)}{2\left(\beta+2\right)+1}}$
if $n=\Omega\left(\left(\beta\sqrt{\log\left(\beta\vee1\right)}\right)^{4\beta+10}\right)$,
and under this condition, the rate $\left(\frac{1}{n}\right)^{\frac{2\left(\beta+2\right)}{2\left(\beta+2\right)+1}}$
is minimax optimal in the case of $y^{'}\left(x\right)=f\left(x,\,y\left(x\right)\right)$.
Relative to the standard smooth class $\mathcal{S}_{\beta+2}$ (where
the least squares procedure needs $n=\Omega\left(\left(\beta\vee1\right)^{2\beta+5}\right)$
for the optimal rate $\left(\frac{1}{n}\right)^{\frac{2\left(\beta+2\right)}{2\left(\beta+2\right)+1}}$),
the requirement on the sample size in the ODE context is much greater.}\textbf{\textit{
Theorem 3.3}}\textit{ also presents upper bounds (with rates no smaller
than $\left(\frac{1}{n}\right)^{\frac{2\left(\beta+2\right)}{2\left(\beta+2\right)+1}}$)
on the mean squared error of $\hat{y}$ under more general situations
where $n=\Omega\left(\left(\beta\sqrt{\log\left(\beta\vee1\right)}\right)^{4\beta+10}\right)$
is not imposed. }
\end{itemize}
Having provided the theoretical implications of (\ref{eq:least squares})
in \textbf{\textit{Theorem 3.3}}, we further ask whether there is
a practical implementation of (\ref{eq:least squares}) with comparable
theoretical guarantees. We exploit the smoothness structures of solutions
in \textbf{\textit{Lemmas 3.1--3.2}} directly by considering some
ambient reproducing kernel Hilbert spaces that contain $\mathcal{Y}$,
and doing so allows us to develop a variant of the classical kernel
ridge regression procedure. For all practical purposes, the performance
guarantee for this procedure (provided in \textbf{\textit{Theorem
3.4}}) turns out quite comparable to the general bounds in \textbf{\textit{Theorem
3.3}}. The practical implementation associated with \textbf{\textit{Theorem
3.4}} provides an example of how our theoretical results can be used
to guide the use of kernel functions in noisy function fitting. 

Besides recovering the ODE solutions, one could also recover the first
derivatives of ODE solutions with a least squares approach. This paper
derives bounds for the covering number of the class consisting of
the first derivative $y^{'}$ of $y\in\mathcal{Y}$ as well. Therefore,
our arguments for Theorems 3.3 and 3.4 can be easily extended for
analyzing the recovery of the first derivatives of ODE solutions. 

\section{A general upper bound on the covering number}

\textbf{General notation}. The $l_{q}-$norm of a $K-$dimensional
vector $\theta$ is denoted by $\left|\theta\right|{}_{q}$, $1\leq q\leq\infty$
where $\left|\theta\right|{}_{q}:=\left(\sum_{j=1}^{K}|\theta_{j}|^{q}\right)^{1/q}$
when $1\leq q<\infty$ and $\left|\theta\right|{}_{q}:=\max_{j=1,...,K}|\theta_{j}|$
when $q=\infty$. Define $\mathbb{P}_{n}:=\frac{1}{n}\sum_{i=1}^{n}\delta_{x_{i}}$
that places a weight $\frac{1}{n}$ on each observation $x_{i}$ for
$i=1,...,n$, and the associated $\mathcal{L}^{2}(\mathbb{P}_{n})-$norm
of the vector $v:=\left\{ v(x_{i})\right\} _{i=1}^{n}$, denoted by
$\left|v\right|_{n}$, is given by $\left[\frac{1}{n}\sum_{i=1}^{n}\left(v(x_{i})\right)^{2}\right]^{\frac{1}{2}}$.
For two functions $f$ and $g$ on $\left[a,\,b\right]^{d}\subseteq\mathbb{R}^{d}$,
we denote the supremum metric by $\left|f-g\right|_{\infty}:=\sup_{x\in\left[a,\,b\right]^{d}}\left|f\left(x\right)-g\left(x\right)\right|$.
For functions $f(n)$ and $g(n)$, $f(n)\succsim g(n)$ means that
$f(n)\geq cg(n)$ for a universal constant $c\in(0,\,\infty)$; similarly,
$f(n)\precsim g(n)$ means that $f(n)\leq c^{'}g(n)$ for a universal
constant $c^{'}\in(0,\,\infty)$; and $f(n)\asymp g(n)$ means that
$f(n)\succsim g(n)$ and $f(n)\precsim g(n)$ hold simultaneously.
As a general rule for this paper, the various $c$ and $C$ constants
(all $\precsim1$) denote positive universal constants that are independent
of the sample size $n$ and the smoothness parameter $\beta$, and
may vary from place to place. \\
\\
\textbf{Definition} (covering numbers). Given a set $\mathbb{T}$,
a set $\left\{ \mathrm{t}^{1},\,\mathrm{t}^{2},...,\mathrm{t}^{N}\right\} \subset\mathbb{T}$
is called a $\delta-$cover of $\mathbb{T}$ with respect to a metric
$\rho$ if for each $\mathrm{t}\in\mathbb{T}$, there exists some
$i\in\left\{ 1,...,N\right\} $ such that $\rho(\mathrm{t},\,\mathrm{t}^{i})\leq\delta$.
The cardinality of the smallest $\delta-$cover is denoted by $N_{\rho}(\delta;\,\mathbb{T})$,
namely, the $\delta-$covering number of $\mathbb{T}$. For example,
$N_{\infty}\left(\delta,\,\mathcal{\mathcal{F}}\right)$ denotes the
$\delta-$covering number of a function class $\mathcal{F}$ with
respect to the supremum metric $\left|\cdot\right|_{\infty}$. \\
\\
Let us define
\begin{eqnarray*}
Y\left(x\right):=\left[\begin{array}{c}
y\left(x\right)\\
y^{'}\left(x\right)\\
\vdots\\
y^{\left(m-1\right)}\left(x\right)
\end{array}\right]^{T} & \textrm{and} & Y_{0}:=\left[\begin{array}{c}
y_{\left(0\right)}\\
y_{\left(1\right)}\\
\vdots\\
y_{\left(m-1\right)}
\end{array}\right]^{T}.
\end{eqnarray*}
We consider the ODE 
\begin{equation}
y^{\left(m\right)}\left(x\right)=f\left(x,\,Y\left(x\right)\right)\label{eq:22-1-1}
\end{equation}
with the initial values 
\[
y_{(0)}=y\left(a_{0}\right),\,y_{(1)}=y^{'}\left(a_{0}\right),\,...,\,y_{(m-1)}=y^{\left(m-1\right)}\left(a_{0}\right)
\]
such that $\left|Y_{0}\right|_{2}\leq C_{0}$ and 
\[
\left(x,\,Y\left(x\right)\right)\in\left\{ \left(x,\,Y\right):\,x\in\left[a_{0},\,a_{0}+a\right],\,\left|Y\right|_{2}\leq b+C_{0}\right\} :=\bar{\Gamma}
\]
(where $a,\,b>0$). Assume $f\left(x,\,Y\right)$ is continuous on
$\bar{\Gamma}$ and satisfies the Picard Lipschitz condition\textit{
\begin{equation}
\left|f\left(x,\,Y\right)-f\left(x,\,\tilde{Y}\right)\right|\leq L\left|Y-\tilde{Y}\right|_{2}\label{eq:lip}
\end{equation}
}for all $\left(x,\,Y\right):=\left(x,\,y,\,...,\,y^{\left(m-1\right)}\right)$
and $\left(x,\,\tilde{Y}\right):=\left(x,\,\tilde{y},\,...,\,\tilde{y}^{\left(m-1\right)}\right)$
in $\bar{\Gamma}$. Then by the Picard local existence theorem (see,
e.g., \cite{Coddington and Levinson 1955}), there is a solution to
(\ref{eq:22-1-1}) in $\mathcal{C}\left[a_{0},\,a_{0}+\alpha\right]$
where $\alpha\leq\min\left\{ a,\,\frac{b}{M}\right\} $ and $M=\max_{\left(x,\,Y\right)\in\bar{\Gamma}}\left|f\left(x,\,Y\right)\right|$. 

The Picard Lipschitz condition (\ref{eq:lip}) is not necessary for
the existence of a solution as suggested by the Cauchy-Peano existence
theorem (see, e.g., \cite{Coddington and Levinson 1955}). However,
(\ref{eq:lip}) allows us to establish a general upper bound on the
covering number of higher order ODEs. In what follows, we consider
the ODE (\ref{eq:22-1-1}) with $\left|Y_{0}\right|_{2}\leq C_{0}$
and $\left(x,\,Y\left(x\right)\right)\in\bar{\Gamma}$.\\
\textbf{}\\
\textbf{Theorem 2.1.} \textit{In terms of (\ref{eq:22-1-1}), suppose
$f$ is continuous on $\bar{\Gamma}$, and satisfies (\ref{eq:lip})
for all $\left(x,\,Y\right)$, $\left(x,\,\tilde{Y}\right)\in\bar{\Gamma}$.
If $f$ ranges over a function class $\mathcal{F}$ with $N_{\infty}\left(\delta,\,\mathcal{F}\right)$
on $\bar{\Gamma}$, and $\max_{\left(x,\,Y\right)\in\bar{\Gamma}}\left|f\left(x,\,Y\right)\right|\leq M$
for all $f\in\mathcal{F}$, then we have 
\begin{equation}
\log N_{\infty}\left(\delta,\,\mathcal{Y}_{k}\right)\leq\log N_{\infty}\left(\frac{\delta}{L_{\max}},\,\mathcal{F}\right)+m\log\left(\frac{2C_{0}L_{\max}}{\delta}+1\right)\label{eq:general_bound}
\end{equation}
where $L_{\max}=\sup_{x\in\left[a_{0},\,a_{0}+\alpha\right]}\left\{ \exp\left(x\sqrt{L^{2}+1}\right)\left[1+\int_{0}^{x}\exp\left(-s\sqrt{L^{2}+1}\right)ds\right]\right\} $
with $\alpha=\min\left\{ a,\,\frac{b}{M}\right\} $, and $\mathcal{Y}_{0}=\mathcal{Y}$
($k=0$) is the class consisting of solutions (to (\ref{eq:22-1-1})
with $f$ ranging over $\mathcal{F}$) on $\left[a_{0},\,a_{0}+\alpha\right]$
and $\mathcal{Y}_{k}$ (the non-negative integer $k\leq m-1$) is
the class consisting of the $k$th derivative $y^{(k)}$ of $y\in\mathcal{Y}_{0}$.
}\\

The proof for Theorem 2.1 is provided in Section \ref{subsec:Theorem-2.1}
of the supplementary materials. In the case where $m=1$, we can use
the sharper constant $L_{\max}=\sup_{x\in\left[0,\,\alpha\right]}\left\{ \exp\left(Lx\right)\left[1+\int_{0}^{x}\exp\left(-Ls\right)ds\right]\right\} $. 

If $m\log\left(\frac{2C_{0}L_{\max}}{\delta}+1\right)\precsim\log N_{\infty}\left(\frac{\delta}{L_{\max}},\,\mathcal{F}\right)$,
Theorem 2.1 implies that the solution class $\mathcal{Y}$ is ``essentially''
no larger than the class $f$ ranges over. Below we discuss the implications
of Theorem 2.1 for ODEs parameterized by a finite $K-$dimensional
vector of coefficients. 

\subsection{Implications on parametric ODEs}

\textbf{Corollary 2.1}. \textit{Consider the ODE 
\begin{equation}
y^{\left(m\right)}\left(x\right)=f\left(x,\,Y\left(x\right);\,\theta\right)\textrm{ with }\left(x,\,Y\left(x\right)\right)\ensuremath{\in}\bar{\Gamma}\textrm{ and initial values }Y_{0},\label{eq:11-2-2}
\end{equation}
where $f$ is parameterized by a finite $K-$dimensional vector of
coefficients 
\[
\theta\in\mathbb{B}_{q}\left(1\right):=\left\{ \theta\in\mathbb{R}^{K}:\,\left|\theta\right|_{q}\leq1\right\} 
\]
with $q\in\left[1,\,\infty\right]$. Suppose $f$ is continuous on
$\bar{\Gamma}$, $\max_{\left(x,\,Y\right)\in\bar{\Gamma}}\left|f\left(x,\,Y;\,\theta\right)\right|\leq M$,
and 
\begin{equation}
\left|f\left(x,\,Y;\,\theta\right)-f\left(x,\,\tilde{Y};\,\theta\right)\right|\leq L\left|Y-\tilde{Y}\right|_{2}\label{eq:15}
\end{equation}
for all $\left(x,\,Y\right),\,\left(x,\,\tilde{Y}\right)\in\bar{\Gamma}$
and $\theta\in\mathbb{B}_{q}\left(1\right)$; moreover, 
\begin{equation}
\left|f\left(x,\,Y;\,\theta\right)-f\left(x,\,Y;\,\theta^{'}\right)\right|\leq L_{K}\left|\theta-\theta^{'}\right|_{q},\label{eq:17}
\end{equation}
for all $\left(x,\,Y\right)\in\bar{\Gamma}$ and $\theta,\,\theta^{'}\in\mathbb{B}_{q}\left(1\right)$.
Let $\mathcal{Y}$ be the class consisting of solutions to (\ref{eq:11-2-2})
with $\theta\in\mathbb{B}_{q}\left(1\right)$ on $\left[a_{0},\,a_{0}+\min\left\{ a,\,\frac{b}{M}\right\} \right]$.
Then we have
\begin{equation}
\log N_{\infty}\left(\delta,\,\mathcal{Y}\right)\leq K\log\left(1+\frac{2L_{\max}L_{K}}{\delta}\right)+m\log\left(\frac{2C_{0}L_{\max}}{\delta}+1\right).\label{eq:16}
\end{equation}
}

The proof for Corollary 2.1 is provided in Section \ref{subsec:Corollary-2.1}
of the supplementary materials. 

In (\ref{eq:16}), the part ``$K\log\left(1+\frac{2L_{\max}L_{K}}{\delta}\right)$''
is related to the ``size'' of the ball that $\theta$ ranges over,
and the part ``$m\log\left(\frac{C_{0}L_{\max}}{\delta}+1\right)$''
is related to the ``size'' of the ball that the initial value $Y_{0}$
ranges over. The part ``$m\log\left(\frac{2C_{0}L_{\max}}{\delta}+1\right)$''
reveals an interesting feature of differential equations: even if
the equation is fixed and known (for example, $y^{'}\left(x\right)=y\left(x\right)$
with solutions $y\left(x\right)=c^{*}e^{x}$), the solution class
is still not a singleton (and has infinitely many solutions) unless
a fixed initial value is given. As a consequence, in the simple example
$y^{'}=y$, the noisy recovery of a solution to this ODE still requires
the estimation of $c^{*}$.\footnote{In terms of the homogeneous higher order linear ODEs

\begin{equation}
a_{0}\left(x\right)y+a_{1}\left(x\right)y^{(1)}+a_{2}\left(x\right)y^{(2)}+\cdots+a_{m}\left(x\right)y^{(m)}=0,\label{eq:linear}
\end{equation}
it is well understood that the solutions of (\ref{eq:linear}) (where
$a_{j}\left(\cdot\right)$s are fixed continuous functions) form a
vector space of dimension $m$. A classical result on the VC dimension
of a vector space of functions (\cite{Dudley 1978,Steele 1978}) implies
that the solution class associated with (\ref{eq:linear}) has VC
dimension at most $m$, when $a_{j}\left(\cdot\right)$s ($j=0,...,m$)
are fixed continuous functions.}

Bound (\ref{eq:16}) suggests that if the class $f$ ranges over is
a ``parametric'' class, then the associated solution class $\mathcal{Y}$
also behaves like a ``parametric'' class. Suppose the class of ODEs
(\ref{eq:11-2-2}) has a fixed initial value. Then the term ``$m\log\left(\frac{C_{0}L_{\max}}{\delta}+1\right)$''
can be dropped while the scaling of ``$K\log\left(1+\frac{2L_{\max}L_{K}}{\delta}\right)$''
can be attained as the following example suggests. Consider the simple
ODE $y^{'}\left(x\right)=-\theta y\left(x\right)$ with $\theta\in\left[0,\,1\right]$,
$x\in\left[0,\,1\right]$ and initial value $y\left(0\right)=1$,
which has the solution $y\left(x\right)=e^{-\theta x}$. It can be
easily verified that $\log N_{\infty}\left(\delta,\,\mathcal{Y}\right)\asymp\log\left(\frac{c}{\delta}+c^{'}\right)$
for some positive universal constants $c$ and $c^{'}$. 

While (\ref{eq:general_bound}) provides a sharp scaling in the example
above, this general bound may be improved in other contexts as we
will see in Section 3.2. 

\subsection{Noisy fitting with plain vanilla methods, issues, and remedies}

Let us consider (\ref{eq:model}) with $y$ a solution to the following
ODE
\begin{equation}
y^{'}\left(x;\,\theta^{*}\right)=f\left(x,\,y\left(x;\,\theta^{*}\right);\,\theta^{*}\right),\qquad y(0;\,\theta^{*})=y_{0}\label{eq:11-2-1}
\end{equation}
where $x\in\left[0,\,1\right]$, $\left|y_{0}\right|\leq C_{0}$,
and $f$ is parameterized by a finite $K-$dimensional vector of coefficients
$\theta\in\mathbb{B}_{q}\left(1\right)$, $q\in\left[1,\,\infty\right]$.
Suppose we have noisy measurements ($Y$) of $y$ sampled at $n$
points $\left\{ x_{i}\right\} _{i=1}^{n}$. If one can solve for $y\left(x;\,\theta^{*},\,y_{0}^{*}\right)$
in an explicit form from (\ref{eq:11-2-1}), then the estimator $\hat{y}_{NLS}=y\left(x;\,\hat{\theta},\,\hat{y}_{0}\right)$
can be formed with the plain vanilla nonlinear least squares (NLS)
\begin{equation}
\left(\hat{\theta},\,\hat{y}_{0}\right)\in\arg\min_{\tilde{\theta}\in\mathbb{B}_{q}\left(1\right),\,\left|\tilde{y}_{0}\right|\leq C_{0}}\frac{1}{2n}\sum_{i=1}^{n}\left(Y_{i}-y\left(x_{i};\,\tilde{\theta},\,\tilde{y}_{0}\right)\right)^{2}.\label{eq:20-1}
\end{equation}
In many ODEs, the solutions are in implicit forms so it is not possible
to write down (\ref{eq:20-1}). In the situation where a quality estimator
$\hat{y}_{0}$ of $y_{0}$ is available, it is natural to exploit
the \textit{Picard iteration} as below:
\begin{equation}
y_{r+1}\left(x;\,\theta\right)=y_{0}+\int_{0}^{x}f\left(s,\,y_{r}\left(s;\,\theta\right);\,\theta\right)ds,\quad\textrm{integer }r\geq0,\,y\left(0;\,\theta\right)=y_{0}.\label{eq:26-1}
\end{equation}
An estimator based on (\ref{eq:26-1}) performs the following steps:
first, we compute
\begin{eqnarray}
\hat{y}_{1}\left(x;\,\theta\right) & = & \hat{y}_{0}+\int_{0}^{x}f\left(s,\,\hat{y}_{0};\theta\right)ds,\nonumber \\
\hat{y}_{2}\left(x;\,\theta\right) & = & \hat{y}_{0}+\int_{0}^{x}f\left(s,\,\hat{y}_{1}\left(s;\,\theta\right);\,\theta\right)ds,\label{eq:26}\\
 & \vdots\nonumber \\
\hat{y}_{R+1}\left(x;\,\theta\right) & = & \hat{y}_{0}+\int_{0}^{x}f\left(s,\,\hat{y}_{R}\left(s;\,\theta\right);\,\theta\right)ds;\nonumber 
\end{eqnarray}
second, we solve the following program
\begin{equation}
\hat{\theta}\in\arg\min_{\theta\in\mathbb{B}_{q}\left(1\right)}\frac{1}{2n}\sum_{i=1}^{n}\left(Y_{i}-\hat{y}_{R+1}\left(x_{i};\theta\right)\right)^{2};\label{eq:20}
\end{equation}
third, we form the estimator $\hat{y}_{Picard}=\hat{y}_{R+1}\left(x;\,\hat{\theta},\,\hat{y}_{0}\right)$. 

For completeness, \textbf{Proposition C.1} in Section \ref{subsec:Proposition-C.1}
and \textbf{Proposition C.2} in Section \ref{subsec:Proposition-C.2}
of the supplementary materials establish upper bounds on the average
squared errors of $\hat{y}_{NLS}$ and $\hat{y}_{Picard}$ with high
probability guarantees, respectively. We choose not to highlight these
propositions here as the abovementioned procedures are not widely
used in practice for reasons discussed in the following. \\
\\
\textbf{Issues with the plain vanilla methods and remedies}. The methods
discussed above have serious practical issues. Clearly it is not possible
to apply (\ref{eq:20-1}) if $y\left(x;\,\theta^{*},\,y_{0}\right)$
cannot be solved in an explicit form from (\ref{eq:11-2-1}), while
the implementation of (\ref{eq:20}) requires a quality estimator
$\hat{y}_{0}$ of $y_{0}$. Even if these methods can be applied,
they suffer numerous computational issues. For both (\ref{eq:20-1})
and (\ref{eq:20}), the objective functions can be highly nonconvex.
The integrals can be hard to compute analytically in (\ref{eq:26});
one could use numerical integrations to approximate the integrals
in (\ref{eq:26}) at a given $\theta$, but doing so would mean that
we have to perform $R+1$ numerical integrations for every single
candidate $\theta$. Such a method is not efficient computationally. 

Because of these issues, researchers have resorted to function fitting
techniques such as splines and alternative smoothing methods to obtain
estimates for $y$ and $y^{'}$. As we have discussed in Section 1,
this strategy typically assumes some smooth class that contains the
solutions and seeks a close enough estimator from the restricted class
of smooth functions. For studying the theoretical behavior of algorithms
used for recovering ODE solutions, \textbf{\textit{Question 1}} raised
in Section 1 is fundamental; moreover, answers to \textbf{\textit{Question
1}} can be useful for studying\textbf{\textit{ Question 2}} (for which
Theorem 2.1 has provided a partial answer). We will show how to exploit
our answers to \textbf{\textit{Question 1}} to improve the general
bound (\ref{eq:general_bound}) when $\mathcal{F}$ corresponds to
a smooth function class. The next section examines these questions
in depth. 

\section{First order ODEs with smooth $f$}

Given that $f$ is often a smooth function in many ODE applications
and prior knowledge on the smoothness structure of $f$ is available,
researchers typically use such information to guide the choice of
$\Pi$ in (\ref{eq:least squares}). This fact motivates us to study
in Section 3.2 the properties of the class of solutions $y$ when
$f$ ranges over the class of functions with various degrees of smoothness.
Our results reveal several sharp contrasts between the standard $(\beta+2)-$degree
smooth class of univariate functions and the $(\beta+2)-$degree smooth
solution classes (of also univariate functions). Before Section 3.2,
we provide a summary of results for the standard smooth class of univariate
functions and their implications in Section 3.1. 

\subsection{A summary of results on standard smooth functions}

Let $p=\left(p_{1},...,\,p_{d}\right)$ and $\left[p\right]=\sum_{k=1}^{d}p_{k}$
where $p_{k}$s are non-negative integers. We write 
\[
D^{p}h\left(z_{1},\,...,\,z_{d}\right):=\partial^{\left[p\right]}h/\partial z_{1}^{p_{1}}...\partial z_{d}^{p_{d}}.
\]
\textbf{Definition} (a typical smooth class). Given a non-negative
integer $\gamma$, we let $\mathcal{S}_{\gamma+1,\,d}\left(\rho,\,\left[\underline{a},\,\overline{a}\right]^{d}\right)$
denote the class of functions such that any function $h\in\mathcal{S}_{\gamma+1,\,d}\left(\rho,\,\left[\underline{a},\,\overline{a}\right]^{d}\right)$
satisfies: (1) $h$ is continuous on $\left[\underline{a},\,\overline{a}\right]^{d}$,
and all partial derivatives of $h$ exist for all $p$ with $\left[p\right]\leq\gamma$;
(2) $\left|D^{p}h\left(X\right)\right|\leq\rho$ for all $X\in\left[\underline{a},\,\overline{a}\right]^{d}$
and all $p$ with $\left[p\right]\leq\gamma$, where $D^{0}h\left(X\right)=h\left(X\right)$
and $\rho$ is a constant independent of the smoothness parameter
$\gamma$; (3) $\left|D^{p}h(X)-D^{p}(X^{'})\right|\leq\rho\left|X-X^{'}\right|_{\infty}$
for all $X,\,X^{'}\in\left[\underline{a},\,\overline{a}\right]^{d}$
and all $p$ with $\left[p\right]=\gamma$. When $\rho=1$, $d=1$,
$\underline{a}=0$ and $\overline{a}=1$, we use the shortform $\mathcal{S}_{\gamma+1}:=\mathcal{S}_{\gamma+1,\,1}\left(1,\,\left[0,\,1\right]\right)$,
the standard smooth class of univariate functions. \\
\textbf{}\\
\textbf{Kolmogorov and Tikhomirov (1959)}. For $\mathcal{S}_{\gamma+1}$,
the standard smooth class of degree $\gamma+1$ on $\left[0,\,1\right]$,
the best upper and lower bounds (to our knowledge) for $\log N_{\infty}\left(\delta,\,\mathcal{S}_{\gamma+1}\right)$
due to Kolmogorov and Tikhomirov \cite{Kolmogorov and Tikhomirov 1959}
take the forms
\begin{eqnarray}
\log N_{\infty}\left(\delta,\,\mathcal{S}_{\gamma+1}\right) & \precsim & \delta^{\frac{-1}{\gamma+1}}+\left(\gamma+1\right)\log\frac{1}{\delta},\label{eq:Kolmogorov_upper}\\
\log N_{\infty}\left(\delta,\,\mathcal{S}_{\gamma+1}\right) & \succsim & \delta^{\frac{-1}{\gamma+1}}.\label{eq:Kolmogorov_lower}
\end{eqnarray}
Kolmogorov and Tikhomirov \cite{Kolmogorov and Tikhomirov 1959} derive
the lower bound on the packing number $M_{\infty}\left(\delta,\,\mathcal{S}_{\gamma+1}\right)$
of $\mathcal{S}_{\gamma+1}$ and (\ref{eq:Kolmogorov_lower}) comes
from the fact that $M_{\infty}\left(2\delta,\,\mathcal{S}_{\gamma+1}\right)\leq N_{\infty}\left(\delta,\,\mathcal{S}_{\gamma+1}\right)$.
Upper bound (\ref{eq:Kolmogorov_upper}) is useful for analyzing the
``local complexity'' associated with $\mathcal{S}_{\gamma+1}$,
as we discuss below.\\
\textbf{}\\
\textbf{Definition} (local complexity). Given a radius $\tilde{r}_{n}>0$
and a function class $\bar{\mathcal{F}}$, define the \textit{local
complexity 
\begin{equation}
\mathcal{G}_{n}(\tilde{r}_{n};\,\bar{\mathcal{F}}):=\mathbb{E}_{\varepsilon}\left[\sup_{h\in\Lambda(\tilde{r}_{n};\,\bar{\mathcal{F}})}\left|\frac{1}{n}\sum_{i=1}^{n}\varepsilon_{i}h(x_{i})\right|\right],\label{eq:13-1}
\end{equation}
}where $\varepsilon_{i}\overset{i.i.d.}{\sim}\mathcal{N}\left(0,\,1\right)$,
$\varepsilon=\left\{ \varepsilon_{i}\right\} _{i=1}^{n}$, $\Lambda(\tilde{r}_{n};\,\bar{\mathcal{F}})=\left\{ h\in\bar{\mathcal{F}}:\,\left|h\right|_{n}\leq\tilde{r}_{n}\right\} $,
and $\left|h\right|_{n}:=\sqrt{\frac{1}{n}\sum_{i=1}^{n}\left(h(x_{i})\right)^{2}}$.
Bounds on local complexity are often used to analyze the distance
between an unknown function from a class and an estimator of this
function constrained in the same class. Therefore, the underlying
$\bar{\mathcal{F}}$ in (\ref{eq:13-1}) takes the following form
\begin{equation}
\bar{\mathcal{F}}:=\left\{ g=g_{1}-g_{2}:\,g_{1},\,g_{2}\in\mathcal{F}\right\} .\label{eq:28}
\end{equation}
Various papers have studied localized forms of the Rademacher complexity
(where $\varepsilon_{i}$ in (\ref{eq:13-1}) is a Rademacher random
variable) and Gaussian complexity (where $\varepsilon_{i}$ in (\ref{eq:13-1})
is a normal random variable); see, e.g., \cite{Bartlett and Mendelson 2002 (main text citation),Koltchinskii 2001,koltchinskii 2006 (main text citation),van de Geer 2000 (main text citation)}.\\
\textbf{}\\
\textbf{Implications on noisy fitting}. In what follows, let us consider
an application of (\ref{eq:13-1}). Suppose that one only observes
noisy measurements $Y_{i}$ of $g^{*}\left(x_{i}\right)$ in the following
form: 

\begin{equation}
Y_{i}=g^{*}(x_{i})+\varepsilon_{i},\quad i=1,...,n,\label{eq:model-1}
\end{equation}
where $\left\{ x_{i}\right\} _{i=1}^{n}$ is a collection of fixed
design points and the unobserved noise terms $\left\{ \varepsilon_{i}\right\} _{i=1}^{n}$
are i.i.d. draws from $\mathcal{N}\left(0,\,1\right)$. The assumption
$\textrm{var}\left(\varepsilon_{i}\right)=1$ is simply to facilitate
the exposition (so that we can focus on the impact of smoothness)
and can be relaxed. The normality assumption can also be relaxed,
for instance, to allow $\varepsilon_{i}$ only subject to a sub-Gaussian
tail condition in both (\ref{eq:13-1}) and (\ref{eq:model-1}). Examples
of sub-Gaussian random variables include bounded variables, mixture
of normal and bounded variables, etc. 

Suppose that $g^{*}(\cdot)\in\mathcal{\mathcal{F}}$. To find a function
$\hat{g}(\cdot)$ that fits $\left(Y_{i},\,x_{i}\right)$, a classical
approach is based on the least squares 
\begin{equation}
\hat{g}\in\arg\min_{\tilde{g}\in\mathcal{F}}\frac{1}{2n}\sum_{i=1}^{n}\left(Y_{i}-\tilde{g}\left(x_{i}\right)\right)^{2}.\label{eq:least squares-1}
\end{equation}
Because $g^{*}\in\mathcal{F}$ ($g^{*}$ is feasible) and $\hat{g}$
is optimal, we have 
\begin{equation}
\frac{1}{n}\sum_{i=1}^{n}\left(\hat{g}(x_{i})-g^{*}(x_{i})\right)^{2}\leq\frac{2}{n}\sum_{i=1}^{n}\varepsilon_{i}\left(\hat{g}(x_{i})-g^{*}(x_{i})\right).\label{eq:44-1}
\end{equation}
Bounding the mean squared error $\mathbb{E}\left[\frac{1}{n}\sum_{i=1}^{n}\left(\hat{g}(x_{i})-g^{*}(x_{i})\right)^{2}\right]$
from above can be reduced to: (1) seeking a sharp enough bound $\mathcal{U}_{n}\left(\tilde{r}_{n};\,\bar{\mathcal{F}}\right)$
on $\mathcal{G}_{n}(\tilde{r}_{n};\,\bar{\mathcal{F}})$, and (2)
seeking $\tilde{r}_{n}>0$ that satisfies 
\begin{equation}
\tilde{r}_{n}^{2}\asymp\mathcal{U}_{n}\left(\tilde{r}_{n};\,\bar{\mathcal{F}}\right).\label{eq:critical}
\end{equation}
We refer the readers to the textbook by \cite{wainright 2019} (Chapter
13, 2019) for a detailed explanation on how this argument works.

It is well known that the local complexity $\mathcal{G}_{n}\left(\tilde{r}_{n};\,\bar{\mathcal{F}}\right)$
in (\ref{eq:13-1}) can be bounded by exploiting the Dudley's entropy
integral (see, e.g., \cite{Ledoux and Talagrand 1991 (main text),van de Geer 2000 (main text citation),wainright 2019}).
In the case where $\mathcal{F}=\mathcal{S}_{\gamma+1}$ in (\ref{eq:28}),
the entropy integral approach gives that
\[
\mathcal{G}_{n}(\tilde{r}_{n};\,\bar{\mathcal{F}})\leq\frac{16}{\sqrt{n}}\int_{\frac{\tilde{r}_{n}^{2}}{4}}^{\tilde{r}_{n}}\sqrt{\log N_{\infty}\left(\delta,\,\bar{\mathcal{F}}\right)}d\delta+\frac{\tilde{r}_{n}^{2}}{4}\leq\underset{\mathcal{U}_{n}\left(\tilde{r}_{n};\,\bar{\mathcal{F}}\right)}{\underbrace{c\left(\tilde{r}_{n}\sqrt{\frac{\gamma+1}{n}}+\frac{1}{\sqrt{n}}\tilde{r}_{n}^{\frac{2\gamma+1}{2\gamma+2}}\right)+\frac{\tilde{r}_{n}^{2}}{4}}}
\]
where the second inequality follows from (\ref{eq:Kolmogorov_upper}).
Setting $\mathcal{U}_{n}\left(\tilde{r}_{n};\,\bar{\mathcal{F}}\right)=\tilde{r}_{n}^{2}$
yields 
\begin{equation}
\tilde{r}_{n}^{2}\asymp\max\left\{ \frac{\gamma\vee1}{n},\,\left(\frac{1}{n}\right)^{\frac{2\left(\gamma+1\right)}{2\left(\gamma+1\right)+1}}\right\} =\begin{cases}
\frac{\gamma\vee1}{n}, & n\precsim\left(\gamma\vee1\right)^{2\gamma+3}\\
\left(\frac{1}{n}\right)^{\frac{2\left(\gamma+1\right)}{2\left(\gamma+1\right)+1}}, & n\succsim\left(\gamma\vee1\right)^{2\gamma+3}.
\end{cases}\label{eq:28-1}
\end{equation}
Therefore, an upper bound on $\mathbb{E}\left[\frac{1}{n}\sum_{i=1}^{n}\left(\hat{g}(x_{i})-g^{*}(x_{i})\right)^{2}\right]$
scales roughly as $\frac{\gamma\vee1}{n}$ when $n\precsim\left(\gamma\vee1\right)^{2\gamma+3}$
and $\left(\frac{1}{n}\right)^{\frac{2\left(\gamma+1\right)}{2\left(\gamma+1\right)+1}}$
when $n\succsim\left(\gamma\vee1\right)^{2\gamma+3}$. 

We can also derive (\ref{eq:28-1}) using the approach due to Mendelson
\cite{Mendelson 2002}, which bounds the local complexity through
eigenvalues of the empirical kernel matrix associated with a reproducing
kernel Hilbert space (RKHS). In particular, we consider the following
RKHS (that contains $\mathcal{S}_{\gamma+1}$) 
\begin{equation}
\mathcal{H}_{\gamma+1}=\left\{ f:\,\left[0,\,1\right]\rightarrow\mathbb{R}\vert f^{(\gamma)}\textrm{ is absolutely continuous and }\int_{0}^{1}\left[f^{(\gamma+1)}\left(t\right)\right]^{2}dt\leq1\right\} \label{eq:sobolev}
\end{equation}
and then apply the Mendelson approach to arrive at (\ref{eq:28-1}).
Suppose $n\geq\gamma+1$. The scaling $\frac{\gamma\vee1}{n}$ is
associated with the space of polynomials of degree $\gamma$, and
the bound $\left(\frac{1}{n}\right)^{\frac{2\left(\gamma+1\right)}{2\left(\gamma+1\right)+1}}$
is associated with the space 
\begin{align*}
\mathcal{H}_{\gamma+1,0} & =\{f:\,\left[0,\,1\right]\rightarrow\mathbb{R}\vert\textrm{ for all non-negative integer }k\leq\gamma,\,f^{(k)}\left(0\right)=0,\\
 & f^{(\gamma)}\textrm{ is abs. cont. and }\int_{0}^{1}\left[f^{(\gamma+1)}\left(t\right)\right]^{2}dt\leq1\}.
\end{align*}
This can be seen from the fact that $f\in\mathcal{H}_{\gamma+1}$
has the expansion $f\left(x\right)=\sum_{k=0}^{\gamma}f^{(k)}(0)\frac{x^{k}}{k!}+\int_{0}^{1}f^{(\gamma+1)}(t)\frac{\left(x-t\right)_{+}^{\gamma}}{\gamma!}dt$,
where $(w)_{+}=w\vee0$. The scheme above is important as it allows
one to implement (\ref{eq:least squares-1}) via spline bases and
kernel functions, very popular techniques for the problems of fitting
and interpolating functions (see, e.g., \cite{Berlinet and Thomas-Agnan 2004,Scholkopf and smola 2002,Wahba 1990}).
For example, if $\mathcal{F}=\mathcal{H}_{2}$ in (\ref{eq:least squares-1}),
we have the cubic spline estimator $\hat{g}(x)=\hat{\alpha}_{0}+\hat{\alpha}_{1}x+\frac{1}{\sqrt{n}}\sum_{i=1}^{n}\hat{\pi}_{i}\mathcal{K}_{1}\left(x,\,x_{i}\right)$,
where $\mathcal{K}_{1}\left(x_{i},\,x_{j}\right)=\int_{0}^{1}\left(x_{i}-t\right)_{+}\left(x_{j}-t\right)_{+}dt$
and 
\begin{align}
\left(\hat{\alpha},\,\hat{\pi}\right) & =\arg\min_{\left(\alpha,\pi\right)\in\mathbb{R}^{2}\times\mathbb{R}^{n}}\left\{ \frac{1}{2n}\left|Y-Z\alpha-\sqrt{n}\mathbb{K}_{1}\pi\right|_{2}^{2}\right\} ,\nonumber \\
\textrm{s.t.} & \pi^{T}\mathbb{K}_{1}\pi\leq1,\quad\forall k=0,...,\beta+1,\label{eq:standard_krr}
\end{align}
with $\alpha=\left(\alpha_{0},\,\alpha_{1}\right)^{T}$, $Z$ an $n\times2$
matrix whose $i$th row is $\left(1,\,x_{i}\right)$, and $\mathbb{K}_{1}$
the $n\times n$ kernel matrix consisting of entries $\mathcal{K}_{1}\left(x_{i},\,x_{j}\right)$.

It is worth pointing out that in the broad statistics, machine learning,
and econometrics community, researchers generally take $O\left(n^{\frac{-2\left(\gamma+1\right)}{2\left(\gamma+1\right)+1}}\right)$
as the upper bound for $\mathbb{E}\left[\frac{1}{n}\sum_{i=1}^{n}\left(\hat{g}(x_{i})-g^{*}(x_{i})\right)^{2}\right]$
in terms of (\ref{eq:least squares-1}). This practice assumes $\log N_{\infty}\left(\delta,\,\mathcal{S}_{\gamma+1}\right)\asymp\delta^{\frac{-1}{\gamma+1}}$
based on (\ref{eq:Kolmogorov_upper}) and (\ref{eq:Kolmogorov_lower}).
In deriving (\ref{eq:Kolmogorov_upper}), Kolmogorov and Tikhomirov
\cite{Kolmogorov and Tikhomirov 1959} partition $\left[0,\,1\right]$
into $s\asymp\delta^{\frac{1}{\gamma+1}}$ intervals and consider
a grid of points $\left(x_{0},\,x_{1},\,...,\,x_{s}\right)$. The
term $\left(\gamma+1\right)\log\frac{1}{\delta}$ in (\ref{eq:Kolmogorov_upper})
comes from counting the number of possible values of $\left(\left\lfloor \frac{f^{(k)}\left(x_{0}\right)}{\delta_{k}}\right\rfloor ,\,k=0,...,\gamma,\,f\in\mathcal{S}_{\gamma+1}\right)$
given the prespecified accuracy $\delta_{k}$ for each $k$th derivative
($\delta_{0}=\delta$). Assuming $\left(\gamma+1\right)\log\frac{1}{\delta}$
is negligible has a similar effect on $\log N_{\infty}\left(\delta,\,\mathcal{S}_{\gamma+1}\right)$
as assuming $f^{(k)}\left(0\right)=0$ for all $k\leq\gamma$ and
$f\in\mathcal{S}_{\gamma+1}$.\footnote{On a related note, we conjecture that the lower bound (\ref{eq:Kolmogorov_lower})
can be sharpened based on the construction in \cite{Kolmogorov and Tikhomirov 1959}.} 

For small $\gamma$, (\ref{eq:Kolmogorov_upper}) and (\ref{eq:Kolmogorov_lower})
would be comparable; nonasymptotically, for large enough (still finite)
$\gamma$, $\left(\gamma+1\right)\log\frac{1}{\delta}$ can be greater
than $\delta^{\frac{-1}{\gamma+1}}$ even with small $\delta$. Figure
1 exhibits the growth of $\left(\gamma\vee1\right)\log\frac{1}{\delta}$
and $\delta^{\frac{-1}{\gamma+1}}$ as $\gamma$ increases for $\delta=0.01$
and $\delta=0.001$. Taking $n^{\frac{-2\left(\gamma+1\right)}{2\left(\gamma+1\right)+1}}$
as the upper bound for the mean squared error of (\ref{eq:least squares-1})
essentially assumes that the sample size $n\succsim\left(\gamma\vee1\right)^{2\gamma+3}$
as in (\ref{eq:28-1}). \\
\textbf{}\\
\textbf{Information theoretic lower bounds and minimax optimality}.
Under the assumption $\left(\gamma+1\right)\log\frac{1}{\delta}\precsim\delta^{\frac{-1}{\gamma+1}}$,
one could use a Yang-Barron version of Fano's method (see \cite{yang and barron 1999})
to show that 

\[
\inf_{\tilde{g}}\sup_{g\in\mathcal{S}_{\gamma+1}}\mathbb{E}\left[\frac{1}{n}\sum_{i=1}^{n}\left[\tilde{g}\left(x_{i}\right)-g\left(x_{i}\right)\right]^{2}\right]\succsim\min\left\{ 1,\,\left(\frac{1}{n}\right)^{\frac{2\left(\gamma+1\right)}{2\left(\gamma+1\right)+1}}\right\} ;
\]
see \cite{wainright 2019} (Chapter 15, 2019) for more details on
this derivation. The bound above implies the rate $\left(\frac{1}{n}\right)^{\frac{2\left(\gamma+1\right)}{2\left(\gamma+1\right)+1}}$
is minimax optimal. Again, we emphasize that this derivation assumes
$\left(\gamma+1\right)\log\frac{1}{\delta}\leq\delta^{\frac{-1}{\gamma+1}}$
in (\ref{eq:Kolmogorov_upper}).\\

Considering $\gamma=\beta+1$, in sum, under the assumption $n\succsim\left(\beta\vee1\right)^{2\beta+5}$,
the minimax optimal rate $\left(\frac{1}{n}\right)^{\frac{2\left(\beta+2\right)}{2\left(\beta+2\right)+1}}$
(which decreases in $\beta$) coincides with the convergence rate
of the least squares. We can view this requirement on the sample size,
$n\succsim\left(\beta\vee1\right)^{2\beta+5}$, a ``curse of smoothness''
for least squares estimations of elements in $\mathcal{S}_{\beta+2}$.
In our subsequent results for ODE solution classes, we will show that
the ``curse of smoothness'' is exacerbated. 

\textbf{}
\begin{figure}

\textbf{\caption{{\small{}``$\ocircle$'': $\left(\gamma\vee1\right)\log\frac{1}{\delta}$;
``$\boxempty$'': $\delta^{\frac{-1}{\gamma+1}}$. (a): $\delta=0.01$;
(b): $\delta=0.001$}}
}

\includegraphics[scale=0.32]{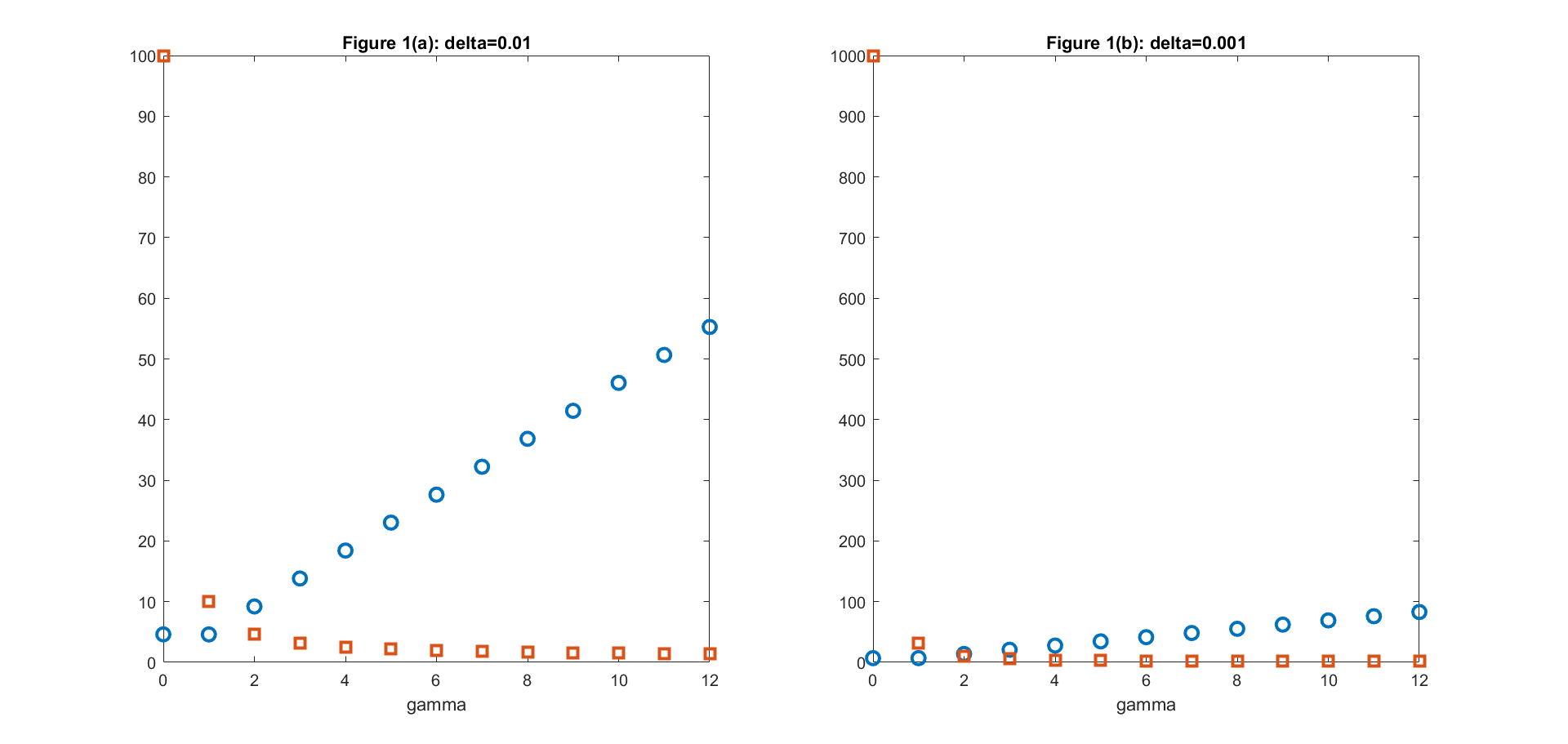}

\end{figure}

\subsection{Results for ODE solution classes}

In this subsection, we consider the ODEs in the forms \textit{
\begin{eqnarray}
\textrm{autonomous ODE: }y^{'}\left(x\right) & = & f\left(y\left(x\right)\right),\qquad y\left(x_{0}\right)=y_{0},\label{eq:11-4}\\
\textrm{nonautonomous ODE: }y^{'}\left(x\right) & = & f\left(x,\,y\left(x\right)\right),\qquad y\left(x_{0}\right)=y_{0}.\label{eq:11-4-1}
\end{eqnarray}
}Extensions of the analyses in this section to higher order ODEs are
certainly possible but are very laborious. 

\subsubsection{Nonasymptotic bounds on derivatives and covering numbers}

We have established a general bound on the covering numbers of ODE
solution classes in Theorem 2.1, which states that the covering number
of the underlying solution class is bounded from above by the product
of the covering number of the class $\mathcal{F}$ that $f$ ranges
over and the covering number of the set where initial values lie.
For the nonautonomous ODE (\ref{eq:11-4-1}), it appears there is
still room to improve upon (\ref{eq:general_bound}). After all, $y$
is a univariate function even though $f$ in (\ref{eq:11-4-1}) is
a bivariate function (on the other hand, improvement upon (\ref{eq:general_bound})
may be more limited for the autonomous ODE (\ref{eq:11-4})). This
observation motivates the following lemmas which build upon or extend
the classical Picard smoothness structure (\ref{eq:lip}) to higher
degree smoothness structures. These results are not only useful for
potentially tightening bounds on the covering numbers but also provide
guidance on the design of function fitting algorithms (as we will
see in Section 3.2.2).\\
\\
\textbf{Lemma 3.1} (autonomous ODE). \textit{Let $\beta$ be a non-negative
integer. In terms of }(\ref{eq:11-4})\textit{, assume that $f$ is
is continuous on $\left[y_{0}-b,\,y_{0}+b\right]$ and $\beta-$times
differentiable, and for all $y,\,\tilde{y}\in\left[y_{0}-b,\,y_{0}+b\right]$,
\begin{eqnarray}
\left|f^{(k)}(y)\right| & \leq & 1,\quad\forall\,0\leq k\leq\beta,\label{eq:a2}\\
\left|f^{\beta}(y)-f^{\beta}(\tilde{y})\right| & \leq & |y-\tilde{y}|,\label{eq:a3}
\end{eqnarray}
where $f^{(0)}=f$. Suppose $y\left(x\right)\in\left[y_{0}-b,\,y_{0}+b\right]$
for $x\in\left[x_{0}-a,\,x_{0}+a\right]$. }

\textit{(i) For all $x,\,x^{'}\in\left[x_{0}-\alpha,\,x_{0}+\alpha\right]$
with $\alpha=a\wedge b$, we have
\[
\left|y^{(k)}(x)\right|\leq\left(k-1\right)!
\]
 for all $0<k\leq\beta+1$ and 
\[
\left|y^{(\beta+1)}(x)-y^{(\beta+1)}(x^{'})\right|\leq\left(\beta+1\right)!|x-x^{'}|.
\]
}

\textit{(ii) There exists an ODE with a solution such that the absolute
value of the $kth$ derivative of this solution equals $\left(k-1\right)!$
for all $k=1,...,\beta+2$. For example, $y^{'}=e^{-y-\frac{1}{2}}$
is such an ODE.}\footnote{Note that $f(y)=e^{-y-\frac{1}{2}}$ is continuous on $\left[-\frac{1}{2},\,\frac{1}{2}\right]$,
$\beta-$times differentiable, and satisfies (\ref{eq:a2})--(\ref{eq:a3})
for all $y\in\left[-\frac{1}{2},\,\frac{1}{2}\right]$.}\textit{}\\
\textit{}\\
\textbf{Lemma 3.2 }(nonautonomous ODE).\textit{ In terms of (\ref{eq:11-4-1}),
assume that $f$ is continuous on $\Upsilon=\left[x_{0}-a,\,x_{0}+a\right]\times\left[y_{0}-b,\,y_{0}+b\right]$
and all partial derivatives $D^{p}$ of $f$ exist for all $p$ with
$\left[p\right]=p_{1}+p_{2}\leq\beta$; $\left|D^{p}f\left(x,\,y\right)\right|\leq1$
for all $p$ with $\left[p\right]\leq\beta$ and $\left(x,\,y\right)\in\Upsilon$,
where $D^{0}f\left(x,\,y\right)=f\left(x,\,y\right)$; and 
\begin{equation}
\left|D^{p}f\left(x,\,y\right)-D^{p}f\left(x^{'},\,\tilde{y}\right)\right|\leq\max\left\{ \left|x-x^{'}\right|,\,\left|y-\tilde{y}\right|\right\} \qquad\forall\left(x,\,y\right),\,\left(x^{'},\,\tilde{y}\right)\in\Upsilon,\label{eq:a7-1}
\end{equation}
for all $p$ with $\left[p\right]=\beta$. Then for all $x,\,x^{'}\in\left[x_{0}-\alpha,\,x_{0}+\alpha\right]$
with $\alpha=a\wedge b$, we have
\[
\left|y^{(k)}(x)\right|\leq2^{k-1}(k-1)!
\]
 for all $0\leq k\leq\beta+1$ and 
\[
\left|y^{(\beta+1)}(x)-y^{(\beta+1)}(x^{'})\right|\leq2^{\beta+1}\left(\beta+1\right)!|x-x^{'}|.
\]
}

The proofs for Lemmas 3.1 and 3.2 are provided in Section \ref{subsec:Lemmas 3.1 and 3.2}
of the supplementary materials. Note that our results obviously hold
if (\ref{eq:a3}) and (\ref{eq:a7-1}) are replaced with stronger
conditions $\left|f^{(\beta+1)}(y)\right|\leq1$ and $\left|D^{p}f\left(x,\,y\right)\right|\leq1$
for all $p$ with $\left[p\right]=\beta+1$, respectively. Recalling
\textit{Question 1}, now we see from Lemma 3.1(ii) that the derivatives
of $y$ may not stay ``nicely'' bounded even if the absolute values
of the derivatives of $f$ are all bounded by $1$. Based on Lemmas
3.1--3.2 and Theorem 2.1, the following two theorems reveal that
our upper bounds for the covering numbers of the $(\beta+2)-$degree
smooth solution classes are greater than those of the standard $(\beta+2)-$degree
smooth class of univariate functions. 

Let us first introduce several definitions used throughout the rest
of this section. To lighten notations, from now on, let us assume
the initial condition $y\left(0\right)=y_{0}$ in (\ref{eq:11-4})--(\ref{eq:11-4-1}),
with $\left|y_{0}\right|\leq C_{0}$. Given $C_{0}>0$, $b>0$, and
$\alpha=\min\left\{ 1,\,b\right\} $, let $\overline{C}:=C_{0}+b$,
$\Xi:=\left[0,\,1\right]\times\left[-\overline{C},\,\overline{C}\right]$,
and $L_{\max}:=\sup_{x\in\left[0,\,\alpha\right]}\left\{ \exp\left(x\right)\left[1+\int_{0}^{x}\exp\left(-s\right)ds\right]\right\} $.
\\
\textbf{}\\
\textbf{Theorem 3.1}\textit{ }(autonomous ODE).\textit{ In terms of
(\ref{eq:11-4}), assume $f$ ranges over $\mathcal{S}_{\beta+1,\,1}\left(1,\,\left[-\overline{C},\,\overline{C}\right]\right)$.
For a given $\frac{\delta}{5}\in\left(0,\,\alpha\right)$, we have
\begin{eqnarray}
\log N_{\infty}\left(\delta,\,\mathcal{Y}\right) & \leq & \min\left\{ \min_{\gamma\in\left\{ 0,...,\beta\right\} }Z_{1}(\delta,\gamma),\,\min_{\gamma\in\left\{ 0,...,\beta\right\} }Z_{2}\left(\frac{\delta}{L_{\max}},\gamma\right)\right\} ,\label{eq:14-4}\\
\log N_{\infty}\left(\delta,\,\mathcal{Y}_{1}\right) & \leq & \min_{\gamma\in\left\{ 0,...,\beta\right\} }Z_{3}(\delta,\gamma),\label{eq:14-4b}
\end{eqnarray}
where $\mathcal{Y}$ is the class of solutions (to (\ref{eq:11-4})
with $f\in\mathcal{S}_{\beta+1,\,1}\left(1,\,\left[-\overline{C},\,\overline{C}\right]\right)$)
on $\left[0,\,\alpha\right]$, $\mathcal{Y}_{1}$ is the class consisting
of the first derivative $y^{'}$ of $y\in\mathcal{Y}$, and
\begin{eqnarray*}
Z_{1}(\delta,\gamma) & = & \log\left(\prod_{i=0}^{\gamma}i!\right)+\frac{\gamma+3}{2}\log\frac{5}{\delta}+\alpha\left(\frac{\delta}{5}\right)^{\frac{-1}{\gamma+2}}\log2+\log\left(4\overline{C}\right),\\
Z_{2}(\delta,\gamma) & = & \frac{\gamma+2}{2}\log\frac{5}{\delta}+2\overline{C}\left(\log2\right)\left(\frac{\delta}{5}\right)^{\frac{-1}{\gamma+1}}+\log4+\log\left(\frac{C_{0}}{\delta}+1\right),\\
Z_{3}(\delta,\gamma) & = & \log\left(\prod_{i=0}^{\gamma}i!\right)+\frac{\gamma+2}{2}\log\frac{5}{\delta}+\alpha\left(\frac{\delta}{5}\right)^{\frac{-1}{\gamma+1}}\log2+\log4.
\end{eqnarray*}
}\textbf{Theorem 3.2}\textit{ }(nonautonomous ODE).\textit{ In terms
of (\ref{eq:11-4-1}), assume $f$ ranges over $\mathcal{S}_{\beta+1,\,2}\left(1,\,\Xi\right)$.
For a given $\frac{\delta}{5}\in\left(0,\,\alpha\right)$, we have
\begin{eqnarray}
\log N_{\infty}\left(\delta,\,\mathcal{Y}\right) & \leq & \min\left\{ \min_{\gamma\in\left\{ 0,...,\beta\right\} }W_{1}(\delta,\gamma),\,\min_{\gamma\in\left\{ 0,...,\beta\right\} }W_{2}\left(\frac{\delta}{L_{\max}},\gamma\right)\right\} ,\label{eq:14-4-2}\\
\log N_{\infty}\left(\delta,\,\mathcal{Y}_{1}\right) & \leq & \min_{\gamma\in\left\{ 0,...,\beta\right\} }W_{3}(\delta,\gamma)\label{eq:14-4-2b}
\end{eqnarray}
where $\mathcal{Y}$ is the class of solutions (to (\ref{eq:11-4-1})
with $f\in\mathcal{S}_{\beta+1,\,2}\left(1,\,\Xi\right)$) on $\left[0,\,\alpha\right]$,
$\mathcal{Y}_{1}$ is the class consisting of the first derivative
$y^{'}$ of $y\in\mathcal{Y}$, and
\begin{eqnarray*}
W_{1}(\delta,\gamma) & = & \log\left(\prod_{i=0}^{\gamma}i!\right)+\frac{\gamma^{2}+\gamma}{2}\log2+\frac{\gamma+3}{2}\log\frac{5}{\delta}+\alpha\left(\frac{\delta}{5}\right)^{\frac{-1}{\gamma+2}}\log4+\log\left(4\overline{C}\right),\\
W_{2}(\delta,\gamma) & = & \frac{\left(\gamma+2\right)\left(\gamma+3\right)}{6}\log\frac{5}{\delta}+20\left(\overline{C}\vee1\right)\left(\log2\right)\left(\frac{\delta}{5}\right)^{\frac{-2}{\gamma+1}}+4\log2+\log\left(\frac{C_{0}}{\delta}+1\right),\\
W_{3}(\delta,\gamma) & = & \log\left(\prod_{i=0}^{\gamma}i!\right)+\frac{\gamma^{2}+\gamma}{2}\log2+\frac{\gamma+2}{2}\log\frac{5}{\delta}+\alpha\left(\frac{\delta}{5}\right)^{\frac{-1}{\gamma+1}}\log4+\log4.
\end{eqnarray*}
}

The proofs for Theorems 3.1 and 3.2 are provided in Section \ref{subsec:Theorems-3.1-and-3.2}
of the supplementary materials. The bounds related to $Z_{2}$ and
$W_{2}$ are based on Theorem 2.1. The bounds related to $Z_{1}$,
$Z_{3}$, $W_{1}$ and $W_{3}$ are based on Lemmas 3.1--3.2, which
imply that $\mathcal{Y}\subseteq\mathcal{AS}_{\gamma+2,\overline{C}}^{\dagger}$
for the autonomous ODEs and $\mathcal{Y}\subseteq\mathcal{S}_{\gamma+2,\overline{C}}^{\dagger}$
for the nonautonomous ODEs, with $\mathcal{AS}_{\gamma+2,\overline{C}}^{\dagger}$
and $\mathcal{S}_{\gamma+2,\overline{C}}^{\dagger}$ defined as follows.\\
\textbf{}\\
\textbf{Definition} (nonstandard smooth classes) 
\begin{itemize}
\item Given $\overline{C}$ and a non-negative integer $\gamma$, we let
$\mathcal{AS}_{\gamma+2,\overline{C}}^{\dagger}$ denote the class
of functions such that any function $h\in\mathcal{AS}_{\gamma+2,\overline{C}}^{\dagger}$
satisfies the following properties: (1) $h$ is continuous on $\left[0,\,\alpha\right]$
and differentiable $\gamma+1$ times; (2) $\left|h(x)\right|\leq\overline{C}$,
and $\left|h^{(k)}(x)\right|\leq\left(k-1\right)!$ for all $k=1,...,\gamma+1$
and $x\in\left[0,\,\alpha\right]$; (3) $\left|h^{(\gamma+1)}(x)-h^{(\gamma+1)}(x^{'})\right|\leq\left(\gamma+1\right)!\left|x-x^{'}\right|$
for all $x,\,x^{'}\in\left[0,\,\alpha\right]$.
\item We let $\mathcal{S}_{\gamma+2,\overline{C}}^{\dagger}$ denote the
class of functions such that any function $h\in\mathcal{S}_{\gamma+2,\overline{C}}^{\dagger}$
satisfies the following properties: (1) $h$ is continuous on $\left[0,\,\alpha\right]$
and differentiable $\gamma+1$ times; (2) $\left|h(x)\right|\leq\overline{C}$,
and $\left|h^{(k)}(x)\right|\leq2^{k-1}\left(k-1\right)!$ for all
$k=1,...,\gamma+1$ and $x\in\left[0,\,\alpha\right]$; (3) $\left|h^{(\gamma+1)}(x)-h^{(\gamma+1)}(x^{'})\right|\leq2^{\gamma+1}\left(\gamma+1\right)!\left|x-x^{'}\right|$
for all $x,\,x^{'}\in\left[0,\,\alpha\right]$.
\end{itemize}

\subsubsection*{Comparing autonomous ODEs with nonautonomous ODEs}

The function $Z_{2}(\delta,\gamma)$ in Theorem 3.1 associated with
the autonomous ODEs comes from a straightforward application of (\ref{eq:Kolmogorov_upper}),
Theorem 2.1 (specialized to $m=1$)\footnote{In this case, we can set $L_{\max}=\sup_{x\in\left[0,\,\alpha\right]}\left\{ \exp\left(x\right)\left[1+\int_{0}^{x}\exp\left(-s\right)ds\right]\right\} $.},
as well as the fact $\mathcal{S}_{\beta+1}\subseteq\mathcal{S}_{\beta}\subseteq\mathcal{S}_{\beta-1}\cdots$. 

In comparison, first notice that the function $W_{2}(\delta,\gamma)$
in Theorem 3.2 associated with the nonautonomous ODEs involves $\frac{\left(\gamma+2\right)\left(\gamma+3\right)}{6}\log\frac{5}{\delta}$,
which has to do with the fact that a function of two variables has
$\gamma+1$ distinct $\gamma$th partial derivatives (if all $\gamma$th
partial derivatives are continuous). While Kolmogorov and Tikhomirov
\cite{Kolmogorov and Tikhomirov 1959} explicitly derived $\left(\gamma+1\right)\log\frac{1}{\delta}$
for the univariate function class, moving from univariate functions
to multivariate functions, this type of log terms appears unmentioned
in \cite{Kolmogorov and Tikhomirov 1959} possibly because $\gamma$
is assumed to be very small. However, the log term $\frac{\left(\gamma+2\right)\left(\gamma+3\right)}{6}\log\frac{5}{\delta}$
can have practical implications, just like the term $\left(\gamma+1\right)\log\frac{1}{\delta}$
in (\ref{eq:Kolmogorov_upper}) discussed in Section 3.1. 

Second, the derivations of $Z_{2}(\delta,\gamma)$ and $W_{2}(\delta,\gamma)$
are based on Theorem 2.1; therefore, $W_{2}(\delta,\gamma)$ involves
$\left(\frac{\delta}{5}\right)^{\frac{-2}{\gamma+1}}$ (since $f$
in (\ref{eq:11-4-1}) is a function of two variables) while $Z_{2}(\delta,\gamma)$
involves $\left(\frac{\delta}{5}\right)^{\frac{-1}{\gamma+1}}$ (since
$f$ in (\ref{eq:11-4}) is a univariate function). On the other hand,
the derivations of $Z_{1}(\delta,\gamma)$ and $W_{1}(\delta,\gamma)$
are based on Lemmas 3.1 and 3.2, and consequently, $\left(\frac{\delta}{5}\right)^{\frac{-1}{\gamma+2}}$
shows up in both $W_{1}(\delta,\gamma)$ and $Z_{1}(\delta,\gamma)$;
after all, even though $f$ in (\ref{eq:11-4-1}) is a function of
two variables, the solution $y$ is a univariate function, and $\left(\frac{\delta}{5}\right)^{\frac{-1}{\gamma+2}}$
in $W_{1}(\delta,\gamma)$ reflects this characteristic. The term
$\log\left(\prod_{i=0}^{\gamma}i!\right)$ in $Z_{1}(\delta,\gamma)$
of Theorem 3.1 and $W_{1}(\delta,\gamma)$ of Theorem 3.2 comes from
the factorial part in Lemmas 3.1 and 3.2. The term $\frac{\gamma^{2}+\gamma}{2}\log2$
in $W_{1}(\delta,\gamma)$ of Theorem 3.2 comes from the exponential
part in Lemma 3.2. Note that these terms do not appear in (\ref{eq:Kolmogorov_upper})
for the standard smooth classes, and reflect the higher complexity
of $\mathcal{AS}_{\beta+2,\overline{C}}^{\dagger}$ and $\mathcal{S}_{\beta+2,\overline{C}}^{\dagger}$. 

Based on the above comparison, for large enough $\beta$, the nonautonomous
ODE (\ref{eq:11-4-1}) benefits more from $W_{1}(\delta,\gamma)$
than the autonomous ODE (\ref{eq:11-4}) from $Z_{1}(\delta,\gamma)$.
Recall that by Lemma 3.1(ii), the bound $\left(k-1\right)!$ on the
absolute value of the $k$th derivative of $y\left(x\right)$ (for
all $k=1,...,\beta+1$) and the bound $\left(\beta+1\right)!\left|x-x^{'}\right|$
on $\left|y^{(\beta+1)}(x)-y^{(\beta+1)}(x^{'})\right|$ are tight
in (\ref{eq:11-4}). While a tight bound is hard to obtain for (\ref{eq:11-4-1}),
it turns out that the bounds $2^{k-1}\left(k-1\right)!$ from Lemma
3.2 only make $W_{1}(\delta,\gamma)$ in Theorem 3.2 differ from $Z_{1}(\delta,\gamma)$
in Theorem 3.1 by an extra term $\frac{\gamma^{2}+\gamma}{2}\log2$.
For any $\gamma\geq5$, $\log\left(\prod_{i=0}^{\gamma}i!\right)\geq\frac{\gamma^{2}+\gamma}{2}\log2$
and therefore the factorial part ($\prod_{i=0}^{\gamma}i!$) dominates
the exponential part ($2^{\frac{\gamma^{2}+\gamma}{2}}$). Furthermore,
note that $\log\left(\prod_{i=0}^{\gamma}i!\right)\leq\gamma^{2}\log\gamma$;
$\frac{\left(\gamma+2\right)\left(\gamma+3\right)}{6}\log\frac{5}{\delta}$
in $W_{2}(\delta,\gamma)$ is on the order of $\gamma^{2}\log\frac{1}{\delta}$,
which dominates $\gamma^{2}\log\gamma$ if $\gamma$ is small relative
to $\frac{1}{\delta}$ (for all practical purposes, this is always
the case).

For large enough $\gamma\in\left\{ 0,...,\beta\right\} $ and a range
of $\delta$, terms like $\log\left(\prod_{i=0}^{\gamma}i!\right)$
can dominate $\left(\frac{\delta}{5}\right)^{\frac{-1}{\gamma+2}}$
in $Z_{1}(\delta,\gamma)$ and $W_{1}(\delta,\gamma)$, and terms
like $\frac{\left(\gamma+2\right)\left(\gamma+3\right)}{6}\log\frac{5}{\delta}$
and even $\frac{\gamma+3}{2}\log\frac{5}{\delta}$ can dominate $\left(\frac{\delta}{5}\right)^{\frac{-2}{\gamma+1}}$
and $\left(\frac{\delta}{5}\right)^{\frac{-1}{\gamma+1}}$ in $W_{2}(\delta,\gamma)$
and $Z_{2}(\delta,\gamma)$. Because of this, we choose $\gamma$s
that minimize the $Z$ functions in Theorem 3.1, and the $W$ functions
in Theorem 3.2, using the fact that $\mathcal{AS}_{\beta+2,\overline{C}}^{\dagger}\subseteq\mathcal{AS}_{\beta+1,\overline{C}}^{\dagger}\subseteq\mathcal{AS}_{\beta,\overline{C}}^{\dagger}\cdots$,
$\mathcal{S}_{\beta+2,\overline{C}}^{\dagger}\subseteq\mathcal{S}_{\beta+1,\overline{C}}^{\dagger}\subseteq\mathcal{S}_{\beta,\overline{C}}^{\dagger}\cdots$,
etc. These previously mentioned log terms would clearly have implications
on the rates of convergence of least squares (\ref{eq:least squares})
and sample size requirement (as we will see in Section 3.2.2).

\subsubsection*{The best lower bound?}

We can easily obtain a lower bound on the covering numbers of solution
classes associated with the nonautonomous ODEs by restricting (\ref{eq:11-4-1})
to a class of separable ODEs, and then apply the construction in \cite{Kolmogorov and Tikhomirov 1959}
which gives (\ref{eq:Kolmogorov_lower}). This result is stated formally
below.\\
\textbf{}\\
\textbf{Lemma 3.3}\textit{ }(lower bound based on separable ODEs).\textit{
Let us consider a special case of (\ref{eq:11-4-1}):}

\begin{equation}
\textrm{separable ODE: }y^{'}\left(x\right)=f\left(x\right),\qquad y\left(0\right)=y_{0}.\label{eq:separate}
\end{equation}
\textit{In terms of (\ref{eq:separate}), suppose $f\in\mathcal{S}_{\beta+1}$.
Let $\mathcal{Y}^{sep}$ be the class of solutions (to (\ref{eq:separate})
with $f\in\mathcal{S}_{\beta+1}$) on $\left[0,\,1\right]$ and $\mathcal{Y}_{1}^{sep}$
be the class consisting of the first derivative $y^{'}$ of $y\in\mathcal{Y}_{sep}$.
We simply have 
\begin{eqnarray}
\log N_{\infty}\left(\delta,\,\mathcal{Y}^{sep}\right) & \succsim & \delta^{\frac{-1}{\beta+2}},\label{eq:14-4-1}\\
\log N_{\infty}\left(\delta,\,\mathcal{Y}_{1}^{sep}\right) & \succsim & \delta^{\frac{-1}{\beta+1}}.\label{eq:14-4-1b}
\end{eqnarray}
}

In the case of (\ref{eq:separate}), $\mathcal{Y}_{1}^{sep}=\mathcal{S}_{\beta+1}$,
and $\mathcal{Y}^{sep}$ and $\mathcal{S}_{\beta+2}$ have similar
size (with differences only in the constants). Let $\mathcal{Y}$
be the class of solutions (to (\ref{eq:11-4-1}) with $f\in\mathcal{S}_{\beta+1,\,2}\left(1,\,\Xi\right)$).
Since $\mathcal{Y}^{sep}\subseteq\mathcal{Y}$, we must have 
\[
\log N_{\infty}\left(\delta,\,\mathcal{Y}\right)\geq\log N_{\infty}\left(\delta,\,\mathcal{Y}^{sep}\right)\succsim\delta^{\frac{-1}{\beta+2}}.
\]
Similarly, letting $\mathcal{Y}_{1}$ be the class consisting of the
first derivative $y^{'}$ of $y\in\mathcal{Y}$, we must have 
\[
\log N_{\infty}\left(\delta,\,\mathcal{Y}_{1}\right)\geq\log N_{\infty}\left(\delta,\,\mathcal{Y}_{1}^{sep}\right)\succsim\delta^{\frac{-1}{\beta+1}}.
\]

The questions remain, whether, one can find a solution class whose
size attains (\ref{eq:14-4-2}). Before answering this question, a
more fundamental question is whether one can find a smooth class whose
size attains (\ref{eq:Kolmogorov_upper}). We conjecture that the
lower bound (\ref{eq:Kolmogorov_lower}) can be sharpened based on
the construction in \cite{Kolmogorov and Tikhomirov 1959}. The problems
only get more challenging in the ODE context and we hope to seek insights
from the community. 

\subsubsection{Implications of the previous results}

In what follows, we apply the results in Section 3.2.1 to study the
theoretical behavior of least squares in noisy function fitting problems.
\\
\textbf{}\\
\textbf{Problem setup}. Let us consider (\ref{eq:model}) where in
\textbf{problem (i):} $y\left(\cdot\right)$ is a solution to (\ref{eq:11-4})
and $f\in\mathcal{S}_{\beta+1,\,1}\left(1,\,\left[-\overline{C},\,\overline{C}\right]\right)$;
\textbf{in problem (ii):} $y\left(\cdot\right)$ is a solution to
(\ref{eq:11-4-1}) and $f\in\mathcal{S}_{\beta+1,\,2}\left(1,\,\Xi\right)$.
We seek rates of convergence for the mean squared errors of the least
squares

\begin{equation}
\hat{y}\in\arg\min_{\tilde{y}\in\mathcal{Y}}\frac{1}{2n}\sum_{i=1}^{n}\left(Y_{i}-\tilde{y}\left(x_{i}\right)\right)^{2}\label{eq:nonpara ls}
\end{equation}
where in problem (i), $\mathcal{Y}$ is the class of solutions (to
(\ref{eq:11-4}) with $f\in\mathcal{S}_{\beta+1,\,1}\left(1,\,\left[-\overline{C},\,\overline{C}\right]\right)$)
on $\left[0,\,\alpha\right]$; in problem (ii), $\mathcal{Y}$ is
the class of solutions (to (\ref{eq:11-4-1}) with $f\in\mathcal{S}_{\beta+1,\,2}\left(1,\,\Xi\right)$)
on $\left[0,\,\alpha\right]$. \\
\textbf{}\\
\textbf{Definition}. Let us define the following functions
\begin{eqnarray*}
\mathcal{A}_{1}\left(\tilde{r}_{n}\right) & = & c_{0}\min_{\gamma\in\left\{ 0,...,\beta\right\} }\left(\frac{\tilde{r}_{n}}{\sqrt{n}}\sqrt{1\vee\log\left(\prod_{i=0}^{\gamma}i!\right)}+\frac{1}{\sqrt{n}}\tilde{r}_{n}^{\frac{2\gamma+3}{2\gamma+4}}\right),\\
\mathcal{A}_{2}\left(\tilde{r}_{n}\right) & = & c_{1}\min_{\gamma\in\left\{ 0,...,\beta\right\} }\left(\tilde{r}_{n}\sqrt{\frac{\gamma\vee1}{n}}+\frac{1}{\sqrt{n}}\tilde{r}_{n}^{\frac{2\gamma+1}{2\gamma+2}}\right),\\
\mathcal{B}_{1}\left(\tilde{r}_{n}\right) & = & c_{0}\min_{\gamma\in\left\{ 0,...,\beta\right\} }\left(\tilde{r}_{n}\sqrt{\frac{1}{n}\log\left(\prod_{i=0}^{\gamma}i!\right)}+\tilde{r}_{n}\sqrt{\frac{\gamma^{2}\vee1}{n}}+\frac{1}{\sqrt{n}}\tilde{r}_{n}^{\frac{2\gamma+3}{2\gamma+4}}\right),\\
\mathcal{B}_{2}\left(\tilde{r}_{n}\right) & = & c_{1}\min_{\gamma\in\left\{ 1,...,\beta\right\} }\left(\tilde{r}_{n}\sqrt{\frac{\gamma^{2}\vee1}{n}}+\frac{1}{\sqrt{n}}\tilde{r}_{n}^{\frac{\gamma}{\gamma+1}}\right),\qquad(\beta>0),
\end{eqnarray*}
and 
\begin{eqnarray*}
\mathcal{M}_{1}^{a}\left(\gamma\right) & = & \max\left\{ \frac{\sigma^{2}}{n}\log\left(\prod_{i=0}^{\gamma}i!\right),\,\frac{\sigma^{2}}{n},\,\left(\frac{\sigma^{2}}{n}\right)^{\frac{2\left(\gamma+2\right)}{2\left(\gamma+2\right)+1}}\right\} ,\\
\mathcal{M}_{2}^{a}\left(\gamma\right) & = & \max\left\{ \frac{\sigma^{2}\left(\gamma\vee1\right)}{n},\,\left(\frac{\sigma^{2}}{n}\right)^{\frac{2\left(\gamma+1\right)}{2\left(\gamma+1\right)+1}}\right\} ,\\
\mathcal{M}_{1}\left(\gamma\right) & = & \max\left\{ \frac{\sigma^{2}}{n}\log\left(\prod_{i=0}^{\gamma}i!\right),\,\frac{\sigma^{2}\left(\gamma^{2}\vee1\right)}{n},\,\left(\frac{\sigma^{2}}{n}\right)^{\frac{2\left(\gamma+2\right)}{2\left(\gamma+2\right)+1}}\right\} ,\\
\mathcal{M}_{2}\left(\gamma\right) & = & \max\left\{ \frac{\sigma^{2}\left(\gamma^{2}\vee1\right)}{n},\,\left(\frac{\sigma^{2}}{n}\right)^{\frac{\gamma+1}{\gamma+2}}\right\} ,
\end{eqnarray*}
and 
\begin{eqnarray}
\tilde{r}_{n,a}^{2} & = & \min\left\{ \mathcal{M}_{1}^{a}\left(\gamma_{1a}^{*}\right),\,\mathcal{M}_{2}^{a}\left(\gamma_{2a}^{*}\right)\right\} ,\quad\textrm{(autonomous ODE)}\label{eq:r_na}\\
\tilde{r}_{n}^{2} & = & \min\left\{ \mathcal{M}_{1}\left(\gamma_{1}^{*}\right),\,\mathcal{M}_{2}\left(\gamma_{2}^{*}\right)\right\} ,\quad\textrm{(nonautonomous ODE)}\label{eq:r_n}
\end{eqnarray}
where $\gamma_{1a}^{*}$, $\gamma_{2a}^{*}$, $\gamma_{1}^{*}$, and
$\gamma_{2}^{*}$ are the minimizers that give $\mathcal{A}_{1}\left(\tilde{r}_{n}\right)$,
$\mathcal{A}_{2}\left(\tilde{r}_{n}\right)$, $\mathcal{B}_{1}\left(\tilde{r}_{n}\right)$
and $\mathcal{B}_{2}\left(\tilde{r}_{n}\right)$, respectively. \\
\textbf{}\\
\textbf{Theorem 3.3}. \textit{Under the }\textbf{problem setup}\textit{, }

\textit{(i) if $\overline{C}\precsim1$ and $\max\left\{ \alpha,\,C_{0}\right\} \succsim\tilde{r}_{n,a}$,
then 
\[
\mathbb{E}\left[\frac{1}{n}\sum_{i=1}^{n}\left(\hat{y}(x_{i})-y(x_{i})\right)^{2}\right]\precsim\tilde{r}_{n,a}^{2}+\exp\left\{ -cn\sigma^{-2}\tilde{r}_{n,a}^{2}\right\} ;\quad\textrm{(autonomous ODE)}
\]
}

\textit{(ii) if $\overline{C}\precsim1$, }$\frac{\sigma^{2}}{n}\precsim1$,\footnote{In fact, $\left(\frac{\sigma^{2}}{n}\right)^{\frac{\gamma+1}{\gamma+2}}$
in $\mathcal{M}_{2}\left(\gamma\right)$ is derived for $\gamma>0$
in Section \ref{subsec:Theorem-3.3} for technical reasons. However,
since we are only interested in (\ref{eq:r_n}), under the condition
$\frac{\sigma^{2}}{n}\precsim1$, having $\gamma=0$ in $\mathcal{M}_{2}\left(\gamma\right)$
will not affect the result in part (ii). }\textit{ and $\max\left\{ \alpha,\,C_{0}\right\} \succsim\tilde{r}_{n}$,
then 
\[
\mathbb{E}\left[\frac{1}{n}\sum_{i=1}^{n}\left(\hat{y}(x_{i})-y(x_{i})\right)^{2}\right]\precsim\tilde{r}_{n}^{2}+\exp\left\{ -cn\sigma^{-2}\tilde{r}_{n}^{2}\right\} .\quad\textrm{(nonautonomous ODE)}
\]
}

\textit{(iii) If $\frac{n}{\sigma^{2}}\succsim\left(\gamma\sqrt{\log\left(\gamma\vee1\right)}\right)^{4\left(\gamma+2\right)+2}$
with $\gamma\in\left\{ 0,...,\beta\right\} $, for both autonomous
and nonautonomous ODEs, we have 
\begin{equation}
\mathbb{E}\left[\frac{1}{n}\sum_{i=1}^{n}\left(\hat{y}(x_{i})-y(x_{i})\right)^{2}\right]\precsim\left(\frac{\sigma^{2}}{n}\right)^{\frac{2\left(\gamma+2\right)}{2\left(\gamma+2\right)+1}}+\exp\left\{ -c\left(\frac{n}{\sigma^{2}}\right)^{\frac{1}{2\left(\gamma+2\right)+1}}\right\} .\label{eq:minimax}
\end{equation}
}

\textit{(iv) For the nonautonomous ODE, if $\frac{n}{\sigma^{2}}\succsim\left(\beta\sqrt{\log\left(\beta\vee1\right)}\right)^{4\left(\beta+2\right)+2}$,
(\ref{eq:minimax}) holds with $\gamma=\beta$ and the rate $\left(\frac{\sigma^{2}}{n}\right)^{\frac{2\left(\beta+2\right)}{2\left(\beta+2\right)+1}}$
is (minimax) optimal.}\\

The proof for Theorem 3.3 is provided in Section \ref{subsec:Theorem-3.3}
of the supplementary materials. 

The bounds in Theorem 3.3 exploit the covering number results in Theorems
3.1 and 3.2. Note that neither $\mathcal{M}_{1}^{a}$ nor $\mathcal{M}_{2}^{a}$
for the autonomous ODEs (respectively, neither $\mathcal{M}_{1}$
nor $\mathcal{M}_{2}$ for the nonautonomous ODEs) outperform each
other in all situations across different sample sizes ($n$) and degrees
of smoothness ($\beta$) for $f$ in the ODEs. These bounds for Theorem
3.3 may seem cumbersome but hard to simplify. \\
\\
\textbf{Remark}. The assumptions $\max\left\{ \alpha,\,C_{0}\right\} \succsim\tilde{r}_{n,a}$
and $\max\left\{ \alpha,\,C_{0}\right\} \succsim\tilde{r}_{n}$ in
Theorem 3.3 simply exclude the trivial case where $b$ and $C_{0}$
are ``too small''. Without these assumptions, as long as $f$ is
bounded from above, we would simply bound the convergence rate of
$\sqrt{\mathbb{E}\left[\frac{1}{n}\sum_{i=1}^{n}\left(\hat{y}(x_{i})-y(x_{i})\right)^{2}\right]}$
with $\max\left\{ \alpha,\,C_{0}\right\} $.\\

Recalling \textit{Question 3}, we now see from Theorem 3.3(iii) that,
the convergence rate for the mean squared error of $\hat{y}$ is roughly
bounded from above by $\left(\frac{\sigma^{2}}{n}\right)^{\frac{2\left(\beta+2\right)}{2\left(\beta+2\right)+1}}$
if $n\succsim\sigma^{2}\left(\beta\sqrt{\log\left(\beta\vee1\right)}\right)^{4\left(\beta+2\right)+2}$.
Note that this sample size requirement is derived from $\mathcal{M}_{1}^{a}\left(\beta\right)$
and $\mathcal{M}_{1}\left(\beta\right)$. Under such a requirement,
the rate $\left(\frac{\sigma^{2}}{n}\right)^{\frac{2\left(\beta+2\right)}{2\left(\beta+2\right)+1}}$
is (minimax) optimal for the nonautonomous ODEs. 

\subsubsection*{Practical fitting via kernel functions}

In Theorem 3.3, we have applied the covering number results from Theorems
3.1 and 3.2 (which are based on Theorem 2.1 and Lemmas 3.1--3.2)
to analyze the theoretical behavior of the least squares. Alternatively,
we can also exploit the factorial bounds in Lemmas 3.1--3.2 directly
by considering the following spaces (that contain $\mathcal{Y}$):

\begin{align}
\mathcal{AH}_{\beta+2}^{\dagger} & =\{h:\,\left[0,\,1\right]\rightarrow\mathbb{R}\vert h^{(k-1)}\textrm{ is abs. cont. and}\nonumber \\
 & \int_{0}^{1}\left[h^{(k)}\left(t\right)\right]^{2}dt\leq\left[\left(k-1\right)!\right]^{2}\;\forall k\in\left\{ 1,...,\beta+2\right\} \},\label{eq:RKHS1}\\
\mathcal{H}_{\beta+2}^{\dagger} & =\{h:\,\left[0,\,1\right]\rightarrow\mathbb{R}\vert h^{(k-1)}\textrm{ is abs. cont. and }\nonumber \\
 & \int_{0}^{1}\left[h^{(k)}\left(t\right)\right]^{2}dt\leq\left[2^{k-1}\left(k-1\right)!\right]^{2}\;\forall k\in\left\{ 1,...,\beta+2\right\} \}.\label{eq:RKHS2}
\end{align}
Importantly, working with the spaces (\ref{eq:RKHS1}) and (\ref{eq:RKHS2})
allows us to implement (\ref{eq:nonpara ls}) via kernel functions.
In particular, let us consider

\begin{equation}
\hat{y}\in\arg\min_{\tilde{y}\in\mathcal{E}}\frac{1}{2n}\sum_{i=1}^{n}\left(Y_{i}-\tilde{y}\left(x_{i}\right)\right)^{2}\label{eq:krr}
\end{equation}
where in problem (i), we let $\mathcal{E}=\mathcal{AH}_{\beta+2}^{\dagger}$;
in problem (ii), we let $\mathcal{E}=\mathcal{H}_{\beta+2}^{\dagger}$. 

For a given $k\in\left\{ 0,...,\beta+1\right\} $, let $\mathbb{K}_{k}\in\mathbb{R}^{n\times n}$
consist of entries in the following form 
\begin{eqnarray*}
\mathcal{K}_{k}\left(x_{i},\,x_{j}\right) & = & x_{i}\wedge x_{j},\quad k=0\\
\mathcal{K}_{k}\left(x_{i},\,x_{j}\right) & = & \int_{0}^{1}\frac{\left(x_{i}-t\right)_{+}^{k}}{k!}\frac{\left(x_{j}-t\right)_{+}^{k}}{k!}dt,\quad k>0
\end{eqnarray*}
with $\left(x\right)_{+}=x\vee0$. For problem (i), solving program
(\ref{eq:krr}) is equivalent to solving the following problem 
\begin{align}
\left(\hat{\alpha},\,\hat{\pi}\right)= & \arg\min_{\left(\alpha,\,\pi\right)\in\mathbb{R}\times\mathbb{R}^{n}}\frac{1}{2n}\left|Y-\alpha1_{n}-\sqrt{n}\mathbb{K}_{\beta+1}\pi\right|_{2}^{2}\nonumber \\
\textrm{s.t.} & \pi^{T}\mathbb{K}_{k}\pi\leq C\left(k!\right)^{2},\quad\forall k=0,...,\beta+1\label{eq:constraints}
\end{align}
(where $1_{n}$ is a column vector of $1$s) and form $\hat{y}\left(\cdot\right)=\hat{\alpha}+\frac{1}{\sqrt{n}}\sum_{i=1}^{n}\hat{\pi}_{i}\mathcal{K}_{\beta+1}\left(\cdot,\,x_{i}\right)$.
For problem (ii), we simply replace (\ref{eq:constraints}) by 
\begin{equation}
\pi^{T}\mathbb{K}_{k}\pi\leq C\left(2^{k}k!\right)^{2}\quad\forall k=0,...,\beta+1.\label{eq:constraints2}
\end{equation}
In this section, we may take $C=1$; more generally in practice, if
we have no prior knowledge on $C$, we can use cross validation to
select $C$.

The typical kernel fitting algorithm such as (\ref{eq:standard_krr})
in Sobolev space taking the form of (\ref{eq:sobolev}) has only one
constraint related to the squared RKHS norm. In comparison, (\ref{eq:constraints})
and (\ref{eq:constraints2}) exploit the unique structures of (\ref{eq:RKHS1})
and (\ref{eq:RKHS2}) as a result of Lemmas 3.1--3.2.\\
\textbf{}\\
\textbf{Theorem 3.4}. \textit{Assume $\overline{C}\precsim1$ and
$\max\left\{ \alpha,\,C_{0}\right\} \succsim\bar{r}_{n}$. Under the
}\textbf{problem setup}\textit{, in terms of (\ref{eq:krr}) and for
both autonomous and nonautonomous ODEs, we have}

\textit{
\begin{equation}
\mathbb{E}\left[\frac{1}{n}\sum_{i=1}^{n}\left(\hat{y}(x_{i})-y(x_{i})\right)^{2}\right]\precsim\bar{r}_{n}^{2}+\exp\left\{ -cn\sigma^{-2}\bar{r}_{n}^{2}\right\} \label{eq:minimax-1}
\end{equation}
where 
\begin{align*}
\bar{r}_{n}^{2} & =\min_{\gamma\in\left\{ 0,...,\beta\right\} }\max\left\{ \frac{\sigma^{2}\left(\left(\gamma\wedge n\right)\vee1\right)}{n},\,\left[2^{\gamma+1}\left(\gamma+1\right)!\right]^{\frac{2}{2\left(\gamma+2\right)+1}}\left(\frac{\sigma^{2}}{n}\right)^{\frac{2\left(\gamma+2\right)}{2\left(\gamma+2\right)+1}}\right\} \\
 & \precsim\min_{\gamma\in\left\{ 0,...,\beta\right\} }\max\left\{ \frac{\sigma^{2}\left(\left(\gamma\wedge n\right)\vee1\right)}{n},\,\left(\frac{\left(\gamma\vee1\right)\sigma^{2}}{n}\right)^{\frac{2\left(\gamma+2\right)}{2\left(\gamma+2\right)+1}}\right\} .
\end{align*}
}

The proof for Theorem 3.4 is provided in Section \ref{subsec:Theorem-3.4}
of the supplementary materials. 

In contrast with (\ref{eq:minimax}) in Theorem 3.3(iii), one theoretical
drawback of bound (\ref{eq:minimax-1}) is that it cannot achieve
the rate $\left(\frac{\sigma^{2}}{n}\right)^{\frac{2\left(\beta+2\right)}{2\left(\beta+2\right)+1}}$
even if $n\succsim\sigma^{2}\left(\beta\sqrt{\log\left(\beta\vee1\right)}\right)^{4\left(\beta+2\right)+2}$.
However, when $n$ is this large, the factor $\beta\vee1$ is much
less substantial than $\frac{\sigma^{2}}{n}$. For all practical purposes,
the guarantees for (\ref{eq:krr}) via kernel fitting are quite comparable
as those in Theorem 3.3(i)-(ii). Figure 2 exhibits the growth of $\mathcal{M}_{1}^{a}\left(\gamma\right)$,
$\mathcal{M}_{2}^{a}\left(\gamma\right)$, $\mathcal{M}_{1}\left(\gamma\right)$,
$\mathcal{M}_{2}\left(\gamma\right)$, and $\mathcal{M}_{3}\left(\gamma\right):=\left(\frac{\left(\gamma\vee1\right)\sigma^{2}}{n}\right)^{\frac{2\left(\gamma+2\right)}{2\left(\gamma+2\right)+1}}$
as $\gamma$ increases for $n=50$ and $n=50000$. These plots suggest
that it is hard for the minimizers $\gamma_{1a}^{*}$, $\gamma_{2a}^{*}$
in (\ref{eq:r_na}), $\gamma_{1}^{*}$, and $\gamma_{2}^{*}$ in (\ref{eq:r_n}),
and the minimizer yielding $\bar{r}_{n}^{2}$, to take large $\gamma$
values unless $n$ is huge. For small $\gamma$s, the factor $\gamma\vee1$
does not contribute much in $\left(\frac{\left(\gamma\vee1\right)\sigma^{2}}{n}\right)^{\frac{2\left(\gamma+2\right)}{2\left(\gamma+2\right)+1}}$.

\begin{figure}
\caption{{\small{}``$*$'': $\mathcal{M}_{1}^{a}\left(\gamma\right)$; ``$+$'':
$\mathcal{M}_{2}^{a}\left(\gamma\right)$; ``$\ocircle$'': $\mathcal{M}_{1}\left(\gamma\right)$;
``$\boxempty$'': $\mathcal{M}_{2}\left(\gamma\right)$; ``$\times$'':
$\mathcal{M}_{3}\left(\gamma\right)$}}

\includegraphics[scale=0.5]{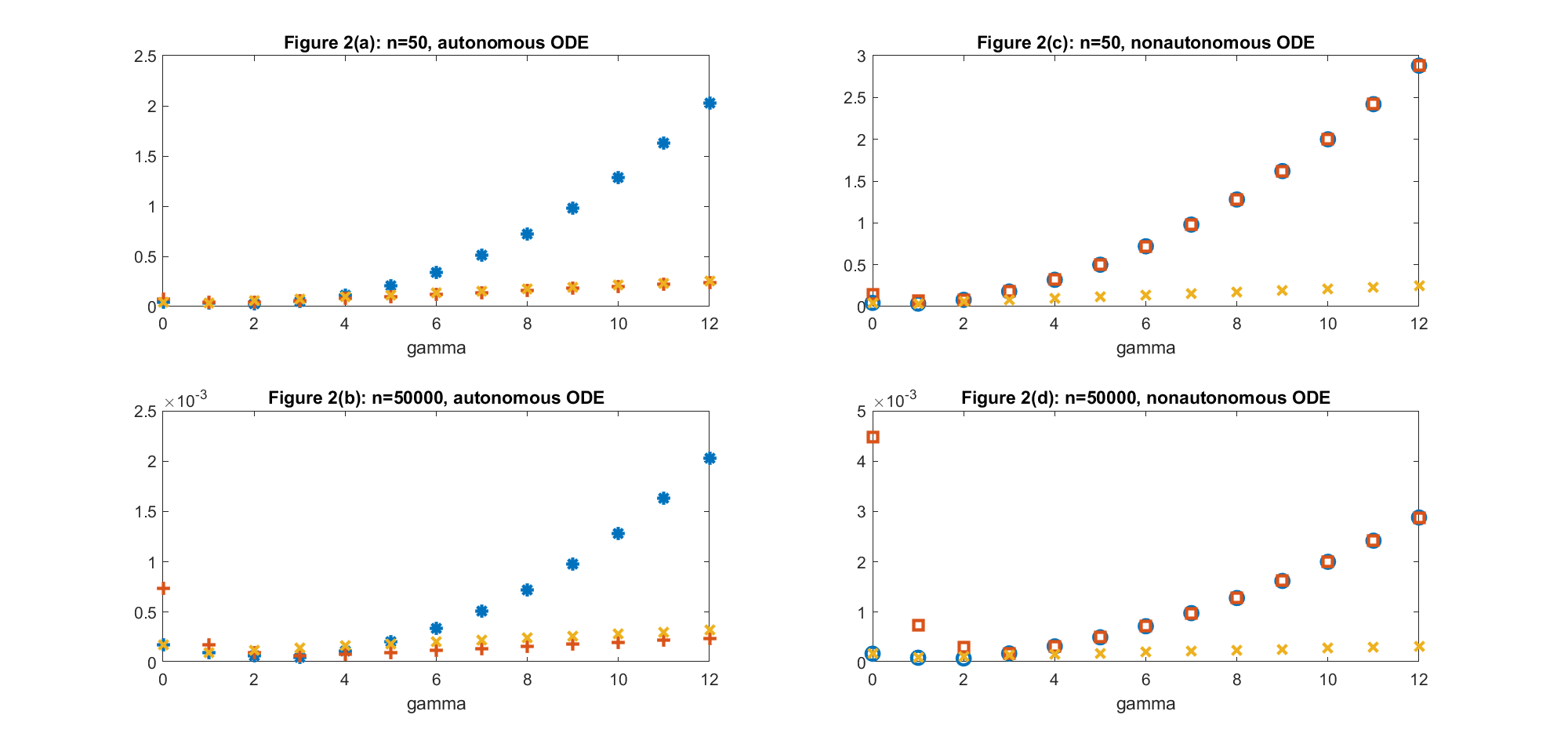}
\end{figure}

\section{Conclusion}

Numerical methods for recovering ODE solutions from data largely rely
on approximating the solutions with basis functions or kernel functions
under a least squares criterion. This strategy assumes some ambient
space of smooth functions that contains the solutions and seeks a
close enough estimator from the restricted class of functions. The
accuracy of these methods hinges on the smoothness of the solutions.
In this paper, we provide a theoretical foundation for these methods
by establishing novel results on the smoothness and covering numbers
of ODE solution classes (as a measure of their ``size'') as well
as the rates of convergence for least squares fitting in noisy settings. 

We choose to focus on ODEs in this paper as this is a natural and
necessary first step before one delves into more complicated dynamic
systems. A potential but perhaps challenging extension would be to
explore the ``size'' of solution classes associated with partial
differential equations. As for other extensions, because of the deep
connections between ODEs and contraction mapping, it may be worthwhile
to extend the analyses in this paper to study the ``size'' of fixed
point classes (which play an important role in reinforcement learning). 

\textbf{Note: The supplement of this paper has its own list of references
at the end. }

\newpage{}

\appendix

\part*{Supplementary materials for ``Classes of ODE solutions: smoothness,
covering numbers, implications for noisy function fitting, and the
curse of smoothness phenomenon''}

by Ying Zhu\footnote{Assistant Professor, Department of Economics, UC San Diego. Corresponding
author. yiz012@ucsd.edu.} and Mozhgan Mirzaei\footnote{PhD, former student at Department of Mathematics, UC San Diego. momirzae@ucsd.edu}

\section{Proofs for the main results}

\subsection{Theorem 2.1\label{subsec:Theorem-2.1}}

\textbf{Proof}. Let $L_{\max}=\sup_{x\in\left[a_{0},\,a_{0}+\alpha\right]}\left\{ \exp\left(x\sqrt{L^{2}+1}\right)\left[1+\int_{a_{0}}^{x}\exp\left(-s\sqrt{L^{2}+1}\right)ds\right]\right\} $.
For a given $\delta>0$, we consider the smallest $\frac{\delta}{L_{\max}}-$covering
of $\mathcal{F}$ with respect to the sup-norm. We also consider the
smallest $\frac{\delta}{L_{\max}}-$covering of $\mathbb{B}_{2}\left(C_{0}\right):=\left\{ \theta\in\mathbb{R}^{m}:\,\left|\theta\right|_{2}\leq C_{0}\right\} $
(where the initial values lie) with respect to the $l_{2}-$norm.
Note that the standard volumetric argument yields
\[
\log N_{2}\left(\frac{\delta}{L_{\max}},\,\mathbb{B}_{2}\left(C_{0}\right)\right)\leq m\log\left(\frac{2C_{0}L_{\max}}{\delta}+1\right).
\]
By Theorem B.2 in Section \ref{sec:Supporting-lemmas-and}, for any
$y\in\mathcal{Y}$ with $Y_{0}\in\mathbb{B}_{2}\left(C_{0}\right)$,
we can find an element (indexed by $i$) from the smallest $\frac{\delta}{L_{\max}}-$covering
of $\mathcal{F}$ and an element (indexed by $i^{'}$) from the smallest
$\frac{\delta}{L_{\max}}-$covering of $\mathbb{B}_{2}\left(C_{0}\right)$
such that 
\[
\left|y^{(k)}\left(x\right)-y_{(i,i^{'})}^{(k)}\left(x\right)\right|\leq\frac{\delta}{L_{\max}}\exp\left(x\sqrt{L^{2}+1}\right)\left[1+\int_{0}^{x}\exp\left(-s\sqrt{L^{2}+1}\right)ds\right]\leq\delta\qquad\forall x\in\left[a_{0},\,a_{0}+\alpha\right]
\]
for all $k=0,...,m-1$, where $y_{(i,i^{'})}$ is a solution to the
ODE associated with $f_{i}$ and the initial value $Y_{0,\,i^{'}}$
from the covering sets, and $y_{(i,i^{'})}^{(k)}$ is the $k$th derivative
of $y_{(i,i^{'})}$ ($k\leq m-1$). Consequently, we obtain a $\delta-$cover
of $\mathcal{Y}_{k}$. We conclude that
\[
\log N_{\infty}\left(\delta,\,\mathcal{Y}_{k}\right)\leq\log N_{\infty}\left(\frac{\delta}{L_{\max}},\,\mathcal{F}\right)+m\log\left(\frac{2C_{0}L_{\max}}{\delta}+1\right).
\]

$\square$

\subsection{Corollary 2.1\label{subsec:Corollary-2.1}}

\textbf{Proof}. Let $N_{q}\left(\delta,\,\mathbb{B}_{q}\left(1\right)\right)$
denote the covering number of $\mathbb{B}_{q}\left(1\right)$ with
respect to the $l_{q}-$norm and $N_{2}\left(\frac{\delta}{L_{\max}},\,\mathbb{B}_{2}\left(C_{0}\right)\right)$
denote the covering number of $\mathbb{B}_{2}\left(C_{0}\right)$
with respect to the $l_{2}-$norm. For a given $\delta>0$, let us
consider the smallest $\frac{\delta}{L_{\max}L_{K}}-$covering $\left\{ \theta^{1},...,\theta^{N}\right\} $
with respect to the $l_{q}-$norm. By (\ref{eq:17}), note that $\left\{ f_{\theta^{1}},\,f_{\theta^{2}},...,\,f_{\theta^{N}}\right\} $
forms a $\frac{\delta}{L_{\max}}-$cover of $\mathcal{F}$ with respect
to the sup-norm. We also consider the smallest $\frac{\delta}{L_{\max}}-$covering
$\left\{ Y_{0,1},...,Y_{0,N^{'}}\right\} $ for the interval $\mathbb{B}_{2}\left(C_{0}\right)$
where the initial values lie. By Theorem B.2 in Section \ref{sec:Supporting-lemmas-and},
for a solution $y$ to the ODE with $f$ parameterized by any $\theta\in\mathbb{B}_{q}\left(1\right)$
(and subject to (\ref{eq:17})) and $Y_{0}\in\mathbb{B}_{2}\left(C_{0}\right)$,
we can find $i\in\left\{ 1,...,N\right\} $ and $i^{'}\in\left\{ 1,...,N^{'}\right\} $
such that
\[
\left|y_{\left(i,i^{'}\right)}\left(x\right)-y\left(x\right)\right|\leq\frac{\delta}{C_{\max}}\exp\left(x\right)\left[1+\int_{0}^{x}\exp\left(-s\right)ds\right]\leq\delta\qquad\forall x\in\left[0,\,\alpha\right]
\]
where $y_{\left(i,i^{'}\right)}$ is a solution to the ODE with $f$
parameterized by $\theta^{i}$ and the initial value being $y_{0,i^{'}}$.
Consequently, we obtain a $\delta-$cover of $\mathcal{Y}$. By the
standard volumetric argument which yields
\[
\log N_{2}\left(\delta,\,\mathbb{B}_{q}\left(1\right)\right)\leq K\log\left(1+\frac{2L_{\max}L_{K}}{\delta}\right),
\]
and 
\[
\log N_{2}\left(\frac{\delta}{L_{\max}},\,\mathbb{B}_{2}\left(C_{0}\right)\right)\leq\log\left(\frac{2C_{0}L_{\max}}{\delta}+1\right),
\]
we conclude that 
\[
\log N_{\infty}\left(\delta,\,\mathcal{Y}\right)\leq K\log\left(1+\frac{2L_{\max}L_{K}}{\delta}\right)+\log\left(\frac{2C_{0}L_{\max}}{\delta}+1\right).
\]

$\square$

\subsection{Lemmas 3.1 and 3.2\label{subsec:Lemmas 3.1 and 3.2}}

In part (i) of Lemma 3.1, by the mean value theorem, (\ref{eq:a2})
with $k=0$ (i.e., $\left|f^{(0)}(y(x))\right|=\left|f(y(x))\right|=\left|y^{'}(x)\right|\leq1$)
and (\ref{eq:a3}) imply that 
\begin{equation}
\left|f^{\beta}(y(x))-f^{\beta}(y(x^{'}))\right|\leq|y(x)-y(x^{'})|\leq|x-x^{'}|;\label{eq:a3-2}
\end{equation}
moreover, (\ref{eq:a2}) with $1\leq k\leq\beta$ implies that 
\begin{equation}
\left|f^{k-1}(y(x))-f^{k-1}(y(x^{'}))\right|\leq|y(x)-y(x^{'})|\leq|x-x^{'}|,\quad\forall\,1\leq k\leq\beta.\label{eq:a2-2}
\end{equation}
In Lemma 3.2, by the mean value theorem, the assumption that $\left|f(x,\,y(x))\right|=\left|y^{'}(x)\right|\leq1$
and (\ref{eq:a7-1}) imply that 
\begin{equation}
\left|D^{p}f\left(x,\,y\left(x\right)\right)-D^{p}f\left(x^{'},\,y\left(x^{'}\right)\right)\right|\leq\max\left\{ \left|x-x^{'}\right|,\,\left|y\left(x\right)-y\left(x^{'}\right)\right|\right\} \leq\left|x-x^{'}\right|\label{eq:a3-2-1}
\end{equation}
for all $p$ with $\left[p\right]=\beta$; moreover, the assumption
that $\left|D^{p}f\left(x,\,y\right)\right|\leq1$ for all $p$ with
$1\leq\left[p\right]\leq\beta$ implies that 
\begin{equation}
\left|D^{p-1}f\left(x,\,y\left(x\right)\right)-D^{p-1}f\left(x^{'},\,y\left(x^{'}\right)\right)\right|\leq\left|x-x^{'}\right|+\left|y\left(x\right)-y\left(x^{'}\right)\right|\leq2\left|x-x^{'}\right|,\label{eq:a2-2-1}
\end{equation}
for all $p$ with $1\leq\left[p\right]\leq\beta$.\textbf{}\\
\\
\textbf{Remark}. The following proofs can easily handle situations
where the bound ``$1$'' on the absolute values of the derivatives
and the Lipschitz conditions for $f$ in Lemmas 3.1 and 3.2 is replaced
with general constants. \\
\\
We first prove parts (i)-(ii) and then part (iii). Let us begin with
some intuitions for the autonomous ODEs. When $\beta=0$, $\left|y^{'}(x)\right|=|f(y(x))|\leq1$
for all $x$ and $y(x)$ on $\left[x_{0}-a,\,x_{0}+a\right]\times\left[y_{0}-b,\,y_{0}+b\right]$;
moreover, we have

\[
\left|y^{'}(x)-y^{'}\left(x^{'}\right)\right|=\left|f(y(x))-f\left(y\left(x^{'}\right)\right)\right|\leq\left|x-x^{'}\right|.
\]
When $\beta=1$, note that 
\[
\left|y^{(2)}(x)\right|=\left|\frac{\partial y^{(1)}(x)}{\partial x}\right|=\left|f^{(1)}(y(x))y^{(1)}(x)\right|\leq1
\]
and 
\[
\begin{aligned}\left|y^{(2)}(x)-y^{(2)}\left(x^{'}\right)\right|= & \left|f^{(1)}(y(x))y^{(1)}(x)-f^{(1)}\left(y\left(x^{'}\right)\right)y^{(1)}\left(x^{'}\right)\right|\\
\leq & \left|f^{(1)}(y(x))y^{(1)}(x)-f^{(1)}\left(y\left(x^{'}\right)\right)y^{(1)}(x)\right|\\
 & +\left|f^{(1)}\left(y\left(x^{'}\right)\right)y^{(1)}(x)-f^{(1)}\left(y\left(x^{'}\right)\right)y^{(1)}\left(x^{'}\right)\right|\\
\leq & \left|y(x)-y\left(x^{'}\right)\right|+\left|y^{(1)}(x)-y^{(1)}\left(x^{'}\right)\right|\\
\leq & 2\left|x-x^{'}\right|.
\end{aligned}
\]
When $\beta=2$, we have 
\[
\begin{aligned}\left|y^{(3)}(x)\right| & =\left|f^{(2)}(y(x))\left(y^{(1)}(x)\right)^{2}+f^{(1)}(y(x))y^{(2)}(x)\right|\\
 & \leq\left|f^{(2)}(y(x))\right|+\left|f^{(1)}(y(x))\right|\leq2
\end{aligned}
\]
and

\[
\left|y^{(3)}(x)-y^{(3)}\left(x^{'}\right)\right|\leq6\left|x-x^{'}\right|.
\]
After trying $\beta=1,...,3$, we observe the following pattern: 

\begin{align}
 & y^{(1)}(x)=f(y(x))\nonumber \\
 & -------->(0)\nonumber \\
 & y^{(2)}(x)=f^{(1)}(y(x)).y^{(1)}(x)\nonumber \\
 & ~~~~~~~~=f^{(1)}(y(x)).f(y(x))\nonumber \\
 & -------->(1,0)\nonumber \\
 & y^{(3)}(x)=f^{(2)}(y(x)).\left(f(y(x))\right)^{2}+\left(f^{(1)}(y(x))\right)^{2}.f(y(x))\label{eq:69}\\
 & --------->(2,0,0),(1,1,0)\nonumber \\
 & y^{(4)}(x)=f^{(3)}(y(x)).\left(f(y(x))\right)^{3}+2f^{(2)}(y(x)).f^{(1)}(y(x)).\left(f(y(x))\right)^{2}+\nonumber \\
 & ~~~~~~~~~~~~2f^{(2)}(y(x)).f^{(1)}(y(x)).\left(f(y(x))\right)^{2}+\left(f^{(1)}(y(x))\right)^{3}.f(y(x))\nonumber \\
 & ---------->(3,0,0,0),2(2,1,0,0),2(2,1,0,0),(1,1,1,0)\nonumber 
\end{align}
In what follows, we derive $b_{k}$s (for all $k=1,...,\beta+1$)
such that 
\[
\left|y^{(k)}(x)\right|\leq b_{k},
\]
and $L_{\beta+1}$ such that 
\[
\left|y^{(\beta+1)}(x)-y^{(\beta+1)}\left(x^{'}\right)\right|\leq L_{\beta+1}\left|x-x^{'}\right|.
\]
\textbf{Proof of Lemma 3.1(i)}. For a $\beta-$times differentiable
function $f,$ define
\begin{align*}
f^{(a_{1},\cdots,a_{k})}(x):=f^{(a_{1})}(y(x)).f^{(a_{2})}(y(x)).\cdots.f^{(a_{k})}(y(x)),
\end{align*}
where $a_{1}+a_{2}+\cdots+a_{k}=k-1$ for all $k=1,...,\beta+1$,
and $a_{i}\geq0$ are all integers. We first show that

\begin{align}
\frac{d}{dx}f^{(a_{1},\cdots,a_{k})}(x)=\sum_{j=1}^{k}f^{(a_{1},\ldots,a_{j-1},a_{j}+1,0,a_{j+1},\ldots,a_{k})}(x).\label{eq:a1}
\end{align}
The equality (\ref{eq:a1}) follows from the derivations below: 
\begin{align*}
 & \frac{d}{dx}(f^{(a_{1})}(y(x)).f^{(a_{2})}(y(x)).\cdots.f^{(a_{k})}(y(x)))\\
 & =\sum_{j=1}^{k}f^{(a_{1})}(y(x)).\ldots f^{(a_{j-1})}(y(x)).\frac{d}{dx}(f^{(a_{j})}(y(x))).f^{(a_{j+1})}(y(x))\ldots f^{a_{k}}(y(x))\\
 & \,\,\,\,\,\,\,\,\,\,\,\,\,\,\,\,\,\,\,\,\,\,\,\,\,\,\,\,\,\,\,\,\,\,\,\,\,\,\,\,\,\,\,\,\,\,\,\,\,\,\,\,\,\,\,\,\,\,\,\,\,\,\,\,\,\,\,\,\,\,\,\,\,\,\,\,\,\,\,\,\,\,\,\,\,\,\,\,\,\,\,\,\,\,\,\,\,\,\,\,\,\,\,\,\,\,\,\,\,\,\,\,\,\,\,\,\,\,\,\,\,\,\,\,\,\,\,\,\,\,\,\,\,\,\,\,\,\,\,\,\,\,\,\,\,\,\,\,\,\,\,\,\,\,\,\,\,\,\,\textrm{by product rule}\\
 & =\sum_{j=1}^{k}f^{(a_{1})}(y(x)).\ldots\left[f^{(a_{j}+1)}(y(x)).y^{'}(x)\right]\ldots f^{a_{k}}(y(x))\\
 & \,\,\,\,\,\,\,\,\,\,\,\,\,\,\,\,\,\,\,\,\,\,\,\,\,\,\,\,\,\,\,\,\,\,\,\,\,\,\,\,\,\,\,\,\,\,\,\,\,\,\,\,\,\,\,\,\,\,\,\,\,\,\,\,\,\,\,\,\,\,\,\,\,\,\,\,\,\,\,\,\,\,\,\,\,\,\,\,\,\,\,\,\,\,\,\,\,\,\,\,\,\,\,\,\,\,\,\,\,\,\,\,\,\,\,\,\,\,\,\,\,\,\,\,\,\,\,\,\,\,\,\,\,\,\,\,\,\,\,\,\,\,\,\,\,\,\,\,\,\,\,\,\,\,\,\,\,\,\,\textrm{by chain rule}\\
 & =\sum_{j=1}^{k}f^{(a_{1})}(y(x)).\ldots\left[f^{(a_{j}+1)}(y(x)).f(y(x))\right]\ldots f^{a_{k}}(y(x)).
\end{align*}

\noindent Now by induction on $k,$ we show that if $f$ is $\beta-$times
differentiable, then for each $1\leq k\leq\beta+1,$ we have

\begin{align}
y^{(k)}(x)=\sum_{i=1}^{(k-1)!}f^{(a_{1}^{i},\cdots,a_{k}^{i})}(x),\hspace{1cm}\textrm{and for all }i\textrm{ in }(a_{1}^{i},\cdots,a_{k}^{i}),\,a_{1}^{i}+\ldots+a_{k}^{i}=k-1.\label{eq:a4}
\end{align}
For the base case, $k=1\Rightarrow y^{'}(x)=f(y(x))=f^{(a_{1})}(x)=\sum_{i=1}^{0!}f^{(a_{1})}(x)$
where $a_{1}=0$. Now assume $k\leq\beta+1$ and the induction hypothesis
holds for $k-1.$ Then

\begin{align*}
y^{(k)}(x) & =\frac{d}{dx}\left(y^{(k-1)}(x)\right)\\
 & =\frac{d}{dx}\left(\sum_{i=1}^{(k-2)!}f^{(a_{1}^{i},\ldots,a_{k-1}^{i})}(x)\right)\textrm{ where }\forall i,\,a_{1}^{i}+\ldots+a_{k-1}^{i}=k-2,\\
 & =\sum_{i=1}^{(k-2)!}\sum_{j=1}^{k-1}f^{\left(a_{1}^{i},\ldots,a_{j-1}^{i},a_{j}+1,0,a_{j+1}^{i},\ldots,a_{k-1}^{i}\right)}(x)~~~~\textrm{ by (\ref{eq:a1})}.
\end{align*}
Notice that $\left(a_{1}^{i},\ldots,a_{j-1}^{i},a_{j}^{i}+1,0,a_{j+1}^{i},\ldots,a_{k-1}^{i}\right)$
has exactly $k$ terms, and adds up to $a_{1}^{i}+\ldots+a_{k-1}^{i}+1+0=k-2+1+0=k-1,$
and in total there are $(k-2)!(k-1)=(k-1)!$ terms. This completes
the induction. By (\ref{eq:a2}), we have $\left|y^{(k)}(t)\right|\leq\left(k-1\right)!$.

To show the second claim $\left|y^{(\beta+1)}(x)-y^{(\beta+1)}(x^{'})\right|\leq\left(\beta+1\right)!|x-x^{'}|$,
we use (\ref{eq:a4}). Note that 
\begin{eqnarray*}
\left|y^{(\beta+1)}(x)-y^{(\beta+1)}(x^{'})\right| & = & \left|\sum_{i=1}^{\beta!}f^{\left(a_{1}^{i},\ldots,a_{\beta+1}^{i}\right)}(x)-\sum_{i=1}^{\beta!}f^{\left(a_{1}^{i},\ldots,a_{\beta+1}^{i}\right)}(x^{'})\right|\\
 & \leq & \sum_{i=1}^{\beta!}\left|f^{\left(a_{1}^{i},\ldots,a_{\beta+1}^{i}\right)}(x)-f^{\left(a_{1}^{i},\ldots,a_{\beta+1}^{i}\right)}(x^{'})\right|,
\end{eqnarray*}
where
\begin{eqnarray*}
 &  & \left|f^{\left(a_{1}^{i},\ldots,a_{\beta+1}^{i}\right)}(x)-f^{\left(a_{1}^{i},\ldots,a_{\beta+1}^{i}\right)}(x^{'})\right|\\
 & = & \left|\sum_{j=1}^{\beta+1}f^{(a_{1})}(y(x)).\ldots f^{(a_{j-1})}(y(x)).\left(f^{(a_{j})}(y(x))-f^{(a_{j})}(y(x^{'}))\right).f^{(a_{j+1})}(y(x^{'}))\ldots f^{(a_{\beta+1})}(y(x^{'}))\right|\\
 & \leq & \sum_{j=1}^{\beta+1}\left|f^{(a_{1})}(y(x))\ldots f^{(a_{j-1})}(y(x))\right|\left|\left(f^{(a_{j})}(y(x))-f^{(a_{j})}(y(x^{'}))\right)\right|\left|f^{(a_{j+1})}(y(x^{'}))\ldots f^{(a_{\beta+1})}(y(x^{'}))\right|\\
 & \leq & \sum_{j=1}^{\beta+1}|x-x^{'}|
\end{eqnarray*}
and the third line in the above follows from (\ref{eq:a3-2}) and
(\ref{eq:a2-2}). Hence, 
\begin{eqnarray*}
\left|y^{(\beta+1)}(x)-y^{(\beta+1)}(x^{'})\right| & = & \sum_{i=1}^{\beta!}\left|f^{\left(a_{1}^{i},\ldots,a_{\beta+1}^{i}\right)}(x)-f^{\left(a_{1}^{i},\ldots,a_{\beta+1}^{i}\right)}(x^{'})\right|\\
 & \leq & \beta!\left(\beta+1\right)|x-x^{'}|=\left(\beta+1\right)!|x-x^{'}|.
\end{eqnarray*}

$\square$\\
\textbf{}\\
\textbf{Proof of Lemma 3.1(ii)}. In view of (\ref{eq:a4}), if we
can find an infinitely differentiable function $f$ on $\mathbb{R}$
such that $\left|f^{(k)}(y(x))\right|=1$ for some $y\left(x\right)$
and all $k$, then $y^{(k)}(x)=(k-1)!$ for all $k=1,...,\beta+2$.
Note that $y^{'}=e^{-y-\frac{1}{2}}$ is such an example. For $y\in\left[-\frac{1}{2},\,\frac{1}{2}\right]$,
$\left|\frac{d^{k}\left(e^{-y-\frac{1}{2}}\right)}{dy^{k}}\right|\leq1$
for all $k$, satisfying (\ref{eq:a2}) and (\ref{eq:a3}). To argue
formally, let us consider $f$ with $f^{2j+1}\left(y\left(x\right)\right)=-1$
and $f^{2j}\left(y\left(x\right)\right)=1$ for all $j=0,1,...$ and
some $y\left(x\right)$. We can verify in (\ref{eq:69}) that $y^{(1)}(x)=0!$,
$y^{(2)}(x)=-1!$, $y^{(3)}(x)=2!$, and $y^{(4)}(x)=-3!$. Now performing
induction on $k,$we show that if $y^{(2j+1)}\left(x\right)=\left(2j\right)!$
($y^{(2j+2)}\left(x\right)=-\left(2j+1\right)!$), then $y^{(2j+2)}\left(x\right)=-\left(2j+1\right)!$
(respectively, $y^{(2j+3)}\left(x\right)=\left(2j+2\right)!$) for
all $j=0,1,...$. Note that whenever $y^{(2j+1)}\left(x\right)=\left(2j\right)!$,
each positive summand ``$1$'' in $y^{(2j+1)}(x)=\sum_{i=1}^{(2j)!}f^{(a_{1}^{i},\cdots,a_{2j+1}^{i})}(x)$
must be due to either any number of even derivatives of $f$ or an
even number of odd derivatives of $f$. Taking the derivative of any
even derivative of $f$ (the derivative of any odd derivative of $f$),
evaluated at $y\left(x\right)$, gives ``$-1$'' (respectively,
``$1$''), by the construction of $f$. This operation clearly makes
each summand in $\sum_{i=1}^{(2j+1)!}f^{(a_{1}^{i},\cdots,a_{2j+2}^{i})}(x)=y^{(2j+2)}(x)$
equal to ``$-1$'' and therefore $y^{(2j+2)}\left(x\right)=-\left(2j+1\right)!$.
The argument for $y^{(2j+3)}\left(x\right)=\left(2j+2\right)!$ from
$y^{(2j+2)}\left(x\right)=-\left(2j+1\right)!$ is very similar. $\square$\\

Now we show Lemma 3.2. As before, let us gain some intuitions first.
For $\beta=0$, $\left|y^{'}(x)\right|=|f(x,\,y(x))|\leq1$ for all
$x$ and $y(x)$ on $\left[x_{0}-a,\,x_{0}+a\right]\times\left[y_{0}-b,\,y_{0}+b\right]$;
moreover, we have

\begin{align*}
\left|y^{'}(x)-y^{'}\left(x^{'}\right)\right| & =\left|f(x,y(x))-f\left(x^{'},y\left(x^{'}\right)\right)\right|\\
 & \leq\left|f(x,y(x))-f\left(x^{'},y(x)\right)\right|+\left|f\left(x^{'},y(x)\right)-f\left(x^{'},y\left(x^{'}\right)\right)\right|\\
 & \leq2\left|x-x^{'}\right|.
\end{align*}
For $\beta=1$, note that
\[
\begin{aligned}\left|y^{(2)}(x)\right| & =\left|\frac{\partial y^{(1)}(x)}{\partial x}\right|=\left|f_{y}^{(1)}(x,y(x))y^{(1)}(x)+f_{x}^{(1)}(x,y(x))\right|\\
 & \leq\left|f_{y}^{(1)}(x,y(x))y^{(1)}(x)\right|+\left|f_{x}^{(1)}(x,y(x))\right|\leq2
\end{aligned}
\]
and 
\[
\begin{aligned}\left|y^{(2)}(x)-y^{(2)}\left(x^{'}\right)\right|\leq & \left|f_{y}^{(1)}(x,y(x))y^{(1)}(x)-f_{y}^{(1)}\left(x^{'},y\left(x^{'}\right)\right)y^{(1)}\left(x^{'}\right)\right|\\
 & +\left|f_{x}^{(1)}(x,y(x))-f_{x}^{(1)}\left(x^{'},y\left(x^{'}\right)\right)\right|\\
\leq & \left|f_{y}^{(1)}(x,y(x))y^{(1)}(x)-f_{y}^{(1)}\left(x^{'},y\left(x^{'}\right)\right)y^{(1)}(x)\right|\\
 & +\left|f_{y}^{(1)}\left(x^{'},y\left(x^{'}\right)\right)y^{(1)}(x)-f_{y}^{(1)}\left(x^{'},y\left(x^{'}\right)\right)y^{(1)}\left(x^{'}\right)\right|\\
 & +\left|f_{x}^{(1)}(x,y(x))-f_{x}^{(1)}\left(x^{'},y\left(x^{'}\right)\right)\right|\\
\leq & 6\left|x-x^{'}\right|.
\end{aligned}
\]
In what follows, we derive $b_{k}$s (for all $k=1,...,\beta+1$)
such that 
\[
\left|y^{(k)}(x)\right|\leq b_{k}
\]
and $L_{\beta+1}$ such that 
\[
\left|y^{(\beta+1)}(x)-y^{(\beta+1)}\left(x^{'}\right)\right|\leq L_{\beta+1}\left|x-x^{'}\right|.
\]
\textbf{Proof} \textbf{of Lemma 3.2}. Writing $\vartheta(x)=(x,y(x))$
and $\eta_{j}(t_{1},t_{2})=\frac{\partial^{a_{j}+b_{j}}f}{(\partial t_{1})^{a_{j}}(\partial t_{2})^{b_{j}}}(t_{1},t_{2})$,
we have
\begin{align*}
 & \frac{d}{dx}\vartheta(x)=(1,y^{'}(x))=(1,f(x,y(x))),\\
 & \frac{\partial\eta_{j}}{\partial t_{1}}(t_{1},t_{2})=\frac{\partial^{a_{j}+b_{j}+1}f}{(\partial t_{1})^{a_{j}+1}(\partial t_{2})^{b_{j}}}(t_{1},t_{2}),\\
 & \frac{\partial\eta_{j}}{\partial t_{2}}(t_{1},t_{2})=\frac{\partial^{a_{j}+b_{j}+1}f}{(\partial t_{1})^{a_{j}}(\partial t_{2})^{b_{j}+1}}(t_{1},t_{2}).
\end{align*}
Hence,
\begin{align}
 & \frac{d}{dx}\left(\frac{\partial^{a_{j}+b_{j}}f}{(\partial t_{1})^{a_{j}}(\partial t_{2})^{b_{j}}}(x,y(x))\right)\nonumber \\
 & =\frac{d}{dx}\left(\eta_{j}\circ\vartheta(x))\right)\nonumber \\
 & =\left(\frac{\partial\eta_{j}}{\partial t_{1}}(x,y(x)),\frac{\partial\eta_{j}}{\partial t_{2}}(x,y(x))\right)\boldsymbol{\cdot}(1,f(x,y(x)))\nonumber \\
 & =\frac{\partial^{a_{j}+b_{j}+1}f}{(\partial t_{1})^{a_{j}+1}(\partial t_{2})^{b_{j}}}(x,y(x))+\frac{\partial^{a_{j}+b_{j}+1}f}{(\partial t_{1})^{a_{j}}(\partial t_{2})^{b_{j}+1}}(x,y(x))f(x,y(x)).\label{eq:a6}
\end{align}
For a $\beta-$times differentiable function $f,$ define
\begin{align*}
f^{(a_{1},b_{1}\cdots,a_{k},b_{k})}(x):=\prod_{i=1}^{k}\frac{\partial^{a_{i}+b_{i}}f}{(\partial t_{1})^{a_{i}}(\partial t_{2})^{b_{i}}}(x,y(x))
\end{align*}
where $\sum_{i=1}^{k}a_{i}+\sum_{i=1}^{k}b_{i}=k-1$ for all $k=1,...,\beta+1$,
and $a_{i},\,b_{i}\geq0$ are all integers. We have

\begin{align*}
 & \frac{d}{dx}f^{(a_{1},b_{1}\cdots,a_{k},b_{k})}(x)\\
= & \sum_{j=1}^{k}\frac{\partial^{a_{1}+b_{1}}f}{(\partial t_{1})^{a_{1}}(\partial t_{2})^{b_{1}}}(x,y(x))\cdot\ldots\cdot\frac{d}{dx}\left(\frac{\partial^{a_{j}+b_{j}}f}{(\partial t_{1})^{a_{j}}(\partial t_{2})^{b_{j}}}(x,y(x))\right)\cdot\cdots\frac{\partial^{a_{k}+b_{k}}f}{(\partial t_{1})^{a_{k}}(\partial t_{2})^{b_{k}}}(x,y(x))\\
= & \sum_{j=1}^{k}\frac{\partial^{a_{1}+b_{1}}f}{(\partial t_{1})^{a_{1}}(\partial t_{2})^{b_{1}}}(x,y(x))\cdot\ldots\cdot\\
 & ~~~~\left[\frac{\partial^{a_{j}+b_{j}+1}f}{(\partial t_{1})^{a_{j}+1}(\partial t_{2})^{b_{j}}}(x,y(x))+\frac{\partial^{a_{j}+b_{j}+1}f}{(\partial t_{1})^{a_{j}}(\partial t_{2})^{b_{j}+1}}(x,y(x))\boldsymbol{\cdot}f(x,y(x))\right]\cdot\ldots\cdot\\
 & ~~~~~\frac{\partial^{a_{k}+b_{k}}f}{(\partial t_{1})^{a_{k}}(\partial t_{2})^{b_{k}}}(x,y(x))
\end{align*}

\noindent where the second equality comes from (\ref{eq:a6}); that
is,

\begin{align}
 & \frac{d}{dx}f^{a_{1},b_{1},\ldots,a_{k},b_{k}}(x)=\nonumber \\
 & \sum_{j=1}^{k}f^{(a_{1},b_{1},\ldots,a_{j-1},b_{j-1},a_{j}+1,b_{j},a_{j+1},b_{j+1},\ldots,a_{k},b_{k})}+\sum_{j=1}^{k}f^{(a_{1},b_{1},\ldots,a_{j-1},b_{j-1},a_{j},b_{j}+1,0,0,a_{j+1},b_{j+1},\ldots,a_{k},b_{k})}.\label{eq:a7}
\end{align}
Given (\ref{eq:a7}), now we show that

\begin{align}
y^{(k)}(x)=\sum_{i=1}^{M_{k-1}}f^{(a_{1}^{i},b_{1}^{i},\ldots,a_{k}^{i},a_{k}^{i})}(x)\label{eq:a12}
\end{align}
satisfies the following properties: (1) the integer $M_{k-1}\leq2^{k-1}(k-1)!$
for each $1\leq k\leq\beta+1$, and (2) $\forall i,\,a_{1}^{i}+b_{1}^{i}+\ldots+a_{k}^{i}+b_{k}^{i}=k-1$.

The base case is obvious. Notice that by (\ref{eq:a7}), differentiating
a term of the form $f^{(a_{1}^{i},b_{1}^{i},\ldots,a_{k}^{i},b_{k}^{i})}(x)$
gives us $2k$ terms of the form $f^{(c_{1}^{i},d_{1}^{i},\ldots,c_{m}^{i},d_{m}^{i})}(x)$
where $m=k$ or $k+1,$ and $c_{1}^{i}+d_{1}^{i}+\ldots+c_{m}^{i}+d_{m}^{i}=1+(a_{1}^{i}+b_{1}^{i}+\ldots+a_{k}^{i}+b_{k}^{i})=1+(k-1)=k.$
So the total number of terms in $y^{(k+1)}\leq2k$(number of terms
in $y^{k}$) $\leq2k(k-1)!2^{k-1}.$ This proves the induction hypothesis.
By the assumption that $\left|f(x,\,y(x))\right|\leq1$, we have $\left|y^{(k)}(x)\right|\leq2^{k-1}(k-1)!$.

To show the second claim $\left|y^{(\beta+1)}(x)-y^{(\beta+1)}(x^{'})\right|\leq2^{\beta+1}\left(\beta+1\right)!|x-x^{'}|$,
we use (\ref{eq:a12}). Note that 
\begin{eqnarray*}
\vert y^{(\beta+1)}(x)-y^{(\beta+1)}(x^{'})\vert & = & \left|\sum_{i=1}^{M_{\beta}}f^{(a_{1}^{i},b_{1}^{i},\ldots,a_{\beta+1}^{i},b_{\beta+1}^{i})}(x)-\sum_{i=1}^{M_{\beta}}f^{(a_{1}^{i},b_{1}^{i},\ldots,a_{\beta+1}^{i},b_{\beta+1}^{i})}(x^{'})\right|\\
 & \leq & \sum_{i=1}^{M_{\beta}}\left|f^{(a_{1}^{i},b_{1}^{i},\ldots,a_{\beta+1}^{i},b_{\beta+1}^{i})}(x)-f^{(a_{1}^{i},b_{1}^{i},\ldots,a_{\beta+1}^{i},b_{\beta+1}^{i})}(x^{'})\right|
\end{eqnarray*}
where the integer $M_{\beta}\leq2^{\beta}\beta!$ and
\begin{align*}
 & \Big|f^{(a_{1}^{i},b_{1}^{i},\ldots,a_{\beta+1}^{i},b_{\beta+1}^{i})}(x)-f^{(a_{1}^{i},b_{1}^{i},\ldots,a_{\beta+1}^{i},b_{\beta+1}^{i})}(x^{'})\Big|\\
= & \Big|\sum_{j=1}^{\beta+1}\frac{\partial^{a_{1}^{i}+b_{1}^{i}}f}{(\partial t_{1})^{a_{1}^{i}}(\partial t_{2})^{b_{1}^{i}}}(x,y(x))\cdot\ldots\\
 & ~~~~~~~\Big|\frac{\partial^{a_{j-1}^{i}+b_{j-1}^{i}}f}{(\partial t_{1})^{a_{j-1}^{i}}(\partial t_{2})^{b_{j-1}^{i}}}(x,y(x))\Big|\cdot\\
 & ~~~~~\left(\frac{\partial^{a_{j}^{i}+b_{j}^{i}}f}{(\partial t_{1})^{a_{j}^{i}}(\partial t_{2})^{b_{j}^{i}}}(x,y(x))-\frac{\partial^{a_{j}^{i}+b_{j}^{i}}f}{(\partial t_{1})^{a_{j}^{i}}(\partial t_{2})^{b_{j}^{i}}}(x^{'},y(x^{'}))\right)\cdot\\
 & ~~~~~~~\Big|\frac{\partial^{a_{j+1}^{i}+b_{j+1}^{i}}f}{(\partial t_{1})^{a_{j+1}^{i}}(\partial t_{2})^{b_{j+1}^{i}}}(x^{'},y(x^{'}))\Big|\cdot\ldots\\
 & ~~~~~~~\frac{\partial^{a_{\beta+1}^{i}+b_{\beta+1}^{i}}f}{(\partial t_{1})^{a_{\beta+1}^{i}}(\partial t_{2})^{b_{\beta+1}^{i}}}(x^{'},y(x^{'}))\Big|\\
\leq & \sum_{j=1}^{\beta+1}\Big|\frac{\partial^{a_{1}^{i}+b_{1}^{i}}f}{(\partial t_{1})^{a_{1}^{i}}(\partial t_{2})^{b_{1}^{i}}}(x,y(x))\Big|\cdot\ldots\\
 & ~~~~~~~\Big|\frac{\partial^{a_{j-1}^{i}+b_{j-1}^{i}}f}{(\partial t_{1})^{a_{j-1}^{i}}(\partial t_{2})^{b_{j-1}^{i}}}(x,y(x))\Big|\cdot\\
 & ~~~~~\underset{\leq2\left|x-x^{'}\right|\textrm{ by (\ref{eq:a3-2-1}) and (\ref{eq:a2-2-1})}}{\underbrace{\left|\left(\frac{\partial^{a_{j}^{i}+b_{j}^{i}}f}{(\partial t_{1})^{a_{j}^{i}}(\partial t_{2})^{b_{j}^{i}}}(x,y(x))-\frac{\partial^{a_{j}^{i}+b_{j}^{i}}f}{(\partial t_{1})^{a_{j}^{i}}(\partial t_{2})^{b_{j}^{i}}}(x^{'},y(x^{'}))\right)\right|}}\cdot\\
 & ~~~~~~~\Big|\frac{\partial^{a_{j+1}^{i}+b_{j+1}^{i}}f}{(\partial t_{1})^{a_{j+1}^{i}}(\partial t_{2})^{b_{j+1}^{i}}}(x^{'},y(x^{'}))\Big|\cdot\ldots\\
 & \ ~~~~~~~\Big|\frac{\partial^{a_{\beta+1}^{i}+b_{\beta+1}^{i}}f}{(\partial t_{1})^{a_{\beta+1}^{i}}(\partial t_{2})^{b_{\beta+1}^{i}}}(x^{'},y(x^{'}))\Big|\leq2\left(\beta+1\right)\left|x-x^{'}\right|.
\end{align*}
Hence, 
\begin{align*}
\left|y^{(\beta+1)}(x)-y^{(\beta+1)}(x^{'})\right|\leq\left(2^{\beta}\beta!\right)\left(2\left(\beta+1\right)\left|x-x^{'}\right|\right) & \leq2^{\beta+1}\left(\beta+1\right)!|x-x^{'}|.
\end{align*}

$\square$

\subsection{Theorems 3.1 and 3.2\label{subsec:Theorems-3.1-and-3.2}}

\subsubsection{Main argument}

In what follows, we show Theorem 3.2 and point out the minor differences
in the proof for Theorem 3.1 at the end. 

\subsubsection*{Theorem 3.2}

\textbf{Proof}. \textbf{Term related to }$W_{1}(\delta,\gamma)$:
In Lemma A.1, we establish the following bound:
\begin{equation}
N_{\infty}\left(\delta,\,\mathcal{S}_{\beta+2,\overline{C}}^{\dagger}\right)\leq\exp\left[\log\left(\overline{C}\prod_{i=0}^{\gamma}i!\right)+\frac{\gamma^{2}+\gamma}{2}\log2+\frac{\gamma+3}{2}\log\frac{5}{\delta}+\alpha\left(\frac{\delta}{5}\right)^{\frac{-1}{\gamma+2}}\log4+\log4\right]\label{eq:w1}
\end{equation}
which is valid for all $\gamma\in\left\{ 0,...,\beta\right\} $. By
Lemma 3.2, $\mathcal{Y}\subseteq\mathcal{S}_{\beta+2,\overline{C}}^{\dagger}$,
and therefore, choosing $\gamma$ that minimizes the RHS of (\ref{eq:w1})
yields 
\begin{equation}
\log N_{\infty}\left(\delta,\,\mathcal{Y}\right)\leq\min_{\gamma\in\left\{ 0,...,\beta\right\} }W_{1}\left(\delta,\gamma\right).\label{eq:76}
\end{equation}

\textbf{Term related to} $W_{2}(\delta,\gamma)$: In Lemma A.2, we
establish the following bound: 
\[
N_{\infty}\left(\delta,\,\mathcal{S}_{\beta+1,\,2}\left(1,\,\Xi\right)\right)\leq\exp\left[\frac{\left(\gamma+2\right)\left(\gamma+3\right)}{6}\log\frac{5}{\delta}+20\left(\overline{C}\vee1\right)\left(\frac{\delta}{5}\right)^{\frac{-2}{\gamma+1}}\log2+4\log2\right]
\]
which is valid for all $\gamma\in\left\{ 0,...,\beta\right\} $. Let
$L_{\max}=\sup_{x\in\left[0,\,\alpha\right]}\left\{ \exp\left(x\right)\left[1+\int_{0}^{x}\exp\left(-s\right)ds\right]\right\} $.
For a given $\delta>0$, let us consider the smallest $\frac{\delta}{L_{\max}}-$covering
$\left\{ f_{1},...,f_{N}\right\} $ of $\mathcal{S}_{\beta+1,\,2}\left(1,\,\Xi\right)$
with respect to the sup-norm such that 
\begin{equation}
\log N_{\infty}\left(\frac{\delta}{L_{\max}},\,\mathcal{S}_{\beta+1,\,2}\left(1,\,\Xi\right)\right)\leq\frac{\left(\gamma+2\right)\left(\gamma+3\right)}{6}\log\frac{5}{\delta}\log\frac{5L_{\max}}{\delta}+20\left(\overline{C}\vee1\right)\left(\frac{\delta}{5L_{\max}}\right)^{\frac{-2}{\gamma+1}}\log2+4\log2\label{eq:w2}
\end{equation}
valid for all $\gamma\in\left\{ 0,...,\beta\right\} $. We also consider
the smallest $\frac{\delta}{L_{\max}}-$covering $\left\{ y_{0,1},...,y_{0,N^{'}}\right\} $
for the interval $\left[-C_{0},\,C_{0}\right]$ where the initial
value lies. Note that 
\[
\log N_{\infty}\left(\frac{\delta}{L_{\max}},\,\left[-C_{0},\,C_{0}\right]\right)\leq\log\left(\frac{C_{0}L_{\max}}{\delta}+1\right).
\]
By Theorem B.1 in Section \ref{sec:Supporting-lemmas-and}, for a
solution $y$ to the ODE associated with any $f\in\mathcal{S}_{\beta+1,\,2}\left(1,\,\Xi\right)$
and $y_{0}\in\left[-C_{0},\,C_{0}\right]$, we can find $i\in\left\{ 1,...,N\right\} $
and $i^{'}\in\left\{ 1,...,N^{'}\right\} $ such that 
\[
\left|y\left(x\right)-y_{(i,i^{'})}\left(x\right)\right|\leq\frac{\delta}{L_{\max}}\exp\left(x\right)\left[1+\int_{0}^{x}\exp\left(-s\right)ds\right]\leq\delta\qquad\forall x\in\left[0,\,\alpha\right]
\]
where $y_{(i,i^{'})}$ is a solution to the ODE associated with $f_{i}$
and the initial value $y_{0,i^{'}}$ from the covering sets. Consequently,
we obtain a $\delta-$cover of $\mathcal{Y}$. Choosing $\gamma$
that minimizes the RHS of (\ref{eq:w2}), we conclude that 
\begin{equation}
\log N_{\infty}\left(\delta,\,\mathcal{Y}\right)\leq\min_{\gamma\in\left\{ 0,...,\beta\right\} }W_{2}\left(\frac{\delta}{L_{\max}},\gamma\right).\label{eq:78}
\end{equation}

Combining (\ref{eq:76}) and (\ref{eq:78}) yields the bound on $\log N_{\infty}\left(\delta,\,\mathcal{Y}\right)$
in Theorem 3.2. We can apply argument almost identical to the above
(related to $W_{1}(\delta,\gamma)$) to derive the bound on $\log N_{\infty}\left(\delta,\,\mathcal{Y}_{1}\right)$.
In particular, we have
\[
\log N_{\infty}\left(\delta,\,\mathcal{Y}_{1}\right)\leq\min_{\gamma\in\left\{ 0,...,\beta\right\} }W_{3}\left(\delta,\gamma\right).
\]

$\square$

\subsubsection*{Theorem 3.1 }

\textbf{Proof}. With slight modifications, the arguments for Theorem
3.1 are almost identical to those for Theorem 3.2. Because most of
these modifications are straightforward, we only point out the main
differences.\\
\\
\textbf{1}. By Lemma 3.1, $\mathcal{Y}\subseteq\mathcal{AS}_{\beta+2,\overline{C}}^{\dagger}$.
Following the argument in the proof for Lemma A.1, we establish the
following bound:
\begin{equation}
N_{\infty}\left(\delta,\,\mathcal{AS}_{\beta+2,\overline{C}}^{\dagger}\right)\leq\exp\left[\log\left(\overline{C}\prod_{i=0}^{\gamma}i!\right)+\frac{\gamma+3}{2}\log\frac{5}{\delta}+\alpha\left(\frac{\delta}{5}\right)^{\frac{-1}{\gamma+2}}\log2+\log4\right]\label{eq:crude-1-1}
\end{equation}
which is valid for all $\gamma\in\left\{ 0,...,\beta\right\} $.\\
\textbf{}\\
\textbf{2. }In applying Theorem 2.1, we exploit the following bound:
\[
N_{\infty}\left(\delta,\,\mathcal{S}_{\beta+1,\,1}\left(1,\,\left[-\overline{C},\,\overline{C}\right]\right)\right)\leq\exp\left[\frac{\gamma+2}{2}\log\frac{5}{\delta}+2\overline{C}\left(\frac{\delta}{5}\right)^{\frac{-1}{\gamma+1}}\log2+\log4\right]
\]
which is valid for all $\gamma\in\left\{ 0,...,\beta\right\} $. $\square$

\subsubsection{Lemma A.1\label{subsec:Lemma-A.1}}

\textbf{Lemma A.1}. \textit{In terms of $\mathcal{S}_{\beta+2,\overline{C}}^{\dagger}$,
we have (\ref{eq:w1}). }\\
\\
\textbf{Proof}. By Lemma 3.2, we have shown that $\mathcal{Y}\subseteq\mathcal{S}_{\beta+2,\overline{C}}^{\dagger}$.
Note that $\mathcal{S}_{\beta+2,\overline{C}}^{\dagger}\subseteq\mathcal{S}_{\beta+1,\overline{C}}^{\dagger}\subseteq\cdots\subseteq\mathcal{S}_{2,\overline{C}}^{\dagger}$.
In what follows, we provide an upper bound on the covering number
of $\mathcal{S}_{\gamma+2,\overline{C}}^{\dagger}$ for any given
$\gamma\in\left\{ 0,...,\beta\right\} $. The argument modifies the
original proof for smooth classes in \cite{Kolmogorov and Tikhomirov 1959-1}.
For every function $h\in\mathcal{S}_{\gamma+2,\overline{C}}^{\dagger}$
and $x,\,x+\Delta\in(0,\alpha),$ we have
\[
h(x+\Delta)=h(x)+\Delta h^{'}(x)+\frac{\Delta^{2}}{2!}h^{''}(x)+\cdots+\frac{\Delta^{\gamma}}{\gamma!}h^{(\gamma)}(x)+\frac{\Delta^{\gamma+1}}{(\gamma+1)!}h^{(\gamma+1)}(z).
\]
Let us define 
\[
R_{h}(x,\Delta):=h(x+\Delta)-h(x)-\Delta h^{'}(x)-\frac{\Delta^{2}}{2!}h^{''}(x)-\cdots-\frac{\Delta^{\gamma}}{\gamma!}h^{(\gamma)}(x)-\frac{\Delta^{\gamma+1}}{(\gamma+1)!}h^{(\gamma+1)}(x)
\]
and we have 
\begin{equation}
\left|R_{h}(x,\Delta)\right|=\frac{\Delta^{\gamma+1}}{(\gamma+1)!}\left|h^{(\gamma+1)}(x)-h^{(\gamma+1)}(z)\right|\leq(2|\Delta|)^{\gamma+2}.\label{eq:rem}
\end{equation}
As a consequence, we obtain

\begin{align*}
h(x+\Delta)=\sum_{k=0}^{\gamma+1}\frac{\Delta^{k}}{k!}h^{(k)}(x)+R_{h}(x,\Delta)\quad\textrm{where }\left|R_{h}(x,\Delta)\right|\leq(2|\Delta|)^{\gamma+2}.
\end{align*}
Let us consider $h^{(i)}\in\mathcal{S}_{\gamma+2-i,1}^{\dagger}$
for $1\leq i\leq\gamma+1$ and the above implies that
\begin{align}
 & h^{(i)}(x+\Delta)=\sum_{k=0}^{\gamma+1-i}\frac{\Delta^{k}}{k!}h^{(i+k)}(x)+R_{h^{(i)}}(x,\Delta)~\textrm{ where }\left|R_{h^{(i)}}(x,\Delta)\right|\leq(2|\Delta|)^{\gamma+2-i}.\label{a30}
\end{align}

To bound $N_{\infty}\left(\delta,\,\mathcal{S}_{\gamma+2,\overline{C}}^{\dagger}\right)$
from above, we fix $\delta>0$ and $x\in(0,\alpha).$ Suppose that
for some $\delta_{0},\ldots,\delta_{\gamma+1}>0$, $h,\,g\,\in\mathcal{S}_{\gamma+2,\overline{C}}^{\dagger}$
satisfy

\[
\left|h^{(k)}(x)-g^{(k)}(x)\right|\leq\delta_{k}\quad\text{ for all }k=0,\dots,\gamma+1.
\]
We can bound $|h(x+\Delta)-g(x+\Delta)|$ from above for some $\Delta$
such that $x+\Delta\in(0,\alpha)$ as follows: 
\begin{eqnarray*}
\left|h(x+\Delta)-g(x+\Delta)\right| & = & \left|\sum_{k=0}^{\gamma+1}\frac{\Delta^{k}}{k!}\left(h^{(k)}(x)-g^{(k)}(x)\right)+R_{h}(x,\Delta)-R_{g}(x,\Delta)\right|\\
 & \leq & \sum_{k=0}^{\gamma+1}\frac{|\Delta|^{k}\delta_{k}}{k!}+2(2|\Delta|)^{\gamma+2}.
\end{eqnarray*}
If $|\Delta|\leq\frac{\left(\frac{\delta}{5}\right)^{\frac{1}{\gamma+2}}}{2}$
and $\delta_{k}=\left(\frac{\delta}{5}\right)^{1-\frac{k}{\gamma+2}}$,
we then have 
\begin{equation}
\left|h(x+\Delta)-g(x+\Delta)\right|\leq\frac{\delta}{5}\left(\sum_{k=0}^{\gamma+1}\frac{1}{k!}+2\right)\leq\delta.\label{eq:a90}
\end{equation}
In other words, by considering a grid of points $\frac{\left(\frac{\delta}{5}\right)^{\frac{1}{\gamma+2}}}{2}-$apart
in $(0,\alpha)$ and covering the $k$th derivative of functions in
$\mathcal{S}_{\gamma+2,\overline{C}}^{\dagger}$ within $\delta_{k}=\left(\frac{\delta}{5}\right)^{1-\frac{k}{\gamma+2}}$
at each point, we can then obtain a $\delta$ cover in the sup-norm
for $\mathcal{S}_{\gamma+2,\overline{C}}^{\dagger}$. Let $x_{1}<\cdots<x_{s}$
be a $\frac{\left(\frac{\delta}{5}\right)^{\frac{1}{\gamma+2}}}{2}-$grid
of points in $(0,1)$ with $s\leq2\alpha\left(\frac{\delta}{5}\right)^{\frac{-1}{\gamma+2}}+2$.
For each $h_{0}\in\mathcal{S}_{\gamma+2,\overline{C}}^{\dagger}$,
let us define

\[
\mathbb{H}\left(h_{0}\right):=\left\{ h\in\mathcal{S}_{\gamma+2,\overline{C}}^{\dagger}:\left\lfloor \frac{h^{(k)}\left(x_{i}\right)}{\delta_{k}}\right\rfloor =\left\lfloor \frac{h_{0}^{(k)}\left(x_{i}\right)}{\delta_{k}}\right\rfloor ,1\leq i\leq s,\,0\leq k\leq\gamma+1\right\} 
\]
where $\left\lfloor x\right\rfloor $ means the largest integer smaller
than or equal to $x$. Our earlier argument implies that the number
of distinct sets $\mathbb{H}\left(h_{0}\right)$ with $h_{0}$ ranging
over $\mathcal{S}_{\gamma+2,\overline{C}}^{\dagger}$ bounds the $\delta-$covering
number of $\mathcal{S}_{\gamma+2,\overline{C}}^{\dagger}$ from above.
Note that for $i=1,\ldots,s$ and $k=0,\ldots,\gamma+1$, $\mathbb{H}\left(h_{0}\right)$
depends on $\left\lfloor h_{0}^{(k)}\left(x_{i}\right)/\delta_{k}\right\rfloor $
only and the number of distinct sets $\mathcal{\mathbb{H}}\left(h_{0}\right)$
is bounded above by the cardinality of 
\begin{equation}
I:=\left\{ \left(\left\lfloor \frac{h^{(k)}\left(x_{i}\right)}{\delta_{k}}\right\rfloor ,1\leq i\leq s\text{ and }0\leq k\leq\gamma+1\right):h\in\mathcal{S}_{\gamma+2,\overline{C}}^{\dagger}\right\} .\label{eq:a50}
\end{equation}

Starting from $x_{1}$, let us count the number of possible values
of the vector 
\[
\left(\left\lfloor \frac{h^{(k)}\left(x_{1}\right)}{\delta_{k}}\right\rfloor ,\,0\leq k\leq\gamma+1\right)
\]
with $h$ ranging over $\mathcal{S}_{\gamma+2,\overline{C}}^{\dagger}$.
Since $\left|h\left(x_{1}\right)\right|\leq\overline{C},\left|h^{(1)}\left(x_{1}\right)\right|\leq1,\left|h^{(k)}\left(x_{1}\right)\right|\leq2^{k-1}(k-1)!$
for $2\leq k\leq\gamma+1$, this number is at most

\begin{equation}
\frac{\overline{C}}{\delta_{0}}\frac{1}{\delta_{1}}\frac{2}{\delta_{2}}\frac{2^{2}2!}{\delta_{3}}\ldots\frac{2^{\gamma-1}(\gamma-1)!}{\delta_{\gamma}}\frac{2^{\gamma}\gamma!}{\delta_{\gamma+1}}\leq\left(\frac{\delta}{5}\right)^{-\frac{\gamma+3}{2}}\overline{C}2^{\frac{\left(\gamma+1\right)\gamma}{2}}\prod_{i=0}^{\gamma}i!\label{a60}
\end{equation}
Now we move to $x_{2}$. Given the values of $\left(\left\lfloor h^{(k)}\left(x_{1}\right)/\delta_{k}\right\rfloor ,\,0\leq k\leq\gamma+1\right)$,
we count the number of possible values of the vector 
\[
\left(\left\lfloor \frac{h^{(k)}\left(x_{2}\right)}{\delta_{k}}\right\rfloor ,\,0\leq k\leq\gamma+1\right).
\]
For each $0\leq k\leq\gamma+1$, we define

\[
A_{k}:=\left\lfloor \frac{h^{(k)}\left(x_{1}\right)}{\delta_{k}}\right\rfloor \quad\textrm{such that }A_{k}\delta_{k}\leq h^{(k)}\left(x_{1}\right)<\left(A_{k}+1\right)\delta_{k}.
\]
Let us fix $0\leq i\leq\gamma+1$. Applying (\ref{a30}) with $x=x_{1}$
and $\Delta=x_{2}-x_{1}$ yields 
\[
\left|h^{(i)}\left(x_{2}\right)-\sum_{k=0}^{\gamma+1-i}\frac{\Delta^{k}}{k!}h^{(i+k)}\left(x_{1}\right)\right|\leq|2\Delta|^{\gamma+2-i}.
\]
Consequently we have
\begin{align*}
 & \left|h^{(i)}\left(x_{2}\right)-\sum_{k=0}^{\gamma+1-i}\frac{\Delta^{k}}{k!}A_{i+k}\right|\\
 & \leq\left|h^{(i)}\left(x_{2}\right)-\sum_{k=0}^{\gamma+1-i}\frac{\Delta^{k}}{k!}h^{(i+k)}\left(x_{1}\right)\right|+\left|\sum_{k=0}^{\gamma+1-i}\frac{\Delta^{k}}{k!}\left(h^{(i+k)}\left(x_{1}\right)-A_{i+k}\right)\right|\\
 & \leq|2\Delta|^{\gamma+2-i}+\sum_{k=0}^{\gamma+1-i}\frac{|\Delta|^{k}}{k!}\delta_{i+k}\\
 & \leq\left(\frac{\delta}{5}\right)^{1-\frac{i}{\gamma+2}}+\left(\frac{\delta}{5}\right)^{1-\frac{i}{\gamma+2}}=2\delta_{i}
\end{align*}
(recalling $|\Delta|=\left|x_{2}-x_{1}\right|=\frac{\left(\frac{\delta}{5}\right)^{\frac{1}{\gamma+2}}}{2}$).
Therefore, given the values of $\left(\left\lfloor h^{(k)}\left(x_{1}\right)/\delta_{k}\right\rfloor ,\,0\leq k\leq\gamma+1\right)$,
$h^{(i)}\left(x_{2}\right)$ takes values in an interval whose length
is no greater than $2\delta_{i}$. Therefore, the number of possible
values of $\left(\left\lfloor \frac{h^{(k)}\left(x_{2}\right)}{\delta_{k}}\right\rfloor ,\,0\leq k\leq\gamma+1\right)$
is at most $2$. The same argument goes through when $x_{1}$ and
$x_{2}$ are replaced with $x_{j}$ and $x_{j+1}$ for any $j=1,...,s-1$.
This result along with (\ref{a60}) gives 
\begin{eqnarray*}
\left|I\right| & \leq & 2^{s}\left(\frac{\delta}{5}\right)^{-\frac{\gamma+3}{2}}\overline{C}2^{\frac{\left(\gamma+1\right)\gamma}{2}}\prod_{i=0}^{\gamma}i!\\
 & \leq & 4^{\alpha\left(\frac{\delta}{5}\right)^{\frac{-1}{\gamma+2}}+1}\left(\frac{\delta}{5}\right)^{-\frac{\gamma+3}{2}}\overline{C}2^{\frac{\left(\gamma+1\right)\gamma}{2}}\prod_{i=0}^{\gamma}i!\\
 & \leq & \exp\left[\log\left(\overline{C}\prod_{i=0}^{\gamma}i!\right)+\frac{\gamma^{2}+\gamma}{2}\log2+\frac{\gamma+3}{2}\log\frac{5}{\delta}+\alpha\left(\frac{\delta}{5}\right)^{\frac{-1}{\gamma+2}}\log4+\log4\right]
\end{eqnarray*}
where $I$ is defined in (\ref{eq:a50}). $\square$

\subsubsection{Lemma A.2\label{subsec:Lemma-A.2}}

\textbf{Lemma A.2}. \textit{In terms of $\mathcal{S}_{\beta+1,\,2}\left(1,\,\Xi\right)$,
we have (\ref{eq:w2}). }\\
\textit{}\\
\textbf{Proof}. Note that $\mathcal{S}_{\beta+1,\,2}\left(1,\,\Xi\right)\subseteq\mathcal{S}_{\beta,\,2}\left(1,\,\Xi\right)\subseteq\cdots\subseteq\mathcal{S}_{1,\,2}\left(1,\,\Xi\right)$.
In what follows, we provide an upper bound on the covering number
of $\mathcal{S}_{\gamma+1,\,2}\left(1,\,\Xi\right)$ for any given
$\gamma\in\left\{ 0,...,\beta\right\} $. As we discussed in the main
paper, while Kolmogorov and Tikhomirov \cite{Kolmogorov and Tikhomirov 1959-1}
explicitly derived $\left(\gamma+1\right)\log\frac{1}{\delta}$ (from
counting the number of distinct values in the first interval on $\left[0,\,1\right]$)
for the univariate function class, moving from univariate functions
to multivariate functions, this type of log terms appears unmentioned
in \cite{Kolmogorov and Tikhomirov 1959-1} possibly because $\gamma$
is assumed to be very small. However, these log terms would clearly
have practical implications. Therefore, in what follows, we carefully
extend the arguments from univariate functions to bivariate functions
and derive the log term explicitly.

Let $p=\left(p_{1},p_{2}\right)$ and $\left[p\right]=p_{1}+p_{2}$
where $p_{1}$ and $p_{2}$ are non-negative integers. We write $D^{p}f\left(w_{1},w_{2}\right)=\partial^{\left[p\right]}f/\partial w_{1}^{p_{1}}\partial w_{2}^{p_{2}}$
with $w=\left(w_{1},\,w_{2}\right)$ and $\Delta^{p}=\Delta_{1}^{p_{1}}\Delta_{2}^{p_{2}}$
with $\Delta=\left(\Delta_{1},\,\Delta_{2}\right)$. For every function
$f\in\mathcal{S}_{\gamma+1,\,2}\left(1,\,\Xi\right)$, $w_{1},\,w_{1}+\Delta_{1}\in(0,1)$,
and $w_{2},\,w_{2}+\Delta_{2}\in\left(-\overline{C},\,\overline{C}\right)$,
we have
\[
f(w+\Delta)=\sum_{k=0}^{\gamma}\sum_{p:\left[p\right]=k}\frac{\Delta^{p}D^{p}f\left(w\right)}{k!}+\underset{R_{0,f}(w,\Delta)}{\underbrace{\sum_{p:\left[p\right]=\gamma}\frac{\Delta^{p}D^{p}f\left(z\right)}{\gamma!}-\sum_{p:\left[p\right]=\gamma}\frac{\Delta^{p}D^{p}f\left(w\right)}{\gamma!}}}
\]
for some $z=\left(z_{1},\,z_{2}\right)\in(0,1)\times\left(-\overline{C},\,\overline{C}\right)$.
Because a function of two variables has $2^{\gamma}$ $\gamma$th
partial derivatives, we have 
\begin{equation}
\left|R_{0,f}(w,\Delta)\right|\leq\frac{\left|2\Delta\right|_{\infty}^{\gamma+1}}{\gamma!}.\label{eq:85}
\end{equation}
Similarly, letting $w+\Delta:=\left(w_{1}+\Delta_{1},\,w_{2}+\Delta_{2}\right)$
and $D^{\tilde{p}}f(w+\Delta)\in\mathcal{S}_{\gamma+1-[\tilde{p}],1}\left(1,\,\Xi\right)$
for $1\leq[\tilde{p}]\leq\gamma$, we have
\begin{align}
 & D^{\tilde{p}}f(w+\Delta)=\sum_{k=0}^{\gamma-[\tilde{p}]}\sum_{p:\left[p\right]=k}\frac{\Delta^{p}D^{p+\tilde{p}}f\left(w\right)}{k!}+\underset{R_{[\tilde{p}],f}(w,\Delta)}{\underbrace{\sum_{p:\left[p\right]=\gamma-[\tilde{p}]}\frac{\Delta^{p}D^{p+\tilde{p}}f\left(\tilde{z}\right)}{\left(\gamma-[\tilde{p}]\right)!}-\sum_{p:\left[p\right]=\gamma-[\tilde{p}]}\frac{\Delta^{p}D^{p+\tilde{p}}f\left(w\right)}{\left(\gamma-[\tilde{p}]\right)!}}}\label{a30-1}
\end{align}
for some $\tilde{z}=\left(\tilde{z}_{1},\,\tilde{z}_{2}\right)\in(0,1)\times\left(-\overline{C},\,\overline{C}\right)$,
where 
\[
\left|R_{[\tilde{p}],f}(w,\Delta)\right|\leq\frac{\left|2\Delta\right|_{\infty}^{\gamma+1-[\tilde{p}]}}{\left(\gamma-[\tilde{p}]\right)!}.
\]

Suppose that for some $\delta_{0},\ldots,\delta_{\gamma}>0$, $f,\,g\,\in\mathcal{S}_{\gamma+1,\,2}\left(1,\,\Xi\right)$
satisfy

\[
\left|D^{p}f\left(w\right)-D^{p}g\left(w\right)\right|\leq\delta_{k}\quad\text{ for all }p\textrm{ with }[p]=k\in\left\{ 0,\dots,\gamma\right\} .
\]
We can bound $|f(w+\Delta)-g(w+\Delta)|$ from above for some $\Delta$
such that $w+\Delta\in(0,1)\times\left(-\overline{C},\,\overline{C}\right)$
as follows: 
\begin{eqnarray*}
\left|f(w+\Delta)-g(w+\Delta)\right| & = & \left|\sum_{k=0}^{\gamma}\sum_{p:\left[p\right]=k}\frac{\Delta^{p}}{k!}\left(D^{p}f\left(w\right)-D^{p}g\left(w\right)\right)+R_{0,f}(w,\Delta)-R_{0,g}(w,\Delta)\right|\\
 & \leq & \sum_{k=0}^{\gamma}\frac{\left|2\Delta\right|_{\infty}^{k}\delta_{k}}{k!}+2\frac{\left|2\Delta\right|_{\infty}^{\gamma+1}}{\gamma!}
\end{eqnarray*}
where the second line follows from (\ref{eq:85}) and the fact that
a function of two variables has $2^{k}$ $k$th partial derivatives.
If $|\Delta|_{\infty}\leq\frac{\left(\frac{\delta}{5}\right)^{\frac{1}{\gamma+1}}}{2}$
and $\delta_{k}=\left(\frac{\delta}{5}\right)^{1-\frac{k}{\gamma+1}}$,
we then have 
\begin{equation}
\left|h(x+\Delta)-g(x+\Delta)\right|\leq\frac{\delta}{5}\left(\sum_{k=0}^{\gamma}\frac{1}{k!}+2\right)\leq\delta.\label{eq:a90-1}
\end{equation}
Let $w_{1,1}<\cdots<w_{1,s_{1}}$ be a $\frac{\left(\frac{\delta}{5}\right)^{\frac{1}{\gamma+1}}}{2}-$grid
of points in $(0,1)$ with $s_{1}\leq2\left(\frac{\delta}{5}\right)^{\frac{-1}{\gamma+1}}+2$
and $w_{2,1}<\cdots<w_{2,s_{2}}$ be a $\frac{\left(\frac{\delta}{5}\right)^{\frac{1}{\gamma+1}}}{2}-$grid
of points in $\left(-\overline{C},\,\overline{C}\right)$ with $s_{2}\leq4\overline{C}\left(\frac{\delta}{5}\right)^{\frac{-1}{\gamma+1}}+2$.
For each $f_{0}\in\mathcal{S}_{\gamma+1,\,2}\left(1,\,\Xi\right)$,
let us define 

\begin{align*}
\mathbb{H}\left(f_{0}\right) & :=\{f\in\mathcal{S}_{\gamma+1,\,2}\left(1,\,\Xi\right):\left\lfloor \frac{D^{p}f\left(w_{1,i},w_{2,j}\right)}{\delta_{k}}\right\rfloor =\left\lfloor \frac{D^{p}f_{0}\left(w_{1,i},w_{2,j}\right)}{\delta_{k}}\right\rfloor ,\\
 & 1\leq i\leq s_{1},\,1\leq j\leq s_{2},\,p\textrm{ with }[p]=k\in\left\{ 0,\dots,\gamma\right\} \}
\end{align*}
where $\left\lfloor x\right\rfloor $ means the largest integer smaller
than or equal to $x$. Our earlier argument implies that the number
of distinct sets $\mathbb{H}\left(f_{0}\right)$ with $f_{0}$ ranging
over $\mathcal{S}_{\gamma+1,\,2}\left(1,\,\Xi\right)$ bounds the
$\delta-$covering number of $\mathcal{S}_{\gamma+1,\,2}\left(1,\,\Xi\right)$
from above. Note that for $i=1,\ldots,s_{1}$, $j=1,\ldots,s_{2}$,
and $p$ with $[p]=k\in\left\{ 0,\dots,\gamma\right\} $, $\mathbb{H}\left(f_{0}\right)$
depends on $\left\lfloor \frac{D^{p}f_{0}\left(w_{1,i},w_{2,j}\right)}{\delta_{k}}\right\rfloor $
only and the number of distinct sets $\mathcal{\mathbb{H}}\left(f_{0}\right)$
is bounded above by the cardinality of 
\begin{eqnarray}
I & = & \{(\left\lfloor \frac{D^{p}f\left(w_{1,i},w_{2,j}\right)}{\delta_{k}}\right\rfloor ,1\leq i\leq s_{1},\,1\leq j\leq s_{2},\,p\textrm{ with }[p]=k\in\left\{ 0,\dots,\gamma\right\} )\nonumber \\
 &  & :f\in\mathcal{S}_{\gamma+1,\,2}\left(1,\,\Xi\right)\}.\label{eq:a50-1}
\end{eqnarray}

Starting from $\left(w_{1,1},w_{2,1}\right)$, let us count the number
of possible values of the vector 
\begin{equation}
\left(\left\lfloor \frac{D^{p}f\left(w_{1,1},w_{2,1}\right)}{\delta_{k}}\right\rfloor ,\,p\textrm{ with }[p]=k\in\left\{ 0,\dots,\gamma\right\} \right)\label{eq:first_point}
\end{equation}
with $f$ ranging over $\mathcal{S}_{\gamma+1,\,2}\left(1,\,\Xi\right)$.
Given a function of two variables has $k+1$ distinct $k$th partial
derivatives (if all $k$th partial derivatives are continuous), this
number is at most 
\begin{eqnarray}
\left(\frac{1}{\delta_{0}}\right)^{1}\left(\frac{1}{\delta_{1}}\right)^{2}\ldots\left(\frac{1}{\delta_{\gamma}}\right)^{\gamma+1} & \leq & \left(\frac{5}{\delta}\right)^{\sum_{j=0}^{\gamma}\left(j+1\right)-\sum_{j=0}^{\gamma}\frac{j\left(j+1\right)}{\gamma+1}}\nonumber \\
 & = & \left(\frac{5}{\delta}\right)^{\frac{\left(\gamma+2\right)\left(\gamma+3\right)}{6}}\label{eq:a60-1}
\end{eqnarray}
The diagram below shows how we move from $\left(w_{1,1},w_{2,1}\right)$
to $\left(w_{1,s_{1}},w_{2,s_{2}}\right)$ to count the the number
of possible values of $\left\lfloor \frac{D^{p}f\left(w_{1,i},w_{2,j}\right)}{\delta_{k}}\right\rfloor $
given that of its previous adjacent pair of points: 
\begin{eqnarray*}
\left(w_{1,1},w_{2,1}\right)\rightarrow & \left(w_{1,2},w_{2,1}\right)\rightarrow\cdots\rightarrow & \left(w_{1,s_{1}},w_{2,1}\right)\\
 &  & \downarrow\\
\left(w_{1,1},w_{2,2}\right)\leftarrow\cdots\leftarrow & \left(w_{1,s_{1}-1},w_{2,2}\right)\leftarrow & \left(w_{1,s_{1}},w_{2,2}\right)\\
\vdots & \vdots & \vdots
\end{eqnarray*}
Note that in each row $j$, we fix $w_{2,j}$ and cycle through all
$w_{1,i}$s; moreover, every time we move from one pair of points
to the next, we only change one coordinate and keep the other fixed.
Because of this construction, we can argue in a similar way as in
the proof for Lemma A.1 that, 
\begin{eqnarray*}
\left|I\right| & \leq & 2^{s_{1}s_{2}}\left(\frac{\delta}{5}\right)^{\frac{\left(\gamma+2\right)\left(\gamma+3\right)}{6}}\\
 & \leq & 2^{\left[2\left(\frac{\delta}{5}\right)^{\frac{-1}{\gamma+1}}+2\right]\left[4\overline{C}\left(\frac{\delta}{5}\right)^{\frac{-1}{\gamma+1}}+2\right]}\left(\frac{5}{\delta}\right)^{\frac{\left(\gamma+2\right)\left(\gamma+3\right)}{6}}\\
 & \leq & \exp\left[\frac{\left(\gamma+2\right)\left(\gamma+3\right)}{6}\log\frac{5}{\delta}+20\left(\overline{C}\vee1\right)\left(\frac{\delta}{5}\right)^{\frac{-2}{\gamma+1}}\log2+4\log2\right]
\end{eqnarray*}
where $I$ is defined in (\ref{eq:a50-1}). $\square$

\subsection{Theorem 3.3\label{subsec:Theorem-3.3}}

\subsubsection*{Preliminary }

In terms of our least squares (\ref{eq:nonpara ls}), the \textit{basic
inequality} 
\[
\frac{1}{2n}\sum_{i=1}^{n}\left(Y_{i}-\hat{y}\left(x_{i}\right)\right)^{2}\leq\frac{1}{2n}\sum_{i=1}^{n}\left(Y_{i}-y\left(x_{i}\right)\right)^{2}
\]
yields 
\begin{equation}
\frac{1}{n}\sum_{i=1}^{n}\left(\hat{y}(x_{i})-y(x_{i})\right)^{2}\leq\frac{2}{n}\sum_{i=1}^{n}\varepsilon_{i}\left(\hat{y}(x_{i})-y(x_{i})\right).\label{eq:44}
\end{equation}
To bound the right-hand-side of (\ref{eq:44}), it suffices to bound
$\mathcal{G}_{n}(\tilde{r}_{n};\,\bar{\mathcal{F}})$ defined in (\ref{eq:13-1})
with $\mathcal{F}=\mathcal{Y}$ in (\ref{eq:28}). 

One way to bound $\mathcal{G}_{n}(\tilde{r}_{n};\,\bar{\mathcal{F}})$
is to seek a sharp enough $\tilde{r}_{n}>0$ that satisfies the \textit{critical
inequality} 
\begin{equation}
\mathcal{G}_{n}\left(\tilde{r}_{n};\,\bar{\mathcal{F}}\right)\precsim\frac{\tilde{r}_{n}^{2}}{\sigma}.\label{eq:14-1-1}
\end{equation}
It is known that the complexity $\mathcal{G}_{n}\left(\tilde{r}_{n};\,\bar{\mathcal{F}}\right)$
can be bounded above by the Dudley's entropy integral (see, e.g.,
\cite{Ledoux and Talagrand 1991,van de Geer 2000,wainwright 2019}).
Let $N_{n}(\delta,\,\Lambda(\tilde{r}_{n};\,\bar{\mathcal{F}}))$
denote the $\delta-$covering number of the set $\Lambda(\tilde{r}_{n};\,\bar{\mathcal{F}})$
in the $\left|\cdot\right|_{n}$ norm. Then by Corollary 13.7 in \cite{wainwright 2019},
the critical radius condition (\ref{eq:14-1-1}) holds for any $\tilde{r}_{n}\in(0,\,\sigma]$
such that 
\begin{equation}
\frac{c}{\sqrt{n}}\int_{\frac{\tilde{r}_{n}^{2}}{4\sigma}}^{\tilde{r}_{n}}\sqrt{\log N_{n}(\delta,\,\Lambda(\tilde{r}_{n};\,\bar{\mathcal{F}}))}d\delta\leq\frac{\tilde{r}_{n}^{2}}{\sigma}.\label{eq:critical-1}
\end{equation}
\textbf{Proof}. In what follows, we show Theorem 3.3(ii) as the argument
for Theorem 3.3(i) is nearly identical; at the end, we show Theorems
3.3(iii)-(iv). 

\textbf{Claim (ii)}. Let us begin with the part $\mathcal{M}_{1}$.
Since $\mathcal{Y}\subseteq\mathcal{S}_{\beta+2,\overline{C}}^{\dagger}$,
we let $\mathcal{F}=\mathcal{S}_{\beta+2,\overline{C}}^{\dagger}$
in (\ref{eq:28}). Because we are working with $\bar{\mathcal{F}}$
in (\ref{eq:28}) in terms of $\mathcal{F}=\mathcal{S}_{\beta+2,\overline{C}}^{\dagger}$,
all the coefficients $2^{k-1}\left(k-1\right)!$ associated with the
derivatives for $k=1,...,\beta+1$, the coefficient $2^{\beta+1}\left(\beta+1\right)!$
associated with the Lipschitz condition, as well as the function value
itself are multiplied by $2$. Slight modifications of the proof for
Lemma A.1 give 
\begin{eqnarray*}
\left|I\right| & \leq & 3^{s}\left(\frac{\delta}{7}\right)^{-\frac{\gamma+3}{2}}\overline{C}2^{\frac{\left(\gamma+1\right)\gamma}{2}+\gamma+2}\prod_{i=0}^{\gamma}i!\\
 & \leq & 9^{\left(\frac{\delta}{7}\right)^{\frac{-1}{\gamma+2}}+1}\left(\frac{\delta}{7}\right)^{-\frac{\gamma+3}{2}}\overline{C}2^{\frac{\left(\gamma+1\right)\gamma}{2}+\gamma+2}\prod_{i=0}^{\gamma}i!\\
 & \leq & \exp\left[\log\left(\overline{C}\prod_{i=0}^{\gamma}i!\right)+\frac{\gamma^{2}+3\gamma}{2}\log2+\frac{\gamma+3}{2}\log\frac{7}{\delta}+\left(\frac{\delta}{7}\right)^{\frac{-1}{\gamma+2}}\log9+\log36\right]
\end{eqnarray*}
valid for all $\gamma=0,...,\beta$. Note that 
\begin{eqnarray*}
 &  & \frac{1}{\sqrt{n}}\int_{\frac{\tilde{r}_{n}^{2}}{4\sigma}}^{\tilde{r}_{n}}\sqrt{\log N_{n}(\delta,\,\Lambda(\tilde{r}_{n};\,\bar{\mathcal{F}}))}d\delta\\
 & \leq & \frac{1}{\sqrt{n}}\int_{0}^{\tilde{r}_{n}}\sqrt{\log N_{\infty}\left(\delta,\,\bar{\mathcal{F}}\right)}d\delta\\
 & \leq & \frac{\tilde{r}_{n}}{\sqrt{n}}\left[\log\left(\overline{C}\prod_{i=0}^{\gamma}i!\right)+\frac{\gamma^{2}+3\gamma}{2}\log2\right]^{\frac{1}{2}}+\left(\frac{\gamma+3}{2n}\right)^{\frac{1}{2}}\int_{0}^{\tilde{r}_{n}}\sqrt{\log\frac{7}{\delta}}d\delta\\
 &  & +\left(\frac{\log9}{n}\right)^{\frac{1}{2}}\int_{0}^{\tilde{r}_{n}}\sqrt{\left(\frac{\delta}{7}\right)^{\frac{-1}{\gamma+2}}}d\delta+\tilde{r}_{n}\sqrt{\frac{\log36}{n}}\\
 & \leq & \underset{\mathcal{T}_{1}\left(\gamma,\tilde{r}_{n}\right)}{\underbrace{c_{0}\left\{ \tilde{r}_{n}\sqrt{\frac{1}{n}\log\left(\prod_{i=0}^{\gamma}i!\right)}+\tilde{r}_{n}\sqrt{\frac{\gamma^{2}}{n}}+\tilde{r}_{n}\sqrt{\frac{1}{n}}+\frac{1}{\sqrt{n}}\tilde{r}_{n}^{\frac{2\gamma+3}{2\gamma+4}}\right\} }}
\end{eqnarray*}
valid for all $\gamma=0,...,\beta$. Therefore, we can take $\min_{\gamma\in\left\{ 0,...,\beta\right\} }\mathcal{T}_{1}\left(\gamma,\tilde{r}_{n}\right)=\mathcal{B}_{1}\left(\tilde{r}_{n}\right)$
and let $\gamma_{1}^{*}$ be the minimizer. Setting $\mathcal{B}_{1}\left(\tilde{r}_{n}\right)\asymp\frac{\tilde{r}_{n}^{2}}{\sigma}$
yields $\mathcal{M}_{1}\left(\gamma_{1}^{*}\right)$.

We now show the part $\mathcal{M}_{2}$. For a given $\delta>0$,
let us consider the smallest $\frac{\delta}{2L_{\max}}-$covering
$\left\{ f^{1},...,f^{N}\right\} $ (w.r.t. the sup-norm) of $\mathcal{S}_{\beta+1,\,2}\left(1,\,\Xi\right)$
and the smallest $\frac{\delta}{2L_{\max}}-$covering $\left\{ y_{0,1},...,y_{0,N^{'}}\right\} $
for the interval $\left[-C_{0},\,C_{0}\right]$ where the initial
value lies. By Theorem B.1 in Section \ref{sec:Supporting-lemmas-and}
and arguments similar to those in the proof for Corollary 2.1, for
any $f,\,\tilde{f}\in\mathcal{S}_{\beta+1,\,2}\left(1,\,\Xi\right)$
and initial values $y_{0},\,\tilde{y}_{0}\in\left[-C_{0},\,C_{0}\right]$,
we can find some $f^{i},\,f^{j}\in\left\{ f^{1},...,f^{N}\right\} $
and $y_{0,i^{'}},\,y_{0,j^{'}}\in\left\{ y_{0,1},...,y_{0,N^{'}}\right\} $
such that 
\begin{eqnarray*}
 &  & \left|y\left(x\right)-\tilde{y}\left(x\right)-\left(y_{(i,i^{'})}\left(x\right)-y_{(j,j^{'})}\left(x\right)\right)\right|\\
 & \leq & \left|y\left(x\right)-y_{(i,i^{'})}\left(x\right)\right|+\left|\tilde{y}\left(x\right)-y_{(j,j^{'})}\left(x\right)\right|\\
 & \leq & \delta
\end{eqnarray*}
where $y$, $\tilde{y}$, $y_{(i,i^{'})}$, and $y_{(j,j^{'})}$ are
solutions to the ODE associated with $\left\{ f,\,y_{0}\right\} $,
$\left\{ \tilde{f},\,\tilde{y}_{0}\right\} $, $\left\{ f^{i},\,y_{0,i^{'}}\right\} $
and $\left\{ f^{j},\,y_{0,j^{'}}\right\} $, respectively. Thus, we
obtain a $\delta-$cover of $\mathcal{\bar{\mathcal{F}}}$ in terms
of $\mathcal{F}=\mathcal{Y}$ in (\ref{eq:28}). The rest of arguments
are very similar to those for $\mathcal{M}_{1}$. For $\beta>0$,
we have 

\begin{eqnarray*}
\frac{1}{\sqrt{n}}\int_{0}^{\tilde{r}_{n}}\sqrt{\log N_{n}(\delta,\,\Lambda(\tilde{r}_{n};\,\bar{\mathcal{F}}))}d\delta & \leq & \frac{1}{\sqrt{n}}\int_{0}^{\tilde{r}_{n}}\sqrt{\log N_{\infty}\left(\delta,\,\bar{\mathcal{F}}\right)}d\delta\\
 & \leq & \underset{\mathcal{T}_{2}\left(\gamma,\tilde{r}_{n}\right)}{\underbrace{c_{1}\tilde{r}_{n}\sqrt{\frac{\gamma^{2}\vee1}{n}}+\frac{1}{\sqrt{n}}\tilde{r}_{n}^{\frac{\gamma}{\gamma+1}}}}
\end{eqnarray*}
valid for all $\gamma=0,...,\beta$. In the second inequality, we
have used (\ref{eq:78}). Therefore, we can take $\min_{\gamma\in\left\{ 0,...,\beta\right\} }\mathcal{T}_{2}\left(\gamma,\tilde{r}_{n}\right)=\mathcal{B}_{2}\left(\tilde{r}_{n}\right)$
and let $\gamma_{2}^{*}$ be the minimizer. Setting $\mathcal{B}_{2}\left(\tilde{r}_{n}\right)\asymp\frac{\tilde{r}_{n}^{2}}{\sigma}$
yields $\mathcal{M}_{2}\left(\gamma_{2}^{*}\right)$.

Now, we can take $\tilde{r}_{n}^{2}=\min\left\{ \mathcal{M}_{1}\left(\gamma_{1}^{*}\right),\,\mathcal{M}_{2}\left(\gamma_{2}^{*}\right)\right\} $.
By Theorem 13.5 in \cite{wainwright 2019}, in terms of (\ref{eq:nonpara ls}),
we have 
\[
\frac{1}{n}\sum_{i=1}^{n}\left(\hat{y}(x_{i})-y(x_{i})\right)^{2}\precsim\tilde{r}_{n}^{2}
\]
with probability at least $1-\exp\left(\frac{-n\tilde{r}_{n}^{2}}{2\sigma^{2}}\right)$.
Integrating the tail bound yields 
\[
\mathbb{E}\left[\frac{1}{n}\sum_{i=1}^{n}\left(\hat{y}(x_{i})-y(x_{i})\right)^{2}\right]\precsim\tilde{r}_{n}^{2}+c_{1}^{'}\exp\left\{ -c_{2}^{'}n\sigma^{-2}\tilde{r}_{n}^{2}\right\} .
\]

\textbf{Claim (i)}. Following almost identical arguments as above,
we have, for the autonomous ODEs, 
\[
\frac{1}{n}\sum_{i=1}^{n}\left(\hat{y}(x_{i})-y(x_{i})\right)^{2}\precsim\tilde{r}_{n,a}^{2}
\]
with probability at least $1-\exp\left(\frac{-n\tilde{r}_{n,a}^{2}}{2\sigma^{2}}\right)$.
Integrating the tail bound yields 
\[
\mathbb{E}\left[\frac{1}{n}\sum_{i=1}^{n}\left(\hat{y}(x_{i})-y(x_{i})\right)^{2}\right]\precsim\tilde{r}_{n,a}^{2}+c_{1}^{'}\exp\left\{ -c_{2}^{'}n\sigma^{-2}\tilde{r}_{n,a}^{2}\right\} .
\]

\textbf{Claim (iii)}. Note that if $\frac{n}{\sigma^{2}}\succsim\left(\gamma\sqrt{\log\left(\gamma\vee1\right)}\right)^{4\left(\gamma+2\right)+2}$
where $\gamma\in\left\{ 0,...,\beta\right\} $, we have 
\begin{eqnarray}
\mathcal{M}_{1}^{a}\left(\gamma\right)\asymp\left(\frac{\sigma^{2}}{n}\right)^{\frac{2\left(\gamma+2\right)}{2\left(\gamma+2\right)+1}}, & \textrm{and } & \tilde{r}_{n,a}^{2}=\mathcal{M}_{1}^{a}\left(\gamma\right),\nonumber \\
\mathcal{M}_{1}\left(\gamma\right)\asymp\left(\frac{\sigma^{2}}{n}\right)^{\frac{2\left(\gamma+2\right)}{2\left(\gamma+2\right)+1}}, & \textrm{and} & \tilde{r}_{n}^{2}=\mathcal{M}_{1}\left(\gamma\right).\label{eq:optimal}
\end{eqnarray}

\textbf{Claim (iv)}. Note that every member in $\mathcal{S}_{\beta+2}$
(the standard smooth class of degree $\beta+2$) can be expressed
as a solution to the ODE (\ref{eq:separate}) and clearly, $\mathcal{S}_{\beta+2}\subseteq\mathcal{Y}$
where $\mathcal{Y}$ is the class of solutions to (\ref{eq:11-4-1})
with $f\in\mathcal{S}_{\beta+1,\,2}\left(1,\,\Xi\right)$. As we have
discussed in Section 3.1, by the existing minimax results on $\mathcal{S}_{\beta+2}$
(see, e.g., \cite{wainwright 2019}, Chapter 15), the minimax lower
bound has a scaling $\left(\frac{\sigma^{2}}{n}\right)^{\frac{2\left(\beta+2\right)}{2\left(\beta+2\right)+1}}$.
When $\frac{n}{\sigma^{2}}\succsim\left(\beta\sqrt{\log\left(\beta\vee1\right)}\right)^{4\left(\beta+2\right)+2}$,
we have (\ref{eq:optimal}) with $\gamma=\beta$ and the rate $\left(\frac{\sigma^{2}}{n}\right)^{\frac{2\left(\beta+2\right)}{2\left(\beta+2\right)+1}}$
is clearly minimax optimal for (\ref{eq:11-4-1}) with $f$ ranging
over $\mathcal{S}_{\beta+1,\,2}\left(1,\,\Xi\right)$. $\square$

\subsection{Theorem 3.4\label{subsec:Theorem-3.4}}

\subsubsection*{Preliminary }

An alternative approach for bounding $\mathcal{G}_{n}(\tilde{r}_{n};\,\bar{\mathcal{F}})$
is developed by \cite{Mendelson 2002-1}. This approach can be formulated
as follows. 

Given a radius $\tilde{r}_{n}>0$ and a function class $\bar{\mathcal{F}}$,
define \textit{
\begin{equation}
\mathcal{G}_{n}(\tilde{r}_{n};\,\bar{\mathcal{F}}):=\mathbb{E}_{\varepsilon}\left[\sup_{h\in\Lambda(\tilde{r}_{n};\,\bar{\mathcal{F}})}\left|\frac{1}{n}\sum_{i=1}^{n}\varepsilon_{i}h(x_{i})\right|\right],\label{eq:13-1-1}
\end{equation}
}where $\varepsilon_{i}\overset{i.i.d.}{\sim}\mathcal{N}\left(0,\,1\right)$,
$\varepsilon=\left\{ \varepsilon_{i}\right\} _{i=1}^{n}$, and 
\[
\Lambda(\tilde{r}_{n};\,\bar{\mathcal{F}})=\left\{ h\in\bar{\mathcal{F}}:\,\left|h\right|_{n}\leq\tilde{r}_{n},\,\left|h\right|_{\mathcal{H}}\leq1\right\} 
\]
where $\bar{\mathcal{F}}$ is defined in (\ref{eq:28}) and $\left|\cdot\right|_{\mathcal{H}}$
is a norm associated with some underlying RKHS. Let $\mathcal{K}$
be the kernel function of this RKHS and $\mu_{1}\geq\mu_{2}\geq\cdots\geq\mu_{n}\geq0$
be the eigenvalues of the matrix $\mathbb{K}$ consisting of entries
$K_{ij}=\frac{1}{n}\mathcal{K}\left(x_{i},\,x_{j}\right)$. Introduce
the critical radius condition 
\begin{equation}
\mathcal{G}_{n}\left(\tilde{r}_{n};\,\bar{\mathcal{F}}\right)\precsim\mathcal{R}\frac{\tilde{r}_{n}^{2}}{\sigma},
\end{equation}
where $\mathcal{R}$ is the radius of the RKHS of interest. For any
$\tilde{r}_{n}>0$, one has 
\begin{equation}
\mathcal{G}_{n}(\tilde{r}_{n};\,\bar{\mathcal{F}})\precsim\sqrt{\frac{1}{n}}\sqrt{\sum_{i=1}^{n}\left(\tilde{r}_{n}^{2}\wedge\mu_{i}\right)}\label{eq:critical-2}
\end{equation}
(see \cite{Mendelson 2002-1}). As a consequence, the critical radius
condition (\ref{eq:critical-2}) is satisfied for any $\tilde{r}_{n}>0$
such that 
\[
\sqrt{\frac{1}{n}}\sqrt{\sum_{i=1}^{n}\left(\tilde{r}_{n}^{2}\wedge\mu_{i}\right)}\precsim\frac{\mathcal{R}\tilde{r}_{n}^{2}}{\sigma}.
\]
Therefore, in terms of
\[
\hat{g}\in\arg\min_{\tilde{g}\in\mathcal{F},\,\left|\tilde{g}\right|_{\mathcal{H}}\leq\mathcal{R}}\frac{1}{2n}\sum_{i=1}^{n}\left(Y_{i}-\tilde{g}\left(x_{i}\right)\right)^{2},
\]
one has 
\[
\frac{1}{n}\sum_{i=1}^{n}\left(\hat{g}(x_{i})-g^{*}(x_{i})\right)^{2}\precsim\mathcal{R}^{2}\tilde{r}_{n}^{2}
\]
with probability at least $1-c^{'}\exp\left(\frac{-c^{''}n\mathcal{R}^{2}\tilde{r}_{n}^{2}}{\sigma^{2}}\right)$.\\
\textbf{}\\
\textbf{Proof}. In what follows, we argue for the nonautonomous ODEs
as the argument for the autonomous ODEs is nearly identical. 

Given that $\mathcal{Y}\subseteq\mathcal{S}_{\beta+2,\overline{C}}^{\dagger}$,
$\mathcal{Y}$ is contained in $\mathcal{H}_{\beta+2}^{\dagger}$,
defined in (\ref{eq:RKHS2}). Clearly, we have 
\begin{align*}
\mathcal{H}_{\beta+2}^{\dagger} & \subseteq\mathcal{H}_{\gamma+2}:=\{h:\,\left[0,\,1\right]\rightarrow\mathbb{R}\vert h^{(\gamma+1)}\textrm{ is abs. cont. with }\\
 & \int_{0}^{1}\left[h^{(\gamma+2)}\left(t\right)\right]^{2}dt\leq\left[2^{\gamma+1}\left(\gamma+1\right)!\right]^{2}\}
\end{align*}
for all $\gamma\in\left\{ 0,...,\beta\right\} $, and 
\begin{equation}
\mathcal{H}_{\gamma+2}=\mathcal{H}_{\gamma+2,1}+\mathcal{H}_{\gamma+2,2}\label{eq:decomp}
\end{equation}
where $\mathcal{H}_{\gamma+2,1}$ is an RKHS of polynomials of degree
$\gamma+1$ and 
\begin{align*}
\mathcal{H}_{\gamma+2,2} & :=\{h:\,\left[0,\,1\right]\rightarrow\mathbb{R}\vert h^{(k)}\left(0\right)=0\textrm{ for all }k\leq\gamma,\\
 & h^{(\gamma+1)}\textrm{ is abs. cont. with }\int_{0}^{1}\left[h^{(\gamma+2)}\left(t\right)\right]^{2}dt\leq\left[2^{\gamma+1}\left(\gamma+1\right)!\right]^{2}\}.
\end{align*}
As a result, we can equip $\mathcal{H}_{\gamma+2}$ with the following
norm 
\[
\left|h\right|_{\mathcal{H}}=\left[\sum_{k=0}^{\gamma+1}\left(h^{(k)}\left(0\right)\right)^{2}+\int_{0}^{1}\left[h^{(\gamma+2)}\left(t\right)\right]^{2}dt\right]^{\frac{1}{2}}.
\]
Note that any $h\in\mathcal{Y}$ has 
\[
\left|h\right|_{\mathcal{H}}\leq c\left[2^{\gamma+1}\left(\gamma+1\right)!\right]^{2}
\]
for some positive universal constant $c$ so we may let $\mathcal{R}_{\gamma}:=c\left[2^{\gamma+1}\left(\gamma+1\right)!\right]^{2}$
be the radius associated with $\mathcal{H}_{\gamma+2}$. Because of
(\ref{eq:decomp}), we can generate $\mathcal{H}_{\gamma+2}$ with
the kernel 
\[
\mathcal{K}\left(x,\,x^{'}\right)=\sum_{k=0}^{\gamma+1}\frac{x^{k}}{k!}\frac{x^{'k}}{k!}+\int_{0}^{1}\frac{\left(x-t\right)_{+}^{\gamma+1}}{\left(\gamma+1\right)!}\frac{\left(x^{'}-t\right)_{+}^{\gamma+1}}{\left(\gamma+1\right)!}dt
\]
where $\left(w\right)_{+}=w\vee0$. As a consequence, we have
\[
\mathcal{G}_{n}(\tilde{r}_{n};\,\bar{\mathcal{F}})\precsim\underset{\mathcal{T}_{3}\left(\gamma,\tilde{r}_{n}\right)}{\underbrace{\max\left\{ \tilde{r}_{n}\sqrt{\frac{\left(\gamma\wedge n\right)\vee1}{n}},\,\frac{1}{\sqrt{n}}\tilde{r}_{n}^{\frac{2\gamma+3}{2\gamma+4}}\right\} }}
\]
valid for all $\gamma\in\left\{ 0,...,\beta\right\} $. Setting $\mathcal{T}_{3}\left(\gamma,\tilde{r}_{n}\right)\asymp\left[2^{\gamma+1}\left(\gamma+1\right)!\right]\tilde{r}_{n}^{2}\sigma^{-1}$
yields 
\[
\mathcal{R}_{\gamma}^{2}\tilde{r}_{n}^{2}\asymp\max\left\{ \frac{\sigma^{2}\left(\left(\gamma\wedge n\right)\vee1\right)}{n},\,\left[2^{\gamma+1}\left(\gamma+1\right)!\right]^{\frac{2}{2\left(\gamma+2\right)+1}}\left(\frac{\sigma^{2}}{n}\right)^{\frac{2\left(\gamma+2\right)}{2\left(\gamma+2\right)+1}}\right\} 
\]
valid for all $\gamma\in\left\{ 0,...,\beta\right\} $. Taking the
minimum of the RHS over all $\gamma\in\left\{ 0,...,\beta\right\} $
yields $\bar{r}_{n}^{2}$ in Theorem 3.4. Hence, 
\[
\frac{1}{n}\sum_{i=1}^{n}\left(\hat{y}(x_{i})-y^{*}(x_{i})\right)^{2}\precsim\bar{r}_{n}^{2}
\]
with probability at least $1-c^{'}\exp\left(\frac{-c^{''}n\bar{r}_{n}^{2}}{\sigma^{2}}\right)$.
Integrating the tail bound yields the claim in Theorem 3.4. $\square$

\section{Supporting lemmas and proofs\label{sec:Supporting-lemmas-and} }

\textbf{Theorem B.1} (Gronwall inequality for first order ODEs). \textit{Consider
the following pair of ODEs: 
\[
y^{'}\left(x\right)=f\left(x,\,y\left(x\right)\right),\quad y\left(0\right)=y_{0},
\]
and 
\[
z^{'}\left(x\right)=g\left(x,\,z\left(x\right)\right),\quad z\left(0\right)=z_{0},
\]
with $\left|y_{0}\right|,\,\left|z_{0}\right|\leq C_{0}$, and $\left(x,\,y\left(x\right)\right),\,\left(x,\,z\left(x\right)\right)\in\bar{\Gamma}:=\left[0,\,1\right]\times\left[-C_{0}-b,\,C_{0}+b\right]$.
Suppose $f$ and $g$ are continuous on $\bar{\Gamma}$; for all $\left(x,\,y\right),\,\left(x,\,\tilde{y}\right)\in\bar{\Gamma}$,
\begin{equation}
\left|f\left(x,\,y\right)-f\left(x,\,\tilde{y}\right)\right|\leq L\left|y-\tilde{y}\right|.\label{eq:22-2-1}
\end{equation}
Assume there is a continuous function $\varphi:\,\left[0,\,a\right]\mapsto[0,\,\infty)$
such that }
\begin{equation}
\left|f\left(x,\,y\left(x\right)\right)-g\left(x,\,y\left(x\right)\right)\right|\leq\varphi\left(x\right).\label{eq:23-3}
\end{equation}
\textit{Then we have}
\[
\left|y\left(x\right)-z\left(x\right)\right|\leq\exp\left(Lx\right)\int_{0}^{x}\exp\left(-Ls\right)\varphi\left(s\right)ds+\exp\left(Lx\right)\left|y_{0}-z_{0}\right|
\]
\textit{for $x\in\left[0,\,\min\left\{ a,\,\frac{b}{M}\right\} \right]$
where $M=\max\left\{ \max_{\left(x,\,y\right)\in\bar{\Gamma}}\left|f\left(x,\,y\right)\right|,\,\max_{\left(x,\,z\right)\in\bar{\Gamma}}\left|g\left(x,\,z\right)\right|\right\} $.
}\\
\textbf{}\\
\textbf{Remark}. Theorem B.1 is a slight modification of Theorem 2.1
in \cite{Howard 1998}, which gives a variant of the Gronwall inequality
\cite{Gronwall 1919}.\\
\\
In the following result, we extend Theorem B.1 to higher order ODEs.
Let 
\begin{eqnarray*}
Y\left(x\right)=\left[\begin{array}{c}
y\left(x\right)\\
y^{'}\left(x\right)\\
\vdots\\
y^{\left(m-1\right)}\left(x\right)
\end{array}\right] & \textrm{and} & Z\left(x\right)=\left[\begin{array}{c}
z\left(x\right)\\
z^{'}\left(x\right)\\
\vdots\\
z^{\left(m-1\right)}\left(x\right)
\end{array}\right]
\end{eqnarray*}
with
\begin{eqnarray*}
Y_{0}:=\left[\begin{array}{c}
y_{\left(0\right)}\\
y_{\left(1\right)}\\
\vdots\\
y_{\left(m-1\right)}
\end{array}\right] & \textrm{and} & Z_{0}:=\left[\begin{array}{c}
z_{\left(0\right)}\\
z_{\left(1\right)}\\
\vdots\\
z_{\left(m-1\right)}
\end{array}\right].
\end{eqnarray*}
In addition, let 
\[
\bar{\Gamma}:=\left\{ \left(x,\,Y\right):\,x\in\left[a_{0},\,a_{0}+a\right],\,\left|Y\right|_{2}\leq b+C_{0}\right\} .
\]
\textbf{Theorem B.2} (Gronwall inequality for higher order ODEs).\textit{
Consider the following pair of ODEs:
\begin{align}
y^{\left(m\right)}\left(x\right)=f\left(x,\,y\left(x\right),\,y^{'}\left(x\right),\,...,y^{\left(m-1\right)}\left(x\right)\right),\nonumber \\
y\left(a_{0}\right)=y_{(0)},\,y^{'}\left(a_{0}\right)=y_{(1)},\,...,\,y^{\left(m-1\right)}\left(a_{0}\right)=y_{(m-1)},\label{eq:22-3}
\end{align}
and 
\begin{align}
z^{\left(m\right)}\left(x\right)=g\left(x,\,z\left(x\right),\,z^{'}\left(x\right),\,...,z^{\left(m-1\right)}\left(x\right)\right),\nonumber \\
z\left(a_{0}\right)=z_{(0)},\,z^{'}\left(a_{0}\right)=z_{(1)},\,...,\,z^{\left(m-1\right)}\left(a_{0}\right)=z_{(m-1)},\label{eq:23-1}
\end{align}
with $\left|Y_{0}\right|_{2},\,\left|Z_{0}\right|_{2}\leq C_{0}$
and $\left(x,\,Y\left(x\right)\right),\,\left(x,\,Z\left(x\right)\right)\in\bar{\Gamma}$.
Suppose $f$ and $g$ are continuous on $\bar{\Gamma}$; and
\begin{equation}
\left|f\left(x,\,Y\right)-f\left(x,\,\tilde{Y}\right)\right|\leq L\left|Y-\tilde{Y}\right|_{2}\label{eq:22-2}
\end{equation}
for all $\left(x,\,Y\right):=\left(x,\,y,\,...,\,y^{\left(m-1\right)}\right)$
and $\left(x,\,\tilde{Y}\right):=\left(x,\,\tilde{y},\,...,\,\tilde{y}^{\left(m-1\right)}\right)$
in $\bar{\Gamma}$. Assume there is a continuous function $\varphi:\,\left[a_{0},\,a_{0}+a\right]\mapsto[0,\,\infty)$
such that }
\begin{equation}
\left|f\left(x,\,Y\left(x\right)\right)-g\left(x,\,Y\left(x\right)\right)\right|\leq\varphi\left(x\right).\label{eq:23}
\end{equation}
\textit{Then we have}
\[
\left|y^{(k)}\left(x\right)-z^{\left(k\right)}\left(x\right)\right|\leq\exp\left(x\sqrt{L^{2}+1}\right)\int_{a_{0}}^{x}\exp\left(-s\sqrt{L^{2}+1}\right)\varphi\left(s\right)ds+\exp\left(x\sqrt{L^{2}+1}\right)\left|Y_{0}-Z_{0}\right|_{2}
\]
\textit{for all $k=0,...,m-1$ and all $x\in\left[a_{0},\,a_{0}+\min\left\{ a,\,\frac{b}{M}\right\} \right]$
where 
\[
M=\max\left\{ \max_{\left(x,\,Y\right)\in\bar{\Gamma}}\left|f\left(x,\,Y\right)\right|,\,\max_{\left(x,\,Z\right)\in\bar{\Gamma}}\left|g\left(x,\,Z\right)\right|\right\} .
\]
}\textbf{Proof}. The following arguments are based on and extend those
in \cite{Howard 1998}. Let $W=\left(w_{j}\right)_{j=0}^{m-1}$ and
\[
F\left(x,\,W\right):=\left[\begin{array}{c}
w_{1}\\
w_{2}\\
\vdots\\
w_{m-1}\\
f\left(x,\,W\right)
\end{array}\right].
\]
We transform (\ref{eq:22-3}) and (\ref{eq:23-1}) into 
\begin{eqnarray*}
Y^{'}\left(x\right) & = & F\left(x,\,Y\left(x\right)\right),\\
Z^{'}\left(x\right) & = & G\left(x,\,Z\left(x\right)\right),
\end{eqnarray*}
where 
\[
F\left(x,\,Y\left(x\right)\right):=\left[\begin{array}{c}
y^{'}\left(x\right)\\
y^{\left(2\right)}\left(x\right)\\
\vdots\\
y^{\left(m-1\right)}\left(x\right)\\
f\left(x,\,y\left(x\right),\,y^{'}\left(x\right),\,...,y^{\left(m-1\right)}\left(x\right)\right)
\end{array}\right]
\]
and 
\[
G\left(x,\,Z\left(x\right)\right):=\left[\begin{array}{c}
z^{'}\left(x\right)\\
z^{\left(2\right)}\left(x\right)\\
\vdots\\
z^{\left(m-1\right)}\left(x\right)\\
g\left(x,\,z\left(x\right),\,z^{'}\left(x\right),\,...,z^{\left(m-1\right)}\left(x\right)\right)
\end{array}\right].
\]
The inequality $\frac{d}{ds}\left|Y\left(s\right)\right|_{2}\leq\left|Y^{'}\left(s\right)\right|_{2}$,
(\ref{eq:22-2}) and (\ref{eq:23}) yield 

\begin{eqnarray*}
\frac{d}{ds}\left|Y\left(s\right)-Z\left(s\right)\right|_{2} & \leq & \left|Y^{'}\left(s\right)-Z^{'}\left(s\right)\right|_{2}\\
 & = & \left|F\left(s,\,Y\left(s\right)\right)-G\left(s,\,Z\left(s\right)\right)\right|_{2}\\
 & \leq & \left|G\left(s,\,Z\left(s\right)\right)-F\left(s,\,Z\left(s\right)\right)\right|_{2}+\left|F\left(s,\,Z\left(s\right)\right)-F\left(s,\,Y\left(s\right)\right)\right|_{2}\\
 & \leq & \varphi\left(s\right)+\left[L^{2}\left|Y\left(s\right)-Z\left(s\right)\right|_{2}^{2}+\sum_{k=1}^{m-1}\left(y^{\left(k\right)}\left(s\right)-z^{\left(k\right)}\left(s\right)\right)^{2}\right]^{\frac{1}{2}}\\
 & \leq & \varphi\left(s\right)+\sqrt{L^{2}+1}\left|Y\left(s\right)-Z\left(s\right)\right|_{2},
\end{eqnarray*}
which is equivalent to 
\[
\frac{d}{ds}\left|Y\left(s\right)-Z\left(s\right)\right|_{2}-\sqrt{L^{2}+1}\left|Y\left(s\right)-Z\left(s\right)\right|_{2}\leq\varphi\left(s\right).
\]
Multiplying both sides above by $\exp\left(-s\sqrt{L^{2}+1}\right)$
gives
\[
\frac{d}{ds}\left[\exp\left(-s\sqrt{L^{2}+1}\right)\left|Y\left(s\right)-Z\left(s\right)\right|_{2}\right]\leq\exp\left(-s\sqrt{L^{2}+1}\right)\varphi\left(s\right).
\]
By the Cauchy-Peano existence theorem, there exist solutions $y$
and $z$ to (\ref{eq:22-3}) and (\ref{eq:23-1}), respectively, on
$\left[a_{0},\,a_{0}+\min\left\{ a,\,\frac{b}{M}\right\} \right]$.
Then integrating both sides of the inequality above from $a_{0}$
to $x$ gives 
\[
\exp\left(-x\sqrt{L^{2}+1}\right)\left|Y\left(x\right)-Z\left(x\right)\right|_{2}-\left|Y_{0}-Z_{0}\right|_{2}\leq\int_{a_{0}}^{x}\exp\left(-s\sqrt{L^{2}+1}\right)\varphi\left(s\right)ds
\]
for all $x\in\left[a_{0},\,a_{0}+\min\left\{ a,\,\frac{b}{M}\right\} \right]$.
The above implies that
\[
\left|y^{(k)}\left(x\right)-z^{\left(k\right)}\left(x\right)\right|\leq\exp\left(x\sqrt{L^{2}+1}\right)\int_{a_{0}}^{x}\exp\left(-s\sqrt{L^{2}+1}\right)\varphi\left(s\right)ds+\exp\left(x\sqrt{L^{2}+1}\right)\left|Y_{0}-Z_{0}\right|_{2}
\]
for all $x\in\left[a_{0},\,a_{0}+\min\left\{ a,\,\frac{b}{M}\right\} \right]$
and all $k=0,...,m-1$. $\square$

\section{Additional results and proofs}

\subsection{Proposition C.1\label{subsec:Proposition-C.1}}

\textbf{Proposition C.1.}\textit{ Suppose $f$ in (\ref{eq:11-2-1})
is continuous on $\left[0,\,1\right]\times\left[-C_{0}-b,\,C_{0}+b\right]$,
$\left|f\left(x,\,y;\,\theta\right)\right|\leq1$ and
\begin{equation}
\left|f\left(x,\,y;\,\theta\right)-f\left(x,\,\tilde{y};\,\theta\right)\right|\leq\left|y-\tilde{y}\right|\label{eq:22}
\end{equation}
for all $\left(x,\,y\right)$ and $\left(x,\,\tilde{y}\right)$ in
$\left[0,\,1\right]\times\left[-C_{0}-b,\,C_{0}+b\right]$, and $\theta\in\mathbb{B}_{q}\left(1\right)$
with $q\in\left[1,\,\infty\right]$; moreover, 
\begin{equation}
\left|f\left(x,\,y;\,\theta\right)-f\left(x,\,y;\,\theta^{'}\right)\right|\leq L_{K}\left|\theta-\theta^{'}\right|_{q},\label{eq:17-1}
\end{equation}
for all $\left(x,\,y\right)\in\left[0,\,1\right]\times\left[-C_{0}-b,\,C_{0}+b\right]$
and $\theta,\,\theta^{'}\in\mathbb{B}_{q}\left(1\right)$. Let us
consider (\ref{eq:model}) where $y\left(\cdot\right)$ is the (unique)
solution to (\ref{eq:11-2-1}) on $[0,\,\alpha)$ where }$\alpha=\min\left\{ 1,\,b\right\} $\textit{.
Suppose we have $n$ design points sampled from the interval }$[0,\,\tilde{\alpha}]$\textit{
with $\tilde{\alpha}<\alpha$. Letting $B_{K}=\left(L_{K}\vee1\right)$,
if 
\begin{eqnarray}
K\log\left(1+\frac{2L_{\max}L_{K}}{\delta}\right) & \succsim & \log\left(\frac{2C_{0}L_{\max}}{\delta}+1\right),\:\forall\delta\precsim B_{K}\sigma\sqrt{\frac{K}{n}}\label{eq:10-1}\\
\max\left\{ \alpha,\,C_{0}\right\}  & \geq & c_{0}B_{K}\sigma\sqrt{\frac{K}{n}},\label{eq:9-1}
\end{eqnarray}
for a sufficiently large positive universal constant $c_{0}$, then
we have
\begin{equation}
\frac{1}{n}\sum_{i=1}^{n}\left(y^{*}(x_{i};\,\hat{\theta},\,\hat{y}_{0})-y^{*}(x_{i};\,\theta^{*},\,y_{0}^{*})\right)^{2}\precsim B_{K}^{2}\frac{\sigma^{2}K}{n}\label{eq:10}
\end{equation}
with probability at least $1-c_{1}\exp\left(-c_{2}B_{K}^{2}K\right)$,
where $\left(\hat{\theta},\,\hat{y}_{0}\right)$ is a solution to
(\ref{eq:20-1}).}\\
\textit{}\\
\textbf{Remark}. Condition (\ref{eq:10-1}) simply restricts $C_{0}$
from being too large, and as a consequence, (\ref{eq:16}) implies
that $\log N_{\infty}\left(\delta,\,\mathcal{Y}\right)\precsim K\log\left(1+\frac{2L_{\max}L_{K}}{\delta}\right)$.
This upper bound implies that $\mathcal{Y}$ is no ``larger'' than
the class of $f$s parameterized by $\theta\in\mathbb{B}_{q}\left(1\right)$.
\\
\\
\textbf{Remark}\textit{. }Condition (\ref{eq:9-1}) in Proposition
C.1 simply excludes the trivial case where $b$ and $C_{0}$ are ``too
small''. Without such a condition, as long as $f$ is bounded from
above, we would simply replace (\ref{eq:10}) with 
\[
\left[\frac{1}{n}\sum_{i=1}^{n}\left(y(x_{i};\,\hat{\theta},\,\hat{y}_{0})-y(x_{i};\,\theta^{*},\,y_{0})\right)^{2}\right]^{\frac{1}{2}}\leq c_{1}\min\left\{ B_{K}\sigma\sqrt{\frac{K}{n}},\,\max\left\{ \alpha,\,C_{0}\right\} \right\} .
\]
\textbf{Proof}. Let $\mathcal{F}=\mathcal{Y}$ in (\ref{eq:28}).
We have 

\begin{eqnarray}
\frac{1}{\sqrt{n}}\int_{0}^{\tilde{r}_{n}}\sqrt{\log N_{n}(\delta,\,\Lambda(\tilde{r}_{n};\,\bar{\mathcal{F}}))}d\delta & \leq & \frac{1}{\sqrt{n}}\int_{0}^{\tilde{r}_{n}}\sqrt{\log N_{\infty}(\delta,\,\Lambda(\tilde{r}_{n};\,\bar{\mathcal{F}}))}d\delta\nonumber \\
 & \leq & \sqrt{\frac{K}{n}}\int_{0}^{\tilde{r}_{n}}\sqrt{\log\left(1+\frac{c\tilde{r}_{n}\left(L_{K}\vee1\right)}{\delta}\right)}d\delta\nonumber \\
 & = & \tilde{r}_{n}\left(L_{K}\vee1\right)\sqrt{\frac{K}{n}}\int_{0}^{\frac{1}{L_{K}\vee1}}\sqrt{\log\left(1+\frac{c}{t}\right)}dt\nonumber \\
 & \leq & \left(L_{K}\vee1\right)\tilde{r}_{n}\sqrt{\frac{K}{n}}\label{eq:24}
\end{eqnarray}
where we have applied a change of variable $t=\frac{\delta}{\tilde{r}_{n}\left(L_{K}\vee1\right)}$
in the third line. Setting $B_{K}\tilde{r}_{n}\sqrt{\frac{K}{n}}\asymp\frac{\tilde{r}_{n}^{2}}{\sigma}$
yields $\tilde{r}_{n}\asymp B_{K}\sigma\sqrt{\frac{K}{n}}$, where
$B_{K}=\left(L_{K}\vee1\right)$. By Theorem 13.5 in \cite{wainwright 2019},
we obtain

\[
\frac{1}{n}\sum_{i=1}^{n}\left(y(x_{i};\,\hat{\theta},\,\hat{y}_{0})-y(x_{i};\,\theta^{*},\,y_{0})\right)^{2}\precsim B_{K}^{2}\frac{\sigma^{2}K}{n}.
\]
with probability at least $1-c_{1}\exp\left\{ -c_{2}B_{K}^{2}K\right\} $.
$\square$

\subsection{Proposition C.2\label{subsec:Proposition-C.2}}

\textbf{Proposition C.2.}\textit{ Suppose the conditions in Proposition
C.1 hold. In terms of $\hat{y}_{R+1}\left(x_{i};\,\hat{\theta},\,\hat{y}_{0}\right)$
where $\hat{\theta}$ is obtained from solving (\ref{eq:20}), if
\begin{equation}
\max\left\{ \tilde{\alpha},\,C_{0}\right\} \geq c_{0}\sigma\tilde{b}\sqrt{\frac{K}{n}}\label{eq:13}
\end{equation}
for a sufficiently large positive universal constant $c_{0}$, then
we have 
\begin{equation}
\left\{ \frac{1}{n}\sum_{i=1}^{n}\left[\hat{y}_{R+1}(x_{i};\,\hat{\theta},\,\hat{y}_{0})-y(x_{i};\,\theta^{*},\,y_{0})\right]^{2}\right\} ^{\frac{1}{2}}\precsim\sigma\left(\tilde{b}\sqrt{\frac{K}{n}}+\frac{1}{1-\tilde{\alpha}}\left|\hat{y}_{0}-y_{0}\right|+\frac{\tilde{\alpha}^{R+1}}{1-\tilde{\alpha}}\max\left\{ C_{0},\,\tilde{\alpha}\right\} \right)\label{eq:22-4}
\end{equation}
with probability at least $1-c_{1}\exp\left(-c_{2}n\right)-c_{3}\exp\left(-c_{4}K\tilde{b}^{2}\right)$,
where }$\tilde{b}=\left(\frac{\tilde{\alpha}L_{K}}{1-\tilde{\alpha}}\vee1\right)$\textit{.}\\
\textit{}\\
\textbf{Remark}. If the integral in (\ref{eq:26-1}) is hard to compute
analytically, numerical integration can be used in (\ref{eq:26}).
This would introduce an additional approximation error $\sigma\left(R+1\right)\cdot Err$,
where $Err$ is an upper bound on the error incurred in each iteration
of (\ref{eq:26}), depending on the smoothness of $f$ and which numerical
method is used. For example, if $f$ is twice differentiable with
bounded first and second derivatives, and the integral is approximated
with the midpoint rule with $T$ slices, then $Err\precsim T^{-2}$.
\\
\\
\textbf{Remark}. The bound (\ref{eq:22-4}) reflects three sources
of errors: $\tilde{b}\sqrt{\frac{K}{n}}$ is due to the estimation
error in $\hat{\theta}$, $\frac{1}{1-\tilde{\alpha}}\left|\hat{y}_{0}-y_{0}\right|$
is due to the estimation error in $\hat{y}_{0}$, and $\frac{\tilde{\alpha}^{R+1}}{1-\tilde{\alpha}}\max\left\{ C_{0},\,\tilde{\alpha}\right\} $
is due to the error from the finite ($R+1$) Picard iterations. \\
\textbf{}\\
\textbf{Proof}. In what follows, we suppress the dependence of $\hat{y}_{R+1}(x_{i};\,\hat{\theta},\,\hat{y}_{0})$
on $\hat{y}_{0}$ ($y(x_{i};\,\theta^{*},\,y_{0})$ on $y_{0}$) and
simply write $\hat{y}_{R+1}(x_{i},\hat{\theta})$ (respectively, $y(x_{i},\theta^{*})$).
By (\ref{eq:44}), we need to bound 
\[
\left|\frac{1}{n}\sum_{i=1}^{n}\varepsilon_{i}\left(\hat{y}_{R+1}(x_{i},\hat{\theta})-y(x_{i},\theta^{*})\right)\right|.
\]
We can write $\hat{y}_{R+1}(x_{i},\hat{\theta})-y(x_{i},\theta^{*})=\sum_{j=1}^{3}T_{j}\left(x_{i}\right)$
where
\begin{eqnarray*}
T_{1}\left(x_{i}\right) & = & \hat{y}_{R+1}(x_{i};\,\hat{\theta})-\hat{y}_{R+1}(x_{i};\,\theta^{*}),\quad\textrm{estimation error due to }\hat{\theta}\\
T_{2}\left(x_{i}\right) & = & \hat{y}_{R+1}(x_{i};\,\theta^{*})-y_{R+1}(x_{i};\,\theta^{*}),\quad\textrm{estimation error due to }\hat{y}_{0}\\
T_{3}\left(x_{i}\right) & = & y_{R+1}(x_{i};\,\theta^{*})-y(x_{i};\,\theta^{*}),\quad\textrm{estimation error due to the finite iterations},
\end{eqnarray*}
where 
\begin{eqnarray*}
y_{0} & = & y_{0},\\
y_{1}\left(x;\,\theta\right) & = & y_{0}^{*}+\int_{0}^{x}f\left(s,\,y_{0};\theta\right)ds,\\
y_{2}\left(x;\,\theta\right) & = & y_{0}^{*}+\int_{0}^{x}f\left(s,\,y_{1}\left(s;\,\theta\right);\theta\right)ds,\\
 & \vdots\\
y_{R+1}\left(x;\,\theta\right) & = & y_{0}+\int_{0}^{x}f\left(s,\,y_{R}\left(s;\,\theta\right);\theta\right)ds.
\end{eqnarray*}
As a result, we have 
\begin{equation}
\left|\frac{1}{n}\sum_{i=1}^{n}\varepsilon_{i}\left(\hat{y}_{R+1}(x_{i},\hat{\theta})-y(x_{i},\theta^{*})\right)\right|\precsim\left|\frac{1}{n}\sum_{i=1}^{n}\varepsilon_{i}T_{1}\left(x_{i}\right)\right|+\sigma\sqrt{\frac{1}{n}\sum_{i=1}^{n}\left[T_{2}\left(x_{i}\right)\right]^{2}}+\sigma\sqrt{\frac{1}{n}\sum_{i=1}^{n}\left[T_{3}\left(x_{i}\right)\right]^{2}}\label{eq:31-1}
\end{equation}
with probability at least $1-c_{1}\exp\left(-c_{2}n\right)$. 

We first analyze $\left|\frac{1}{n}\sum_{i=1}^{n}\varepsilon_{i}T_{1}\left(x_{i}\right)\right|$
in (\ref{eq:31-1}). In (\ref{eq:28}), let 
\begin{equation}
\mathcal{F}=\left\{ g_{\theta}\left(x\right)=\hat{y}_{R+1}\left(x;\,\theta\right):\,\theta\in\mathbb{B}_{q}\left(1\right),\,x\in\left[0,\,\tilde{\alpha}\right]\right\} \label{eq:41}
\end{equation}
where $\hat{y}_{R+1}\left(s;\,\theta\right)$ is constructed in the
following fashion: 
\begin{eqnarray*}
\hat{y}_{0} & = & \hat{y}_{0},\\
\hat{y}_{1}\left(x;\,\theta\right) & = & \hat{y}_{0}+\int_{0}^{x}f\left(s,\,\hat{y}_{0};\theta\right)ds,\\
\hat{y}_{2}\left(x;\,\theta\right) & = & \hat{y}_{0}+\int_{0}^{x}f\left(s,\,\hat{y}_{1}\left(s;\,\theta\right);\theta\right)ds,\\
 & \vdots\\
\hat{y}_{R+1}\left(x;\,\theta\right) & = & \hat{y}_{0}+\int_{0}^{x}f\left(s,\,\hat{y}_{R}\left(s;\,\theta\right);\theta\right)ds.
\end{eqnarray*}
For any $\theta,\,\theta^{'}\in\mathbb{B}_{q}\left(1\right)$, at
the beginning, we have 
\begin{equation}
\left|\hat{y}_{1}\left(x;\theta\right)-\hat{y}_{1}\left(x;\theta^{'}\right)\right|\leq\tilde{\alpha}L_{K}\left|\theta-\theta^{'}\right|_{q}\quad\forall x\in\left[0,\,\tilde{\alpha}\right]\label{eq:31}
\end{equation}
where the inequality follows from (\ref{eq:17}). For the second iteration,
we have 
\begin{eqnarray*}
\left|\hat{y}_{2}\left(x;\,\theta\right)-\hat{y}_{2}\left(x;\,\theta^{'}\right)\right| & \leq & \left|\int_{0}^{x}f\left(s,\,\hat{y}_{1}\left(s;\,\theta\right);\theta\right)ds-\int_{0}^{x}f\left(s,\,\hat{y}_{1}\left(s;\,\theta^{'}\right);\theta\right)ds\right|+\\
 &  & \left|\int_{0}^{x}f\left(s,\,\hat{y}_{1}\left(s;\,\theta^{'}\right);\theta\right)ds-\int_{0}^{x}f\left(s,\,\hat{y}_{1}\left(s;\,\theta^{'}\right);\theta^{'}\right)ds\right|\\
 & \leq & \underset{(i)}{\underbrace{\tilde{\alpha}\left(\tilde{\alpha}L_{K}\left|\theta-\theta^{'}\right|_{q}\right)}}+\underset{(ii)}{\underbrace{\tilde{\alpha}L_{K}\left|\theta-\theta^{'}\right|_{q}}}\\
 & \leq & \tilde{\alpha}^{2}L_{K}\left|\theta-\theta^{'}\right|_{q}+\tilde{\alpha}L_{K}\left|\theta-\theta^{'}\right|_{q}\quad\forall x\in\left[0,\,\tilde{\alpha}\right]
\end{eqnarray*}
where (i) and (ii) in the second inequality follow from (\ref{eq:22})
with (\ref{eq:31}) and (\ref{eq:17}), respectively. Continuing with
this pattern until the $\left(R+1\right)$th iteration, we obtain
\begin{eqnarray}
\left|\hat{y}_{R+1}\left(x;\,\theta\right)-\hat{y}_{R+1}\left(x;\,\theta^{'}\right)\right| & \leq & \left(L_{K}\left|\theta-\theta^{'}\right|_{q}\right)\sum_{i=1}^{R+1}\tilde{\alpha}^{i}\nonumber \\
 & \leq & \frac{\tilde{\alpha}L_{K}}{1-\tilde{\alpha}}\left|\theta-\theta^{'}\right|_{q}\quad\forall x\in\left[0,\,\tilde{\alpha}\right].\label{eq:30}
\end{eqnarray}

In particular, (\ref{eq:30}) holds for $x\in\left\{ x_{1},\,x_{2},\,...,\,x_{n}\right\} $.
Consequently,
\[
\left\{ \frac{1}{n}\sum_{i=1}^{n}\left[\hat{y}_{R+1}(x_{i};\,\theta)-\hat{y}_{R+1}(x_{i};\,\theta^{'})\right]^{2}\right\} ^{\frac{1}{2}}\leq\frac{\tilde{\alpha}L_{K}}{1-\tilde{\alpha}}\left|\theta-\theta^{'}\right|_{q}.
\]
Let $\tilde{b}=\left(\frac{\tilde{\alpha}L_{K}}{1-\tilde{\alpha}}\vee1\right)$.
For a given $\delta>0$, let us consider the smallest $\frac{\delta}{2\tilde{b}}-$covering
$\left\{ \theta^{1},...,\theta^{N}\right\} $ (with respect to the
$l_{q}-$norm), and by (\ref{eq:30}), for any $\theta,\,\theta^{'}\in\mathbb{B}_{q}\left(1\right)$,
we can find some $\theta^{i}$ and $\theta^{j}$ from the covering
set $\left\{ \theta^{1},...,\theta^{N}\right\} $ such that 
\begin{eqnarray*}
 &  & \left|\hat{y}_{R+1}\left(x;\,\theta\right)-\hat{y}_{R+1}\left(x;\,\theta^{'}\right)-\left(\hat{y}_{R+1}\left(x;\,\theta^{i}\right)-\hat{y}_{R+1}\left(x;\,\theta^{j}\right)\right)\right|\\
 & \leq & \left|\hat{y}_{R+1}\left(x;\,\theta\right)-\hat{y}_{R+1}\left(x;\,\theta^{i}\right)\right|+\left|\hat{y}_{R+1}\left(x;\,\theta^{'}\right)-\hat{y}_{R+1}\left(x;\,\theta^{j}\right)\right|\\
 & \leq & \delta.
\end{eqnarray*}
Thus, $\left\{ g_{\theta^{1}},\,g_{\theta^{2}},...,\,g_{\theta^{N}}\right\} \times\left\{ g_{\theta^{1}},\,g_{\theta^{2}},...,\,g_{\theta^{N}}\right\} $
forms a $\delta-$cover of $\bar{\mathcal{F}}$ in terms of $\mathcal{F}$
defined in (\ref{eq:41}). Consequently, we have 
\begin{eqnarray*}
\frac{1}{\sqrt{n}}\int_{0}^{\tilde{r}_{n}}\sqrt{\log N_{n}(\delta,\,\Lambda(\tilde{r}_{n};\,\bar{\mathcal{F}}))}d\delta & \leq & \frac{1}{\sqrt{n}}\int_{0}^{\tilde{r}_{n}}\sqrt{\log N_{\infty}(\delta,\,\Lambda(\tilde{r}_{n};\,\bar{\mathcal{F}}))}d\delta\\
 & \leq & \sqrt{\frac{K}{n}}\int_{0}^{\tilde{r}_{n}}\sqrt{2\log\left(1+\frac{c\tilde{b}\tilde{r}_{n}}{\delta}\right)}d\delta\\
 & \precsim & \tilde{b}\tilde{r}_{n}\sqrt{\frac{K}{n}}.
\end{eqnarray*}
Setting $\tilde{b}\tilde{r}_{n}\sqrt{\frac{K}{n}}\asymp\frac{\tilde{r}_{n}^{2}}{\sigma}$
yields $\tilde{r}_{n}\asymp\sigma\tilde{b}\sqrt{\frac{K}{n}}$ and
therefore,
\begin{equation}
\left|\frac{1}{n}\sum_{i=1}^{n}\varepsilon_{i}T_{1}\left(x_{i}\right)\right|\precsim\sigma^{2}\tilde{b}^{2}\frac{K}{n}\label{eq:32-2}
\end{equation}
with probability at least $1-c_{1}\exp\left(-c_{2}\tilde{b}^{2}K\right)$. 

To analyze $\sqrt{\frac{1}{n}\sum_{i=1}^{n}\left[T_{2}\left(x_{i}\right)\right]^{2}}$
in (\ref{eq:31-1}), note that (\ref{eq:22}) implies

\begin{eqnarray*}
\left|\hat{y}_{1}\left(x_{i};\,\theta^{*}\right)-y_{1}\left(x_{i};\,\theta^{*}\right)\right| & \leq & \left|\hat{y}_{0}-y_{0}\right|+\tilde{\alpha}\left|\hat{y}_{0}-y_{0}\right|,\\
\left|\hat{y}_{2}\left(x_{i};\,\theta^{*}\right)-y_{2}\left(x_{i};\,\theta^{*}\right)\right| & \leq & \left|\hat{y}_{0}-y_{0}\right|+\tilde{\alpha}\left|\hat{y}_{0}-y_{0}\right|+\tilde{\alpha}^{2}\left|\hat{y}_{0}-y_{0}\right|,\\
 & \vdots\\
\left|\hat{y}_{R+1}\left(x_{i};\,\theta^{*}\right)-y_{R+1}\left(x_{i};\,\theta^{*}\right)\right| & \leq & \left|\hat{y}_{0}-y_{0}\right|+\tilde{\alpha}\left|\hat{y}_{0}-y_{0}\right|+\cdots+\tilde{\alpha}^{R+1}\left|\hat{y}_{0}-y_{0}\right|\leq\frac{1}{1-\tilde{\alpha}}\left|\hat{y}_{0}-y_{0}\right|
\end{eqnarray*}
for all $i=1,...,n$. As a result, we have 
\[
\sqrt{\frac{1}{n}\sum_{i=1}^{n}\left[T_{2}\left(x_{i}\right)\right]^{2}}\leq\frac{1}{1-\tilde{\alpha}}\left|\hat{y}_{0}-y_{0}\right|.
\]

For $\sqrt{\frac{1}{n}\sum_{i=1}^{n}\left[T_{3}\left(x_{i}\right)\right]^{2}}$
in (\ref{eq:31-1}), standard argument for the Picard-Lindelöf Theorem
implies that 
\[
\sup_{x\in\left[0,\,\tilde{\alpha}\right]}\left|y_{r}(x;\,\theta^{*})-y_{r+r^{'}}(x;\,\theta^{*})\right|\leq\tilde{\alpha}^{r}\frac{1-\tilde{\alpha}^{r^{'}}}{1-\tilde{\alpha}}\sup_{x\in\left[0,\,\tilde{\alpha}\right]}\left|y_{1}(x;\,\theta^{*})-y_{0}\right|
\]
for any non-negative integers $r$ and $r^{'}$; as a result, we have

\[
\sqrt{\frac{1}{n}\sum_{i=1}^{n}\left[T_{3}\left(x_{i}\right)\right]^{2}}\leq\frac{\tilde{\alpha}^{R+1}}{1-\tilde{\alpha}}\sup_{x\in\left[0,\,\tilde{\alpha}\right]}\left|y_{1}(x;\,\theta^{*})-y_{0}\right|\asymp\frac{\tilde{\alpha}^{R+1}}{1-\tilde{\alpha}}\max\left\{ C_{0},\,\tilde{\alpha}\right\} .
\]

$\square$

\end{document}